\documentclass{gtart}

\usepackage{amsmath, amssymb, latexsym, amscd, graphicx, epstopdf}
\usepackage{euscript, setspace}
\usepackage{subfigure}
\usepackage{epsfig}
\usepackage{psfrag}
\usepackage{wrapfig}
\usepackage{floatflt}

\usepackage{mathptmx}

\def\bkC{{\rm \kern.20em \vrule width.05em height1.5ex depth-.03ex \kern-.27em C}}

\def\bksC{{\rm \kern.24em \vrule width.05em height1ex depth-.05ex 
\kern-.26em C}}

\def\bkE{{\rm I\kern-.22em E}}
 
\def\bkH{{\rm I\kern-.22em H}}
 
\def\bkN{{\rm I\kern-.17em N}}

\def\bkQ{{\rm \kern.24em \vrule width.05em height1.4ex depth-.05ex 
\kern-.26em Q}}

\def\bkR{{\rm I\kern-.17em R}}
\def\RR{\bkR}
\def\bkZ{{\rm Z\kern-.32em Z}}
\def\Z{\bkZ}
\def\bksZ{{\rm Z\kern-.22em Z}}

\def\tri{\mathcal{T}}

\def\l{\lambda}

\def\m{\mu}

\def\v{{\mathfrak{v}}}
\def\w{{\mathfrak{w}}}

\def\N{PF}

\DeclareMathOperator{\im}{im}

\theoremstyle{plain}

\newtheorem{thm}{Theorem}[section]
\newtheorem*{thm*}{Theorem}
\newtheorem{lem}[thm]{Lemma}
\newtheorem*{lem*}{Lemma}
\newtheorem{cor}[thm]{Corollary}
\newtheorem*{cor*}{Corollary}
\newtheorem{pro}[thm]{Proposition}
\newtheorem*{pro*}{Proposition}
\newtheorem{definition}[thm]{Definition}

\theoremstyle{remark}
\newtheorem{rem}[thm]{Remark}
\newtheorem*{rem*}{Remark}

\numberwithin{equation}{section}

\def\normalt{t}
\def\normalq{q}

\def\simplex{\Delta}
\def\normal{\mathfrak{n}}
\def\normalc{\mathfrak{c}}

\def\pseudo{P}

\def\normala{\mathfrak{a}}
\def\normalb{\mathfrak{b}}

\def\normald{\mathfrak{d}}
\def\normalq{\mathfrak{q}}
\def\normalt{\mathfrak{t}}

\DeclareMathOperator{\Int}{int}

\addtocounter{MaxMatrixCols}{4}

\begin{document}
\title{Normal surfaces in topologically finite 3--manifolds}
\author{Stephan Tillmann}

\begin{abstract}
The concept of a normal surface in a triangulated, compact 3--manifold was generalised by Thurston to 
a spun-normal surface in a non-compact 3--manifold with ideal triangulation. This paper defines a boundary curve map which takes a spun-normal surface to an element of the direct sum of the first homology groups of the vertex linking surfaces. The boundary curve map is used to study the topology of a spun-normal surface as well as to determine the dimension of the projective solution space of the algebraic equations arising from the quadrilateral coordinates of spun-normal surfaces.
\end{abstract}

\maketitle

\section*{Introduction}

Spun-normal surfaces made their first appearance in unpublished work by Thurston. He described essential surfaces in the figure eight knot complement as spun-normal surfaces with respect to the ideal triangulation by two regular ideal hyperbolic tetrahedra. There have since been interesting applications of spun-normal surfaces; for instance in work by Weeks \cite{we}, Kang \cite{k}, Kang and Rubinstein \cite{kr}. They are also part of the repertoire of many people working in the field. There does not however appear to be a basic reference for fundamental properties of these surfaces. This paper attempts to fill this gap.

\begin{figure}[t]
\psfrag{v}{{\small $\v$}}
\psfrag{B}{{\small $B_\v$}}
    \begin{center}
      \includegraphics[width=12cm]{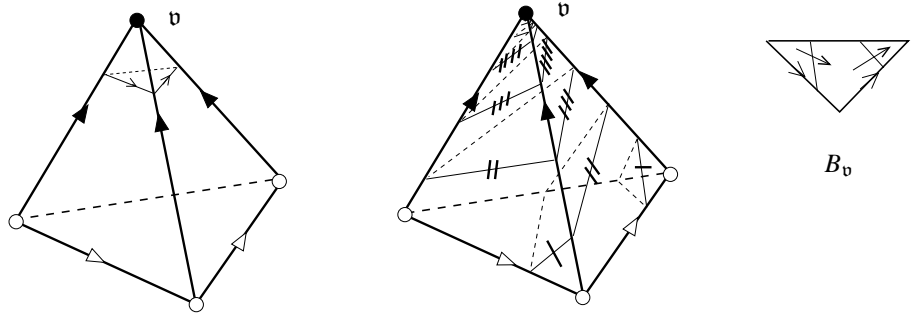}
    \end{center}
\caption{A spun-normal surface in a non-orientable pseudo-manifold with boundary. The triangulation has a single 3--simplex with two faces identified as indicated by the arrows. There are two vertices; the link of $\v$ is a M\"obius band, and the link of the other vertex is a disc. The shown spun-normal surface spins into $\v;$ it is a disc minus a point on its boundary and it meets a suitable vertex linking surface, $B_\v,$ in a 2--sided, non-separating arc.}
\label{fig:spun-surface}
\end{figure}

Let $M$ be the interior of a compact 3--manifold with non-empty boundary. Then $M$ admits an ideal triangulation and this gives the end-compactification of $M$ the structure of a triangulated, closed pseudo-manifold $\pseudo$ (see Proposition \ref{pro:ideal triangulation exists}). The complement of the 0--skeleton in $\pseudo$ is identified with $M,$ and spun-normal surfaces in $M$ can be studied using the triangulation of $\pseudo.$ The remainder of this introduction restricts to this setting; Figure \ref{fig:spun-surface} shows a spun-normal surface in a pseudo-manifold with boundary and indicates the more general situation.

A spun-normal surface meets each 3--simplex in a union of pairwise disjoint normal discs such that there are at most finitely many normal quadrilaterals but possibly infinitely many normal triangles which accumulate at the 0--skeleton of $P.$ A spun-normal surface is therefore properly embedded in $M.$ Two spun-normal surfaces are regarded as equivalent if there is an isotopy taking one to the other whilst leaving all simplices invariant (normal isotopy). One may choose a small open neighbourhood of each 0--simplex $\v$ in $\pseudo$ such that its boundary, $B_\v,$ is a closed (possibly non-orientable) normal surface. The topology of spun-normal surfaces is analysed in Section~\ref{Normal surface theory}, and the main observations are the following:
\begin{enumerate}
\item A spun-normal surface without vertex linking components is up to normal isotopy uniquely determined by its quadrilateral discs (see Lemma \ref{lem:uniqueness}).
\item If two spun-normal surfaces meet every 3--simplex in quadrilateral discs of the same type, then their \emph{geometric sum} is well-defined (see Lemma \ref{lem:geometric sum is well-defined}).
\item For each spun-normal surface $S$ and each 0--simplex $\v$ there is a well-defined element $\partial_\v(S) \in H_1(B_\v; \mathbb{Z})$ which determines $S$ in a neighbourhood of $\v$ uniquely up to normal isotopy (see Lemma \ref{lem:boundary curves} and Lemma \ref{lem:controlled spinning}).
\item A spun-normal surface $S$ is topologically finite unless there is some 0--simplex $\v$ such that $\partial_\v(S)\neq 0$ and $\chi(B_\v)<0$ (see Corollary \ref{cor:topology of spun-normal}).
\item (Kneser--Haken finiteness) If $S$ is a spun-normal surface no two components of which are normally isotopic, then the number of components of $S$ which are not vertex linking is at most $12t,$ where $t$ is the number of 3--simplices in the triangulation (see Lemma \ref{lem:Kneser--Haken finiteness} which includes a stronger bound).
\end{enumerate}

Section \ref{sec:Matching equations} analyses algebraic properties of a spun-normal surface $S$ derived from the quadrilaterals in its cell structure. Recording the number of quadrilaterals of each type gives the \emph{normal $Q$-coordinate of $S.$} The \emph{$Q$--matching equations} described by Tollefson \cite{to} for compact 3--manifolds result from the fact that a properly embedded surface meets a small regular neighbourhood of an edge in a collection of discs, and that each of these discs is uniquely determined by its intersection with the boundary of the neighbourhood. The equations are shown to give a necessary and sufficient condition on a collection of quadrilateral discs to be (up to normal isotopy) contained in the cell structure of some spun-normal surface. A spun-normal surface without vertex linking components is thus uniquely determined by its normal $Q$--coordinate (see Theorem~\ref{thm:normal unique}). 

The triangulation of $\pseudo$ is denoted by $\tri$ and the solution space of the $Q$--matching equations by $Q(\tri).$ This is a real vector subspace of $\RR^{3t}.$ The projective solution space $PQ(\tri)$ is the intersection of $Q(\tri)$ with the unit simplex; it is a convex rational polytope. Rational points are therefore dense in $PQ(\tri)$ and correspond to the projective classes not merely of spun-normal surfaces (which are embedded), but also of \emph{immersed} and \emph{branched immersed spun-normal surfaces} (see Proposition \ref{pro:normal branched immersions}).

Section~\ref{Boundary stuff} introduces the \emph{boundary curve map} 
$$\partial \co Q(\tri) \to \oplus_\v H_1(B_\v, \mathbb{R}),$$
which is linear and generalises the above map $\partial_\v.$ It is thus possible to compute $\partial_\v(S)$ directly from the normal $Q$--coordinate of $S.$ The definition of $\partial$ is motivated by the boundary map used by Kang and Rubinstein \cite{kr} in the case where each vertex linking surface is a torus or Klein bottle, and their methods are generalised to prove the following:

\begin{thm}\label{thm:kr}
Let $\pseudo$ be a closed 3--dimensional pseudo-manifold with triangulation $\tri.$ Then 
\begin{enumerate}
\item $Q(\tri)$ has dimension $v_o-e+3t=\chi(\pseudo)+2t-v_n;$ 
\item $\partial \co Q(\tri) \to \oplus_\v H_1(B_\v, \mathbb{R})$ is onto, and its restriction to integer lattice points in $Q(\tri)$ has image of finite index in $\oplus_\v H_1(B_\v, \mathbb{Z});$ and
\item $PQ(\tri)$ is non-empty and hence of dimension $\dim Q(\tri)-1.$
\end{enumerate}
Here $t, e, v$ is the number of 3--simplices, 1--simplices, 0--simplices respectively; so $\chi(\pseudo)=v-e+t.$ Moreover, $v=v_o+v_n,$ where $v_n$ is the number of 0--simplices with non-orientable linking surface.
\end{thm}

Note that $\chi(\pseudo)=0=v_n$ if $\pseudo$ is a manifold, and $\chi(\pseudo)=v$ if the link of each vertex is a torus or Klein bottle. The above dimension of $Q(\tri)$ and target of $\partial$ correct the main result of \cite{kr} in the case of pseudo-manifolds with non-orientable vertex links.

The set of projective classes of (embedded) spun-normal surfaces has a natural compactification, $\N (\tri) \subset PQ(\tri),$ points of which correspond to projective classes of certain transversely measured singular codimension-one foliations of the complement of the 0--skeleton, $M.$ The restriction $\partial\co \N (\tri) \to\oplus_\v H_1(B_\v, \mathbb{R})$ determines the (possibly singular) foliation of each $B_\v.$ If $\pseudo$ is orientable, then the intersection pairing on $\oplus_\v H_1(B_\v, \mathbb{R})$ defines a bi-linear skew-symmetric form on $Q(\tri);$ for any $N, L \in Q(\tri)$:
\begin{equation*}
\langle N , L\rangle =  \sum \iota \Big( \partial_\v(N) , \partial_\v(L) \Big).
\end{equation*}
This is used to analyse the structure of the components of $\N (\tri).$
\begin{thm}\label{thm:dim of N(M)}
Let $\pseudo$ be a closed 3--dimensional pseudo-manifold with triangulation $\tri.$ The set $\N (\tri)$ is a finite union of convex rational polytopes. Each maximal convex polytope $R$ in this union satisfies $\dim R \le t-1.$ Moreover, if $\pseudo$ is orientable, then:
\begin{equation}\label{bounds for orientable}
     \chi(\pseudo)-1 \le \dim R \le \chi(\pseudo) + \dim (R \cap \ker \partial), 
\end{equation}
and if $\pseudo$ is non-orientable, then:
\begin{equation}\label{bounds for non-orientable}
     \chi(\pseudo)-v_n-1 \le \dim R \le 2 \chi(\pseudo) -v_n+ \dim (R \cap \ker \partial), 
\end{equation}
where $\dim \emptyset = -1$ throughout.
\end{thm}

The lower bounds given in equations (\ref{bounds for orientable}) and (\ref{bounds for non-orientable}) are shown to be sharp by the examples given in Section \ref{sec:Examples}, the complement of the figure eight knot and the Gieseking manifold.

This article does not address many topics which are standard in algorithmic topology. For instance, work by Jaco and Oertel \cite{JO} implies that there is an algorithm to determine whether a topologically finite manifold is large, i.e.\thinspace whether it contains a closed, incompressible, 2--sided surface distinct from a 2--sphere or a vertex linking surface. The extension of such algorithms to spun-normal surfaces (for instance, to determine whether there are essential discs or annuli) is left for future research; it presupposes the ability to put incompressible surfaces into spun-normal form and to decide whether a spun-normal surface is incompressible --- neither of which is addressed in this paper.


\subsection*{Acknowledgements}

The author thanks Daryl Cooper for many enlightening conversations about spun-normal surfaces. He also thanks Steven Boyer, Craig Hodgson and Hyam Rubinstein for helpful comments. Thanks to Saul Schleimer for pointing out that there does not seem to be a complete treatment of the geometric sum operation in the literature. This work was supported by a Postdoctoral Fellowship by the CRM/ISM in Montr\'eal and under the Australian Research Council's Discovery funding scheme (project number DP0664276).

\begin{spacing}{0.5}
\tableofcontents
\end{spacing}

\section{Spun-normal surfaces}
\label{Normal surface theory}

This section defines spun-normal surfaces in a 3-dimensional pseudo-manifold (possibly with boundary) and analyses their basic properties. It is also shown that a topologically finite 3--manifold is homeomorphic to the complement of the 0--skeleton in a closed, triangulated pseudo-manifold.


\subsection{Pseudo-manifold}
\label{sec:3-dimensional pseudo-manifolds}

Let $\widetilde{\simplex} = \{ \widetilde{\simplex}^3_1,...,\widetilde{\simplex}^3_t\}$ be a disjoint union of 3--simplices, each of which is a regular Euclidean 3--simplex of edge length one. The subsimplices of the elements in $\widetilde{\simplex}$ are referred to as simplices in $\widetilde{\simplex}.$ Let $\Phi$ be a collection of Euclidean isometries between 2--simplices in $\widetilde{\simplex}$ with the property that for each 2--simplex $\widetilde{\simplex}^2$ there is at most one isometry having $\widetilde{\simplex}^2$ as its range or domain (but not both). There is a natural quotient map $p \co \widetilde{\simplex} \to \widetilde{\simplex} / \Phi.$ Let $\pseudo = \widetilde{\simplex} / \Phi$ and note that $\pseudo$ inherits a singular PL cone-metric. Since $p$ restricts to an embedding of the interior of every 2--simplex and every 3--simplex, the only possible non--manifold points in $\pseudo$ are at images of vertices and barycentres of 1--simplices.

\begin{definition}[Pseudo-manifold]
The quotient space $\pseudo = \widetilde{\simplex} / \Phi$ is a \emph{pseudo-manifold} if $p$ restricts to an embedding of the interior of every 1--simplex in $\widetilde{\simplex}.$
\end{definition}

A pseudo-manifold is \emph{closed} if each 2--simplex in $\widetilde{\simplex}$ is the range or domain of an element in $\Phi,$ otherwise it is a \emph{pseudo-manifold with boundary.} A pseudo-manifold is \emph{orientable} if all 3--simplices in $\widetilde{\simplex}$ can be oriented such that all elements in $\Phi$ are orientation reversing. For the remainder of this paper, it is assumed that $\pseudo$ is a pseudo-manifold and a quotient map $p \co \widetilde{\simplex} \to \pseudo$ is fixed; additional hypotheses will be added.


\subsection{Triangulation}
\label{sec:Triangulations}

The image under $p$ of an $n$-simplex in $\widetilde{\simplex}$ is called an \emph{$n$-singlex} in $\pseudo.$ The resulting combinatorial cell-decomposition of $\pseudo$ is termed a \emph{triangulation} of $\pseudo$ and denoted by $\tri=(\widetilde{\simplex}, \Phi, p).$ If $p$ restricts to an embedding of each simplex, then $\tri$ is a triangulation in the traditional sense; it will then be referred to as a \emph{simplicial triangulation}. The image in $\pseudo$ of the set of $i$--simplices in $\widetilde{\simplex}$ is denoted by $\pseudo^{(i)}$ and termed the \emph{$i$--skeleton of $\pseudo.$}

The number of 3--singlices in $\pseudo$ is $t,$ and we denote by $f$ the number of 2--singlices, $e$ the number of 1--singlices and $v$ the number of 0--singlices.


\subsection{Ideal triangulation}
\label{sec:Ideal triangulations}

A manifold is \emph{topologically finite} if it is homeomorphic to the interior of a compact 3--manifold. If $\pseudo$ is closed, then removing the 0--skeleton yields a topologically finite 3--manifold $M=\pseudo \setminus \pseudo^{(0)}$ and $\pseudo$ is called the \emph{end-compactification} of $M.$ This motivates the following terminology. An \emph{ideal $i$--simplex,} $i \in \{1,2,3\},$ is an $i$--simplex with its vertices removed; the vertices of the $i$--simplex are referred to as the \emph{ideal vertices} of the ideal $i$--simplex. Similarly for singlices. The restriction $(\widetilde{\simplex}\setminus \widetilde{\simplex}^{(0)}, \Phi|_{\widetilde{\simplex}^{(2)}\setminus \widetilde{\simplex}^{(0)}}, p|_{\widetilde{\simplex}\setminus \widetilde{\simplex}^{(0)}})$ of the triangulation of $\pseudo$ to $\pseudo \setminus \pseudo^{(0)}$ is an \emph{ideal triangulation} of $\pseudo \setminus \pseudo^{(0)}.$ A 0--simplex of $\pseudo$ is referred to as an \emph{ideal vertex} of the ideal triangulation or of $M.$

If $M$ is the interior of a compact 3--manifold $\overline{M}$ with non--empty boundary and $M$ is homeomorphic to $\pseudo \setminus \pseudo^{(0)}$ for some pseudo-manifold $P,$ then $M$ is said to admit an ideal triangulation. The following result is implicit in Matveev \cite{M} and shows that every topologically finite 3-manifold arises in this way from a closed pseudo-manifold:

\begin{pro}[Topologically finite manifold has ideal triangulation]\label{pro:ideal triangulation exists}
If $M$ is the interior of a compact 3--manifold $\overline{M}$ with non--empty boundary, then $M$ admits an ideal triangulation. The ideal vertices of the ideal triangulation are in one--to--one correspondence with the boundary components of $\overline{M}.$
\end{pro}

\begin{proof}
This follows from the following results in \cite{M}. Theorem 1.1.13 due to Casler asserts that $\overline{M}$ possesses a special spine $\Sigma;$ it has the property that $\overline{M}$ is homeomorphic to a regular neighbourhood of $\Sigma$ in $M.$ Theorem 1.1.26 implies that $\Sigma$ is dual to an ideal triangulation of $M$ with the property that the ideal vertices are in one--to--one correspondence with the boundary components of $\overline{M}.$
\end{proof}


\subsection{Normal discs, arcs and corners}

A \emph{normal corner} is an interior point of a 1--simplex. A \emph{normal arc} is a properly embedded straight line segment on a 2--simplex with boundary consisting of normal corners. A \emph{normal disc} is a properly embedded disc in a 3--simplex whose boundary consists of normal arcs no two of which are contained on the same face of the 3--simplex; moreover, the normal disc is the cone over its boundary with cone point the barycentre of its normal corners. It follows that the boundary of a normal disc consists of either three or four normal arcs, and it is accordingly called a \emph{normal triangle} or a \emph{normal quadrilateral}. A normal disc is uniquely determined by its intersection with the 1--skeleton.

\begin{figure}[t]
\begin{center}
  \includegraphics[width=10cm]{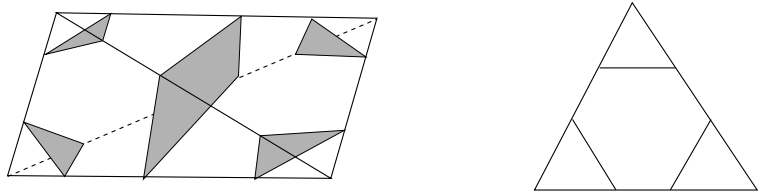}
\end{center}
    \caption{Normal discs in a 3--simplex and normal arcs in a 2--simplex}
    \label{fig:Q_normal_discs}
\end{figure}

A \emph{normal isotopy} is an isotopy of a simplex which preserves all its subsimplices. There are exactly three normal isotopy classes of normal arcs in a 2--simplex. A normal arc is called $\v$--type if it separates the vertex $\v$ from the other two vertices of the 2--simplex. There are exactly seven normal isotopy classes of normal discs in a 3--simplex (four of normal triangles and three of normal quadrilaterals). Given a normal triangle $\normalt$ in a 3--simplex, there is a unique complementary region containing precisely one 0--simplex $\v;$ $\normalt$ is said to be \emph{dual} to $\v.$

A piecewise linear arc in a 2--simplex is termed an \emph{elementary arc} if it is normally isotopic to a normal arc. A piecewise linear disc in a 3--simplex is termed an \emph{elementary disc} if it is normally isotopic to a normal disc. The following lemma provides a useful normal isotopy for any collection of elementary cells.

\begin{lem}[Straightening]
Given an arbitrary collection of pairwise disjoint elementary arcs in a 2--simplex, there is a normal isotopy keeping the 1--skeleton fixed and taking each elementary arc to a normal arc. Given a arbitrary collection of pairwise disjoint elementary discs in a 3--simplex, there is a normal isotopy keeping the 1--skeleton fixed and taking each elementary disc to a normal disc.
\end{lem}

\begin{proof}
Let $\simplex^2$ be a 2--simplex containing an arbitrary collection of pairwise disjoint elementary arcs. There is a homeomorphism $\simplex^2 \to \simplex^2$ which fixes the boundary and takes each elementary arc to a normal arc. The first part of the lemma now follows from the fact that a homeomorphism of a closed ball which fixes the boundary is isotopic to the identity. Given an arbitrary collection of pairwise disjoint elementary discs in a 3--simplex $\simplex^3,$ there is a homeomorphism $h\co \simplex^3 \to \simplex^3$ which fixes the 1--skeleton and takes each elementary disc to a normal disc. The first part of the lemma shows that the restriction of $h$ to $\partial \simplex^3$ is isotopic to the identity. One can therefore isotope $h$ in a neighbourhood of $\partial \simplex^3$ such that the boundary of each elementary disc is straight and the result follows as above.
\end{proof}

The image under $p$ of a normal disc (respectively arc, corner) in $\widetilde{\simplex}$ is called a \emph{normal disc (respectively arc, corner) in $\pseudo$.} Note that a normal disc in $\pseudo$ may have fewer arcs or corners in its boundary than its pre-image in $\widetilde{\simplex}.$ Two normal discs in $\pseudo$ are said to be of the same \emph{type} if they are the images of normally isotopic normal discs in $\widetilde{\simplex}.$ A \emph{normal isotopy of $\pseudo$} is an isotopy of $\pseudo$ which leaves all singlices invariant. 

The above definitions apply verbatim to define normal cells in ideal simplices and ideal singlices, and normal isotopies in $\pseudo \setminus \pseudo^{(0)}.$


\subsection{Normal surfaces and normal curves}

This subsection collects some well--known facts about normal surfaces; please consult Jaco and Rubinstein \cite{JR89} and Thompson \cite{AT94} for details and proofs.

\begin{definition}[Normal surface]
A subset $S$ of $\pseudo$ is a \emph{normal surface} in $\pseudo$ if it meets every 3--singlex in $\pseudo$ in a (possibly empty) finite union of pairwise disjoint normal discs.
\end{definition}

Note that a normal surface is a properly embedded (nor necessarily connected) compact surface in $(\pseudo\setminus \pseudo^{(0)}) \subset \pseudo.$ 

\begin{lem}\label{lem: normal arcs determined by 1--skeleton}
Let $A$ be a finite union of pairwise disjoint normal arcs in $\simplex^2.$ Then $A$ is uniquely determined by its intersection with the 1--skeleton up to normal isotopy.
\end{lem}

A curve on a compact, triangulated surface is \emph{normal} if it meets each 2--singlex in a finite collection of pairwise disjoint normal arcs. (A curve will always be connected, whilst a surface may not be connected.)

\begin{lem}
Every non-empty normal curve in $\partial \simplex^3$ is simple and closed. Moreover, it bounds some properly embedded disc in $\simplex^3.$
\end{lem}

Consequences of Lemma \ref{lem: normal arcs determined by 1--skeleton} are the following:

\begin{cor}\label{cor:curves on simplex determined by 1--skeleton}
A finite union of pairwise disjoint normal curves on $\partial \simplex^3$ is uniquely determined by its intersection with the 1--skeleton up to normal isotopy.
\end{cor}

\begin{cor}
Let $S$ be a normal surface in $\pseudo.$ Then $S$ is uniquely determined by its intersection with the 1--skeleton of $\pseudo$ up to normal isotopy.
\end{cor}

The \emph{length} of a normal curve is the number of intersections with the 1--skeleton (equivalently: the number of normal arcs it is composed of).

\begin{lem}
Every normal curve $\gamma$ bounds a disc in $\partial \simplex^3$ containing one or two vertices. In the former case, $\gamma$ has length three, in the latter, its length is a multiple of four.
\end{lem}

\begin{lem}
If a non-empty normal curve $\gamma$ does not meet some edge of $\simplex^3,$ then $\gamma$ has length three or four.
\end{lem}


\subsection{Spun-normal surfaces}

The definition of a normal surface is now extended; the following notions are defined with respect to the PL cone-metric on $\pseudo.$

\begin{definition}[Normal subset]\label{defn:normal subset}
Let $\pseudo$ be a pseudo-manifold (possibly with boundary). A subset $S$ of $\pseudo$ is \emph{normal} if
\begin{enumerate}
\item $S$ intersects each 3--singlex in $\pseudo$ in a (possibly empty) countable union of pairwise disjoint normal discs;
\item the set of accumulation points of $S \cap \pseudo^{(1)}$ is contained in $\pseudo ^{(0)};$
\item if $\{ x_i\} \subset S$ has accumulation point $x\in \pseudo,$ then $x \in S$ or $x \in \pseudo ^{(0)}.$
\end{enumerate}
\end{definition}

Note that any normal subset of $\pseudo$ is contained in $\pseudo\setminus \pseudo^{(0)}.$ The following two facts follow from the observation that whenever a normal subset $S$ of $\pseudo$ meets a 3--singlex $\simplex^3$ in a point $x \in \partial \simplex^3,$ then it meets $\simplex^3$ in a normal disc containing $x.$

\begin{lem}
Let $S$ be a normal subset of $\pseudo.$ Then $p^{-1}(S)$ is a normal subset of $\widetilde{\simplex}.$
\end{lem}

\begin{lem}[Normal subset is surface]
Let $S$ be a normal subset of $\pseudo.$ Then every point on $S$ has a small closed neighbourhood in $S$ which is homeomorphic to a disc.
\end{lem}

For each 0--singlex $\v\in \pseudo^0$ choose a small, open neighbourhood $N_\v$ with the property that the neighbourhoods are pairwise disjoint and $\overline{\partial N_\v \setminus \partial \pseudo}= B_\v$ is a normal surface in $\pseudo.$ Note that all normal discs in $B_\v$ are normal triangles no two of which are normally isotopic. A normal subset in $\pseudo$ is called a \emph{vertex linking surface} if it is normally isotopic to $B_\v$ for some 0--singlex $\v.$ Let $\pseudo^c = \pseudo \setminus \cup N_\v.$

\begin{rem}
Condition (2) in Definition \ref{defn:normal subset} is not redundant. For instance, it rules out that countably many pairwise disjoint copies of $B_\v$ accumulate on $B_\v.$
\end{rem}

Each surface $B_\v$ consists of normal triangles and hence inherits an induced triangulation $\tri_\v.$ A path in $B_\v$ is called \emph{normal} if it intersects each triangle in $\tri_\v$ in a (possibly empty) disjoint union of normal arcs. 

\begin{lem}[Vertex linking surface]\label{lem:vertex linking surface}
If $S$ is a normal subset of $\pseudo$ consisting only of normal triangles, then each connected component of $S$ is a vertex linking surface.
\end{lem}

\begin{proof}
Without loss of generality, it may be assumed that $S$ is connected. Given a normal triangle in $S$ dual to $\v \in \pseudo^{(0)},$ the normal triangles it meets along its boundary arcs are all dual to $\v.$ Whence all normal triangles in $S$ are dual to a fixed 0--singlex $\v.$ Since $S \cap \pseudo^c \cap \pseudo^{(1)}$ consists of finitely many points, one may (up to normal isotopy) assume that $S$ is contained in $N_\v.$ Since $S$ meets every 3--singlex in a union of pairwise disjoint normal discs, it follows that it meets every 3--singlex incident with $\v$ in at least one normal triangle contained in $N_\v.$ This forces $S$ to be normally isotopic to $B_\v.$
\end{proof}

\begin{lem}[Finiteness]\label{lem:finiteness}
A normal subset $S$ of $\pseudo$ contains at most finitely many normal quadrilaterals.
\end{lem}

\begin{proof}
If $S$ contains infinitely many normal quadrilaterals, then it contains infinitely many quadrilaterals of the same type. Using condition (2), there are finitely many possibilities for the limit points of the normal corners of these quadrilaterals. In each case, one can construct a sequence of points on $S$ which accumulates on a point not on $S$ but in the interior of a 1--simplex in $\pseudo,$ contradicting the third condition.
\end{proof}

\begin{lem}
Every connected component of a normal subset of $\pseudo$ is a properly embedded surface in $\pseudo\setminus \pseudo^{(0)}.$
\end{lem}

\begin{proof}
It suffices to assume that $S$ is connected and contains a normal quadrilateral. Then $S$ contains finitely many normal quadrilaterals. If $S$ contains only finitely many normal triangles, then $S$ is a normal surface and therefore a closed, properly embedded surface in $\pseudo\setminus \pseudo^{(0)}.$ Hence assume that $S$ contains infinitely many normal triangles. We may assume that $\pseudo^c$ contains all normal quadrilaterals contained in $S$ and that $\partial \pseudo^c \cap S$ is transverse. Conditions (2) and (3) in Definition~\ref{defn:normal subset} imply that $\pseudo^c \cap \pseudo^{(1)} \cap S$ contains finitely many points. Since normal triangles are flat, it follows that $S \cap \partial \pseudo^c$ consists of finitely many pairwise disjoint simple closed curves and arcs. Hence $S \cap \pseudo^c$ is a properly embedded surface in $\pseudo^c.$ Since $N_\v$ can be made arbitrarily small, the conclusion follows.
\end{proof}

It follows from Lemma \ref{lem:vertex linking surface} and Lemma \ref{lem:finiteness} that a normal subset may have infinitely many connected components at most finitely many of which are not vertex linking. 

\begin{definition}[Spun-normal surface]
A normal subset of $\pseudo$ is termed a \emph{spun-normal surface} in $\pseudo$ if it has finitely many connected components. If a spun-normal surface contains infinitely many normal triangles dual to $\v \in \pseudo^{(0)},$ then it is said to \emph{spin into $\v.$} In particular, a normal surface is a spun-normal surface which spins into no vertex.
\end{definition}

\begin{lem}[Uniqueness]\label{lem:uniqueness}
Let $S_0$ and $S_1$ be normal subsets in $\pseudo$ with the property that for any type of normal disc, $S_0$ and $S_1$ contain the same (not necessarily finite) number of normal discs of this type. Then $S_0$ and $S_1$ are normally isotopic.
\end{lem}

\begin{proof}
There is a normal isotopy in $\widetilde{\simplex}$ which takes $p^{-1}(S_0)$ to $p^{-1}(S_1)$ and, in particular, takes the intersection of $p^{-1}(S_0)$ with the 1--skeleton to the intersection of $p^{-1}(S_1)$ with the 1--skeleton. From the latter, a normal isotopy taking $p^{-1}(S_0)$ to $p^{-1}(S_1)$ can be constructed which descends to $\pseudo.$
\end{proof}


\subsection{Regular exchange and geometric sum: Normal arcs}

Given two normal subsets $S$ and $F$ of $\pseudo,$ one may perturb $S$ by an arbitrarily small normal isotopy so that $S$ and $F$ have no common point of intersection with the 1--skeleton, and one can attempt to produce a new normal subset from $S \cup F$ by performing a canonical cut and paste operation along $S \cap F.$ Such a \emph{geometric sum of $S$ and $F$} was first used by Haken. Since a normal disc is uniquely determined by its boundary, it is convenient to start with the geometric sum of normal arcs.

\begin{definition}[Normal subset of $\simplex^2$]\label{defn:normal subset 2}
A subset $A$ of a 2--simplex $\simplex$ is \emph{normal} if
\begin{enumerate}
\item $A$ is a (possibly empty) countable union of pairwise disjoint normal arcs;
\item the set of accumulation points of $A \cap \simplex^{(1)}$ is contained in $\simplex^{(0)};$
\item if $\{ x_i\} \subset A$ has accumulation point $x\in \simplex,$ then $x \in A$ or $x \in \simplex^{(0)}.$
\end{enumerate}
\end{definition}

\begin{definition}[General position]
Two normal subsets of a 2--simplex $\simplex^2$ are said to be in \emph{general position} if their intersection is contained in the interior of $\simplex^2.$
\end{definition}

\begin{figure}[t]
\begin{center}
  \subfigure[Geometric sum of normal arcs of different types]{
      \includegraphics[height=2.2cm]{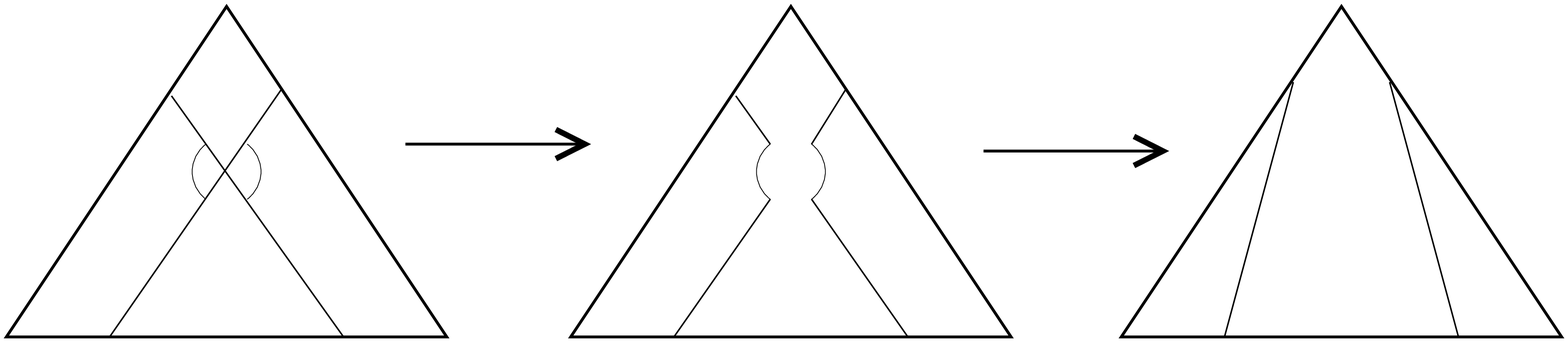}
    } 
    \\
\subfigure[Geometric sum of normally isotopic arcs]{
      \includegraphics[height=2.2cm]{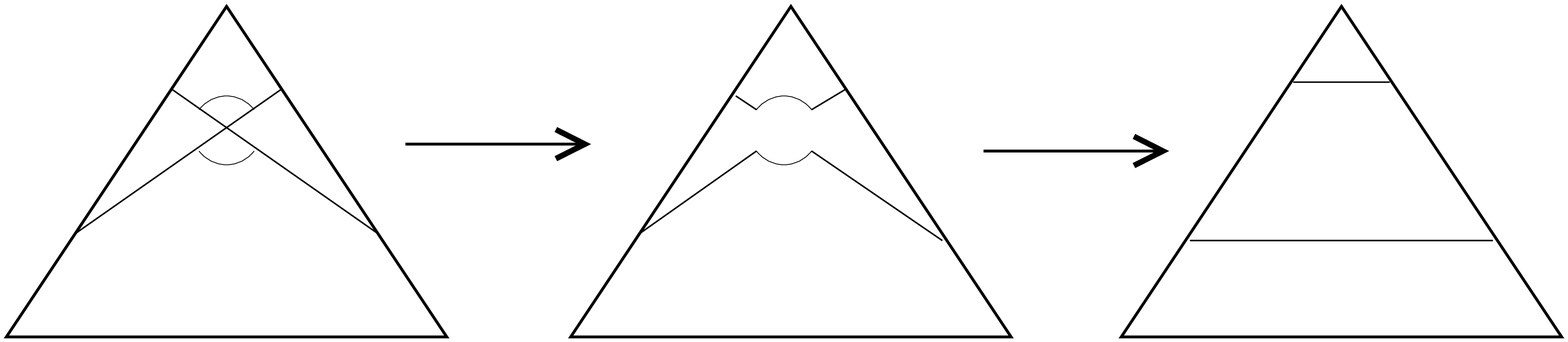}
    } 
\end{center}
    \caption{Regular exchange of normal arcs}
     \label{fig:regular exchange}
\end{figure}

Let $\normala$ be a normal arc on a 2--simplex $\simplex^2$ dual to vertex $\v.$ Denote by $C_\normala$ the complementary region of $\simplex^2 \setminus \normala$ containing $\v.$ Given normal arcs $\normala, \normalb$ on $\simplex^2$ in general position, denote by $D$ a small (piecewise linear) disc in $\Int \simplex^2$ with centre $\normala \cap \normalb.$ Then define $\normala \uplus \normalb$ to be the union of two normal arcs obtained as follows. Remove from $\normala \cup \normalb$ the intersection with $D$ and adjoin to it two arcs on $\partial D$ (termed \emph{circular arcs}) as follows. If $\normala$ and $\normalb$ are of different types, then add the intersection of $\partial D$ with $(C_\normala \cup C_\normalb) \setminus (C_\normala \cap C_\normalb).$ Otherwise take the intersection with the complement thereof. This procedure is termed a \emph{regular exchange at the intersection point,} and the position of the circular arcs is referred to as the \emph{switch condition at the intersection point}. The result is a disjoint union of two elementary arcs which can therefore be straightened to a disjoint union of two normal arcs which is denoted by $\normala \uplus \normalb$ and termed the \emph{geometric sum of $\normala$ and $\normalb$.} This is illustrated in  Figure \ref{fig:regular exchange}.

Let $A$ and $B$ be two normal subsets of a 2--simplex $\simplex^2$ which are in general position. For each intersection point $p\in A\cap B,$ one can choose a sufficiently small disc $D$ with centre $p$ so that if $p = \normala \cap \normalb$ for $\normala \subset A$ and $\normalb \subset B,$ then $D \cap (A\cup B) = D \cap (\normala \cup \normalb),$ and any two such discs are disjoint. One can then perform a regular exchange at each intersection point and the result, $\mathcal{C},$ is (a-priori) a collection of embedded arcs (possibly open or half-open) and circles in $\simplex^2.$ This is termed \emph{the regular exchange at $A\cap B$} since $\mathcal{C}$ is uniquely determined up to normal isotopy fixing the 1--skeleton.

The definition of regular exchange extends to any two \emph{elementary} arcs meeting transversely in a single point in the interior of $\simplex^2.$ Since the resulting normal arcs are uniquely determined by the relative position of the intersection points of the elementary arcs with $\partial \simplex^2,$ it follows that performing a regular exchange between two elementary arcs and then straightening gives the same result as first straightening the elementary arcs, then taking the geometric sum. The following argument is similar to the one given in lecture notes by Cameron Gordon \cite{G}.

\begin{lem}[Geometric sum of normal arcs]\label{lem:regular exchange for arcs}
Let $A$ and $B$ be two normal subsets of a 2--simplex $\simplex^2$ which are in general position. Then the regular exchange at $A\cap B$ yields a family of pairwise disjoint elementary arcs. It can be straightened to a normal subset of $\simplex^2;$ this is denoted by $A \uplus B$ and called the \emph{geometric sum of $A$ and $B$.} Moreover, if $| \cdot |_\v \in \mathbb{N} \cup \{\infty\}$ denotes the number of normal arcs dual to vertex $\v$ of $\simplex^2,$ then $|A \uplus B|_\v = |A|_\v + |B|_\v.$
\end{lem}

\begin{proof}
Choose discs for the regular exchange at $A \cap B,$ and hence two circular arcs at each intersection point determined by the switch condition. These discs and arcs will be fixed throughout the proof; the set of circular arcs is denoted by $C.$ Since $A$ and $B$ are normal subsets, the normal subset $A_0$ consisting of all normal arcs $\normala \subset A$ with the property that $\normala$ meets an arc in $B$ not normally isotopic to $\normala$ is finite. Define $B_0\subset B$ analogously. Assume that some arc in $A_0$ meets an arc in $B_0$ which is not normally isotopic to it. The properties that will be used are the following: (1) each connected component of $A_0$ (resp.\thinspace $B_0$) is an elementary arc which is made up of \emph{straight subarcs} of arcs in $A \cup B$ and circular arcs in $C,$ (2) the arcs in $A_0$ (resp.\thinspace $B_0$) are pairwise disjoint, (3) whenever an elementary arc in $A_0$ meets an elementary arc in $B_0,$ then they meet in a single point, and the corresponding regular exchange can be done by deleting subarcs from $A \cup B$ and inserting arcs from $C.$

Denote the vertices of $\simplex^2$ by $x,y,z$ and the opposite edges by $s_x,s_y,s_z$ respectively. Label the normal corners of elementary arc $\normala$ dual to $\v$ by $\v.$ For each side of $\simplex^2,$ this gives a finite word in the labels. Consider the side $s_z.$ If some arc dual to $x$ meets an arc dual to $y$, then on $s_z,$ read from $x$ to $y,$ the word has the form $w_0(x,y)\cdot y \cdot w_1(x,y) \cdot x \cdot w_2(x,y).$  Consequently, if the word is of the form $x^n y^m,$ then there are no intersections between arcs dual to $x$ and arcs dual to $y.$ Assume that there is some point of intersection. Then the word contains the subword $yx.$ Since arcs in $A_0$ are pairwise disjoint, the two labels do not correspond to endpoints of two normal arcs in $A_0.$ Similar for $B_0.$ Without loss of generality, assume that $y$ is the label of an endpoint of $\normala \subset A_0$ and $x$ is the label of an endpoint of $\normalb \subset B_0.$ Perform the regular exchange at $\normala\cap\normalb$ using circular arcs from $C.$ This yields two elementary arcs, one containing a portion of $\normalb$ between $s_y$ and the intersection point as well as the corner $\normala \cap s_y.$ This elementary arc, denoted by $\normalb_0$ has the property that $\normalb_0 \cap (B_0 \setminus \normalb)=\emptyset$ as well as $\normalb_0 \cap A_0 = (\normalb \cap A_0) \setminus (\normala\cap\normalb).$ Moreover, the switch condition at any point of $\normalb_0 \cap A_0$ coincides with the switch condition at that point with respect to $\normalb,$ since $\normalb_0$ and $\normalb$ are normally isotopic. In particular, it can be realised using circular arcs from $C.$ Moreover, any elementary arc met by $\normalb_0$ has at most one point of intersection with it since none of the subarcs in which $\normalb_0$ and $\normalb$ differ meet any other arc. The same discussion applies to the other elementary arc, $\normala_0,$ arising from the regular exchange. Let $A_1$ be the collection of elementary arcs obtained from $A_0$ by deleting $\normala$ and adding $\normala_0,$ and define $B_1$ likewise. Then $|A_0|_\v + |B_0|_\v = |A_1|_\v+|B_1|_\v$ for each $\v \in \{x,y,z\},$ $|A_1|=|A_0|,$ $|B_1|=|B_0|,$ $|A_1\cap B_1|=|A_0\cap B_0|-1$ and $A_1$ and $B_1$ satisfy the above properties (1)--(3). It follows that this process can be iterated to yield two normal subsets $A_k$ and $B_k$ with the property that $|A_k|=|A_0|,$ $|B_k|=|B_0|,$ and whenever arcs in $A_k$ and $B_k$ meet, then they are normally isotopic. Moreover, the number of arcs dual to $x$ (resp.\thinspace $y,z$) in $A_k\cup B_k$ equals the number of normal arcs in $A_0$ dual to $x$ (resp.\thinspace $y,z$) plus the number of normal arcs in $B_0$ dual to $x$ (resp.\thinspace $y,z$).

Letting $A' = (A \setminus A_0) \cup A_k$ and $B' = (B \setminus B_0) \cup B_k$ gives sets with the property that each point in $A'\cap B'$ is the intersection point of normally isotopic arcs; $A'$ and $B'$ can be viewed as obtained by performing all regular exchanges at all points in $A\cap B$ corresponding to intersection points of arcs of different normal types. Let $A_x$ and $B_x$ be the sets of all elementary arcs dual to $x$ in $A'$ and $B'$ respectively. The elementary arcs in $A_x$ can be labelled $\normala_0, \normala_1, ...$ with the property that $C_{\normala_0} \supset C_{\normala_1} \supset \dots.$ Similarly for $B_x.$ Since $A$ and $B$ are in general position, the intersection points of $A_x\cup B_x$ with the side $s_y$ (resp.\thinspace $s_z$) can be labelled $y_1, y_2, ...$ (reps.\thinspace $z_1, z_2, ....$) such that $(y_i, x) \subset (y_{i-1}, x)$ (resp.\thinspace $(z_i, x) \subset (z_{i-1}, x)$) for each $i > 1.$ It now follows inductively that performing all regular exchanges at $A_x \cap B_x$ yields a union of pairwise disjoint elementary arcs $\{\normal_i\}$ with $\normal_i$ having endpoints $y_i$ and $z_i.$ The resulting number of normal arcs dual to $x$ is thus uniquely determined by the intersection points with the 1--skeleton. Whence the claimed relationship.

If follows that performing all regular exchanges at $A \cap B$ yields a union of pairwise disjoint elementary; it can be pulled straight to a union of pairwise disjoint normal arcs, $A \uplus B.$ The other parts of Definition \ref{defn:normal subset 2} are verified as follows. We have $(A \uplus B) \cap \simplex^{2(1)}=(A \cap \simplex^{2(1)})\cup (B \cap \simplex^{2(1)}),$ whence the set of limit points of $(A \uplus B) \cap \simplex^{2(1)}$ is contained in $\simplex^{2(0)}.$ It follows from Definition \ref{defn:normal subset 2} that for each vertex $\v$ there is a neighbourhood $N_\v$ with the property that $N_\v$ only meets normal arcs dual to $\v$ in $A\cup B.$ It follows from the construction that $N_\v$ can be chosen such that it meets only normal arcs dual to $\v$ in $A \uplus B.$ Thus, $A \uplus B$ is a normal subset of $\simplex^2.$
\end{proof}

\begin{figure}[t]
\begin{center}
    \subfigure[Geometric sum of immersed may not be normal]{
      \includegraphics[height=2.2cm]{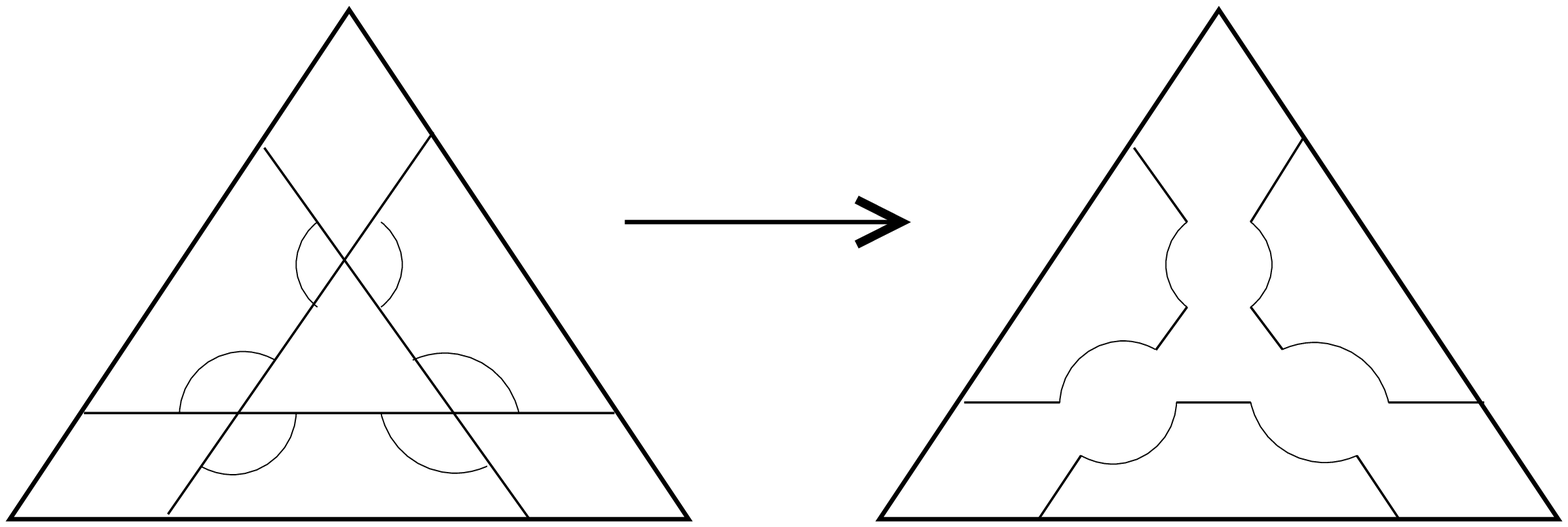}
    } 
    \\
    \subfigure[Geometric sum is associative]{
      \includegraphics[height=2.2cm]{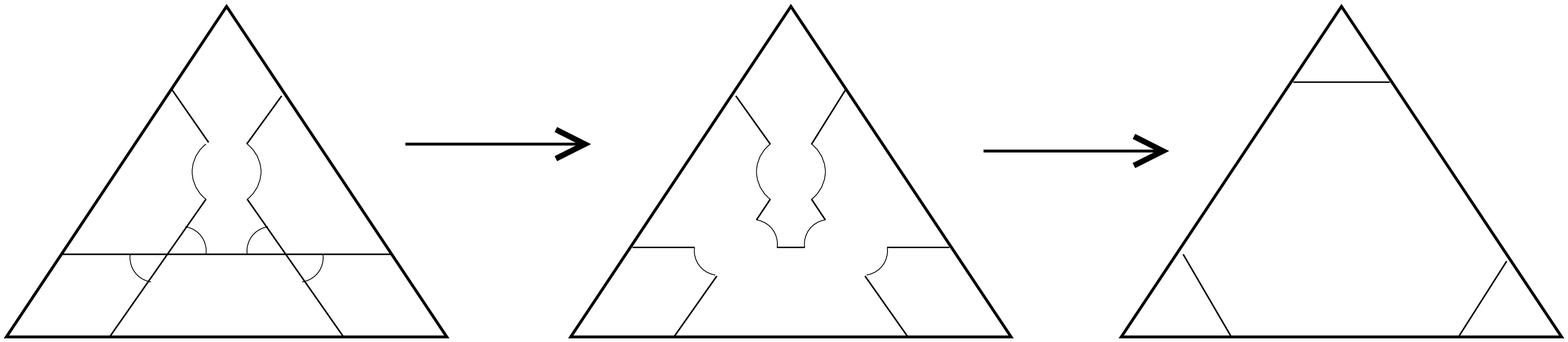}
    } 
\end{center}
    \caption{Regular exchange is associative}
     \label{fig:regular exchange is associative}
\end{figure}

The method of pre-assigning switch conditions at intersection points fails for certain immersed sets of normal arcs. Assume $A=\{\normala_i\}$ is a set consisting of three normal arcs on a 2--simplex $\simplex^2$ which are pairwise in general position and no two of which are normally isotopic. Then for any distinct $\normala_i, \normala_j \in A,$ the regular exchange at $\normala_i \cap \normala_j$ is defined. However, if all these regular exchanges are performed, then the result may not be elementary; Figure \ref{fig:regular exchange is associative}(a) shows a constellation where this fails. However, we have $(\normala_0 \uplus \normala_1) \uplus \normala_2 = \normala_0 \uplus (\normala_1 \uplus \normala_2) = (\normala_0 \uplus \normala_2) \uplus \normala_1$ even though intersection points are possibly resolved differently. See Figure \ref{fig:regular exchange is associative}(b) and notice that the last picture is invariant under rotation by $\frac{2\pi}{3}$ whilst the first and second are not. 

\begin{lem}
Let $A,B,C$ be three normal subsets of a 2--simplex which are pairwise in general position. Then
$A \uplus (B \uplus C) = (A \uplus B) \uplus C.$
\end{lem}

\begin{proof}
Let $A'$ be the set of all normal arcs in $A$ which meet normal arcs of different types in $B$ or $C;$ define $B'$ and $C'$ similarly. Then the previous lemma applies to show that $A' \uplus (B' \uplus C')$ and $(A' \uplus B') \uplus C'$ are normal subsets with the property that for every vertex $\v,$ $|A' \uplus (B' \uplus C')|_\v = |A'|_\v + |B' \uplus C'|_\v = |A'|_\v + |B' |_\v+| C'|_\v = |A' \uplus B'|_\v + | C'|_\v = |(A' \uplus B') \uplus C'|_\v.$ Whence $A' \uplus (B' \uplus C')$ and $(A' \uplus B') \uplus C'$ are identical, since they have identical intersection with the 1--skeleton.

It remains analyse arcs which are pairwise in general position and dual to a common vertex. Here the switch conditions at intersection points are uniquely determined regardless of the order in which regular exchanges are performed. The lemma follows.
\end{proof}

Given two normal subsets $A$ and $B$ on a triangulated, compact surface, let $A\uplus B = \cup ((A\cap \simplex^2) \uplus (B\cap \simplex^2)),$ where the union is taken over all 2--singlices in the triangulation. Then $A\uplus B$ is a well-defined normal subset.

\begin{lem}\label{lem:algebra for finite union of normal curves}
Let $A$ and $B$ be two finite unions of pairwise disjoint normal curves on $\partial \simplex^3$ with the property that for each face $\simplex^2,$ the sets $A\cap \simplex^2$ and $B\cap \simplex^2$ are in general position. Then $C=A \uplus B$ is a finite collection of normal curves. Moreover, if $A$ and $B$ only contain curves of length three or four, and each curve of length four does not meet the same edge $\simplex^1$ of $\simplex^3,$ then $C$ contains only curves of length three or four and each curve of length four does not meet $\simplex^1.$
\end{lem}

\begin{proof}
The statement follows from Corollary \ref{cor:curves on simplex determined by 1--skeleton} and Lemma \ref{lem:regular exchange for arcs}. By assumption, one may normally isotope $A$ to a normal subset $A'$ which is disjoint from $B.$ Then the intersection of $A'\cup B$ with the 1--skeleton is normally isotopic to the intersection of $A\cup B$ with the 1--skeleton. It follows that $A'\cup B$ is normally isotopic to $A \uplus B.$
\end{proof}


\subsection{Regular exchange and geometric sum: Normal discs}

\begin{definition}[General position]
Two normal subsets $S$ and $F$ of a 3--simplex $\simplex^3$ are said to be in \emph{general position} if $S \cap F \cap \simplex^{3(1)}= \emptyset.$ In particular, for every face $\simplex^2$ of $\simplex^3,$ the normal subsets $S \cap \simplex^2$ and $F \cap \simplex^2$ are in general position
\end{definition}

Let $\normald_0$ and $\normald_1$ be two normal discs in a 3--simplex $\simplex^3$ which are in general position. A geometric sum of $\normald_0$ and $\normald_1$ should restrict to the geometric sum of $\partial \normald_0$ and $\partial \normald_1$ and yield two normal discs. Whence each normal curve in $\partial \normald_0\uplus \partial \normald_1$ should have length three or four. It follows that $\normald_0$ and $\normald_1$ cannot be quadrilateral discs of different types since otherwise $\partial \normald_0\uplus \partial \normald_1$ is a single curve of length eight.

Hence assume that if both $\normald_0$ and $\normald_1$ are normal quadrilaterals, then they are normally isotopic. Also assume that $\alpha=\normald_0 \cap \normald_1\neq \emptyset.$ Recall the definition of a normal disc as the cone to the barycentre of its vertices. Any normal triangle is a flat Euclidean triangle, and a normal quadrilateral is made up of (at most) four Euclidean triangles. The incidence between two normal discs does not change under a normal isotopy of their union, so to determine the possibilities for $\alpha,$ it may be assumed that if $\normald_0$ is a normal quadrilateral, then it is flat. It follows that if not both $\normald_0$ and $\normald_1$ are normal quadrilaterals, then $\alpha$ is a properly embedded arc in each disc and necessarily has endpoints on distinct faces of $\simplex^3.$ If both discs are normal quadrilaterals, then $\alpha$ cannot contain a connected component not meeting $\partial \simplex^3$ for otherwise $\normald_0$ would separate the vertices of $\normald_1$ from the barycentre of $\normald_1.$ Moreover, $\alpha$ meets $\partial \simplex^3$ in either two or four points. It follows that $\alpha$ is either a single arc, the disjoint union of two arcs, or a cross; each properly embedded in $\simplex^3.$ In particular, the above definition allows a piecewise linear analogue of saddle tangencies. 

Let $N_\alpha$ be a regular neighbourhood of $\alpha$ in $\simplex^3.$ A \emph{regular exchange of $\normald_0$ and $\normald_1$ at $\alpha$} consists of deleting the portion of $\normald_0\cup \normald_1$ in $N_\alpha,$ and adding discs contained in the complement of $\normald_0\cup \normald_1$ on $\partial(N_\alpha)$ as follows. If $\normald_0$ and $\normald_1$ are normal quadrilaterals, then let $X_0$ and $X_1$ be be the connected components of $\simplex^2 \setminus (\normald_0\cup\normald_1)$ containing 1--simplices. If $\normald_0$ and $\normald_1$ are a normal quadrilateral and a normal triangle, let $X_0$ (resp.\thinspace $X_1$) be the component containing a 1--simplex (resp.\thinspace the vertex dual to the triangle). If $\normald_0$ and $\normald_1$ are normal triangles of the same type, let $X_0$ (resp.\thinspace $X_1$) be the component containing a 2--simplex (resp.\thinspace the vertex dual to the triangles). If $\normald_0$ and $\normald_1$ are normal triangles of the different types, let $X_0$ and $X_1$ be the components containing the vertices dual to the triangles. Then let $C = (X_0 \cup X_1) \cap \partial N_\alpha.$ It follows that $((\normald_0 \cup \normald_1)\setminus N_\alpha ) \cup C$ is a union of two disjoint elementary discs which can be straightened to give two disjoint normal discs, denoted by $\normald_0 \uplus \normald_1,$ and the restriction of $C$ to each face coincides with the switch condition at each endpoint of $\alpha$ on that face. For each connected component $\alpha'$ of $\alpha,$ the placement of the discs in $\partial N_{\alpha'}$ is termed the \emph{switch condition at $\alpha'.$}

\begin{figure}[t]
\begin{center}
  \subfigure[]{
      \includegraphics[height=5cm]{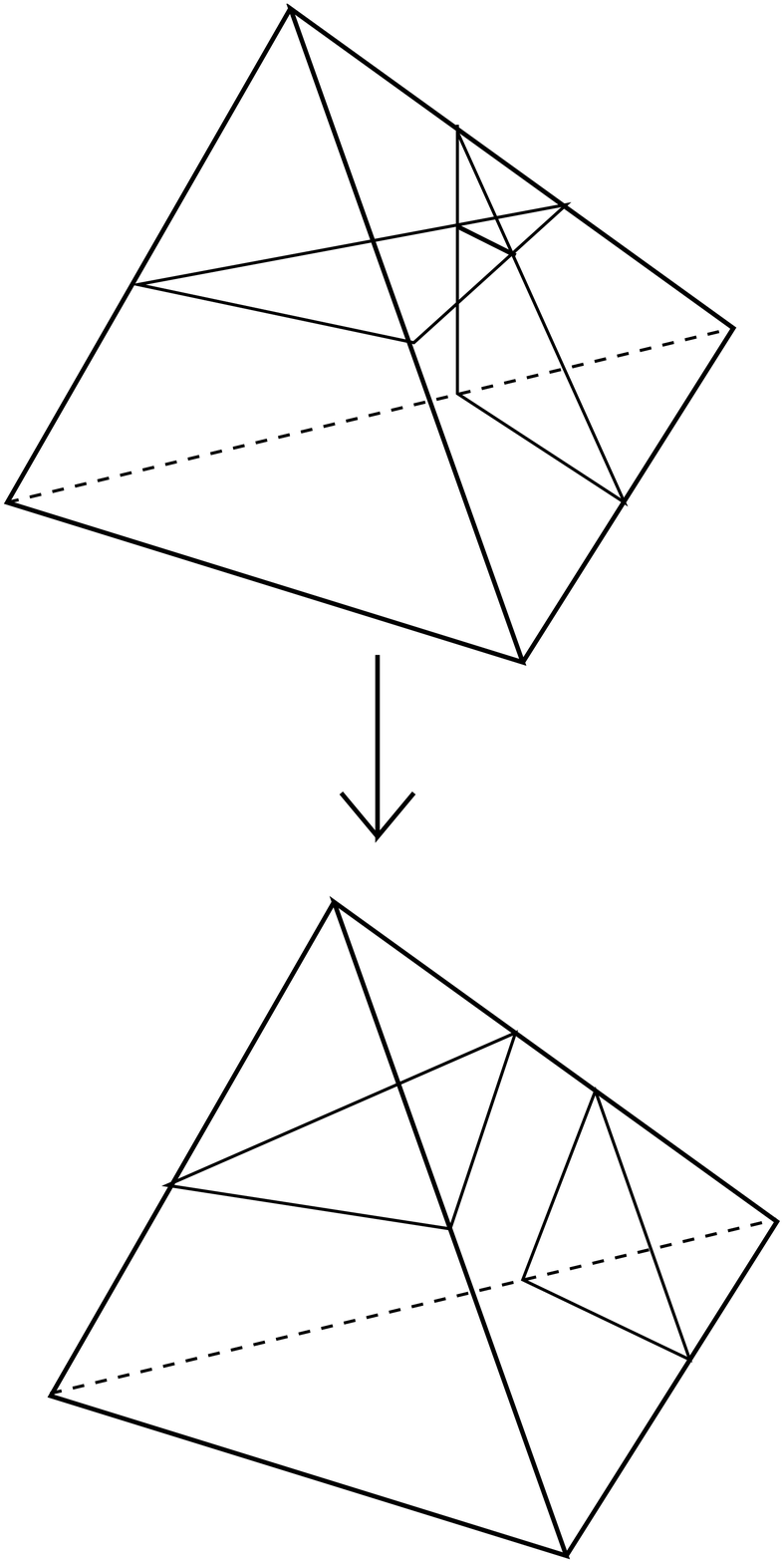}
    } 
    \quad
\subfigure[]{
      \includegraphics[height=5cm]{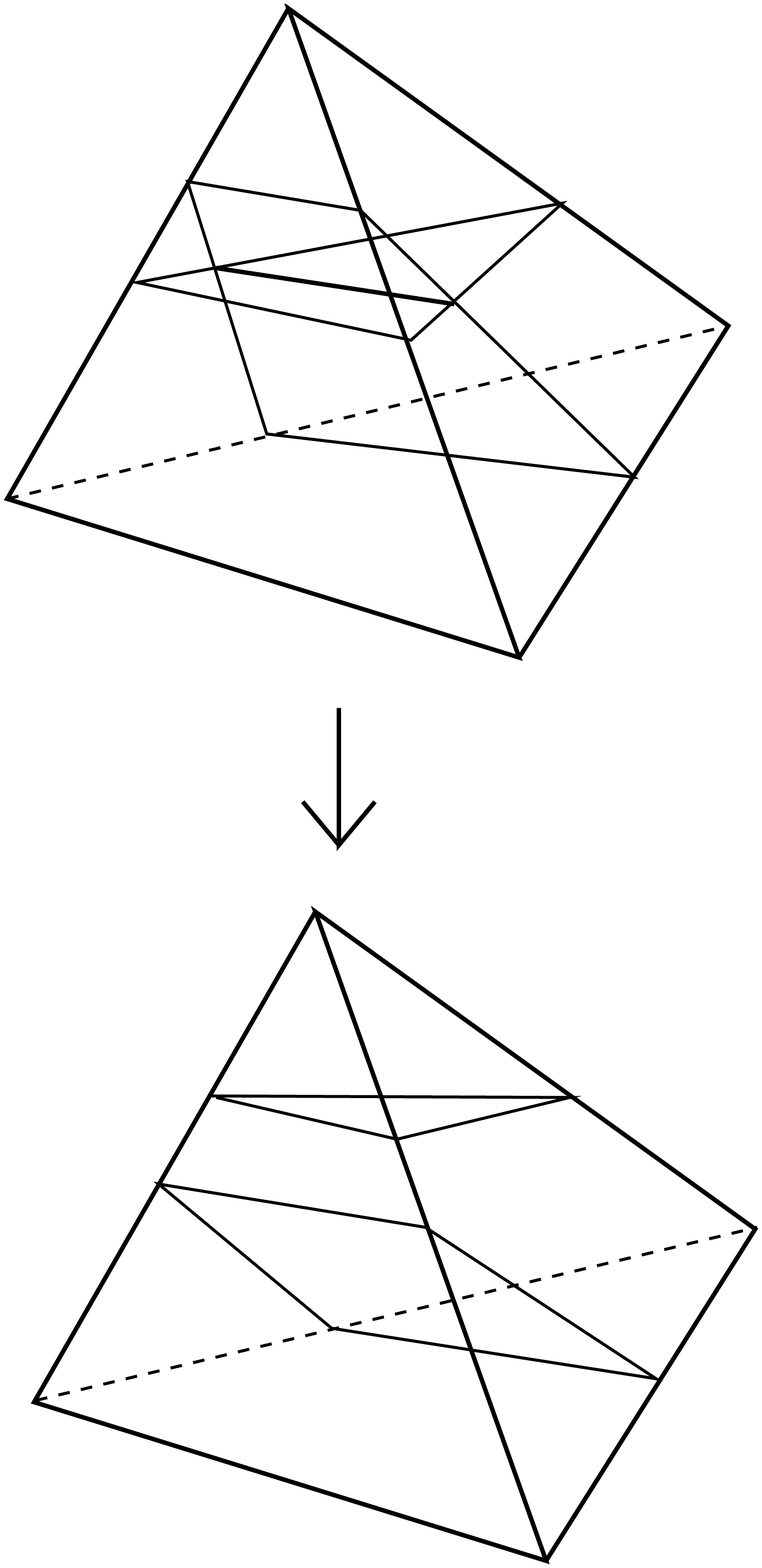}
    } 
    \quad
\subfigure[]{
      \includegraphics[height=5cm]{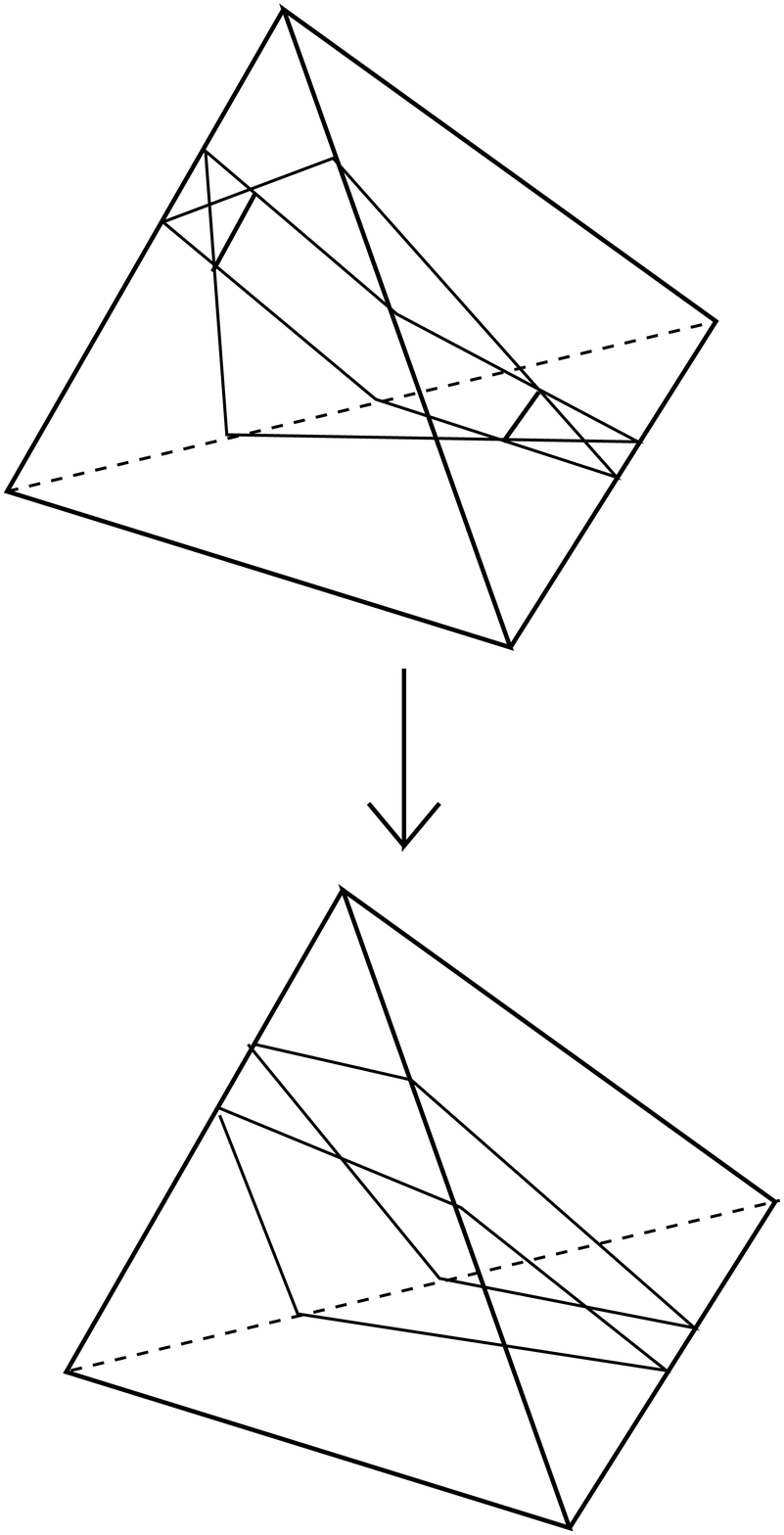}
    } 
    \quad
\subfigure[]{
      \includegraphics[height=5cm]{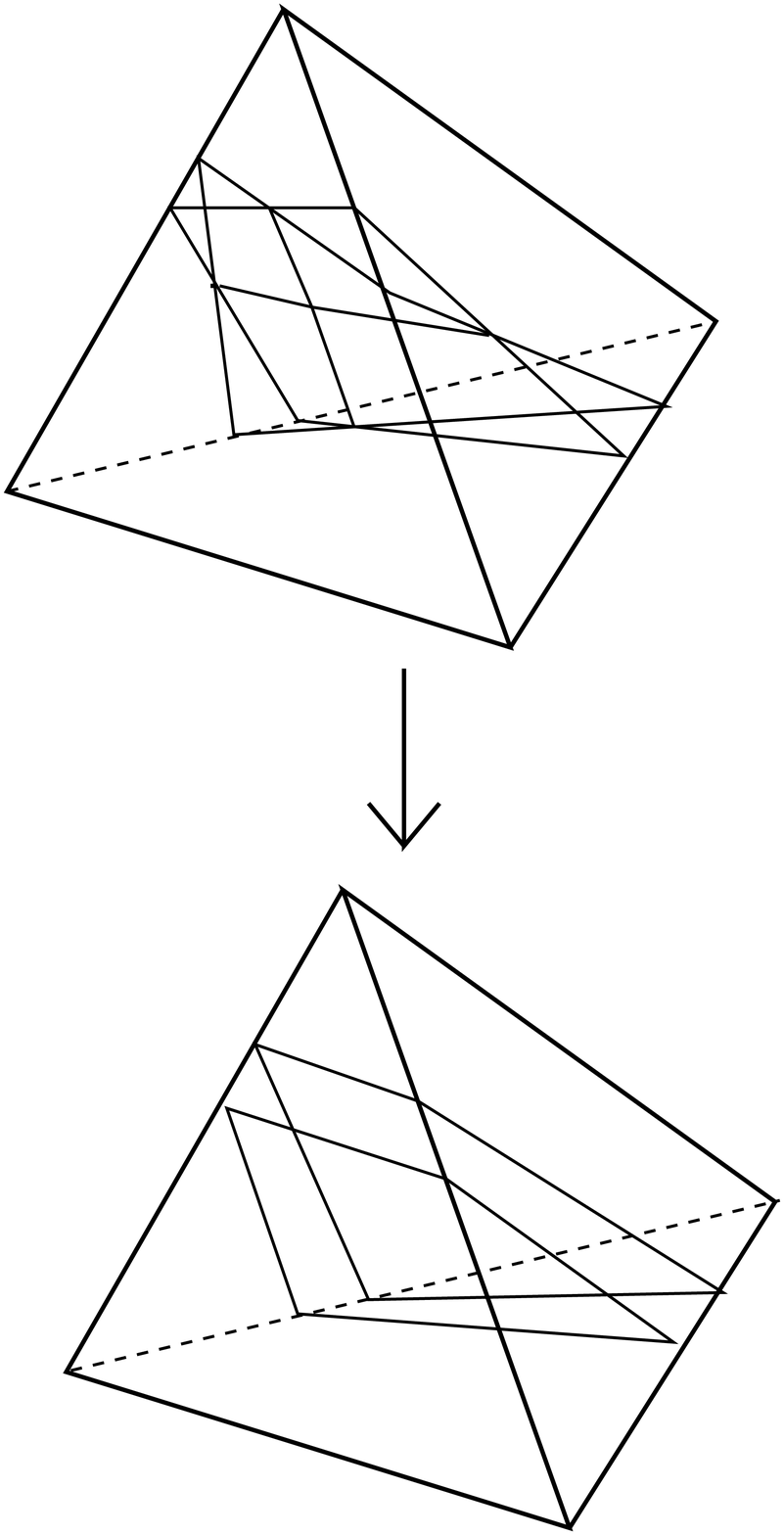}
    } 
\end{center}
    \caption{Regular exchange of normal discs (not shown are the two possible intersections of two quadrilateral discs in one arc)}
     \label{fig:regular exchange of discs}
\end{figure}

\begin{definition}[Compatible normal subsets]
Two normal subsets $S$ and $F$ of $\pseudo$ are \emph{compatible} if they are in general position and if the normal quadrilaterals in $S\cap \simplex^3$ and $F\cap \simplex^3$ are all of the same type for each 3--singlex $\simplex^3$ in $\pseudo.$
\end{definition}

The above definition of regular exchange extends to a regular exchange at $S\cap F$ where $S$ and $F$ are compatible normal subsets. It remains to show that the result (after straightening) is a normal subset. For the proof of the below lemma, it is necessary to extend the above definition of regular exchange to elementary discs $\normald_0$ and $\normald_1$ in the 3--simplex $\simplex^3$ subject to (1) $\normald_0\cap \normald_1 \cap \simplex^{3(1)}=\emptyset,$ (2) if one of the discs has three corners, then they have precisely one arc of intersection with endpoints on different faces of $\simplex^3,$ and (3) if both elementary discs meet the 1--skeleton in four normal corners, then the discs are normally isotopic and have intersection consisting of either one or two arcs or a cross and no two endpoints of $\normald_0\cap \normald_1$ are on the same face of $\simplex^3.$ This generalisation is straight forward; it should be noted that in the realm of elementary discs, saddle tangencies can be eliminated by an arbitrarily small isotopy fixing $\partial \simplex^3.$

\begin{lem}[Geometric sum of compatible normal subsets]\label{lem:regular exchange for discs}
Let $S$ and $F$ be compatible normal subsets of a 3--simplex $\simplex^3.$ Then the regular exchange at $S\cap F$ yields a family of pairwise disjoint elementary discs. It can be straightened to a normal subset of $\simplex^3;$ this is denoted by $S \uplus F.$ Moreover, the number of discs in $S \uplus F$ of any type equals the sum of the number (possibly infinite) of the discs of that type in $S$ and $F.$
\end{lem}

\begin{proof}
The structure of the proof is as in Lemma \ref{lem:regular exchange for arcs}, and only an outline is given. Let $S_0$ be the set of all normal discs in $S$ meeting normal discs in $F$ of a different normal type; define $F_0$ analogously. Then $S_0$ and $F_0$ are both finite (possibly empty) sets of normal discs.

Let $\simplex^1$ be an edge of $\simplex^3$ and denote its vertices by $\v$ and $\w.$ Choose an ordering of the faces $\simplex^2_1$ and $\simplex^2_2$ incident with $\simplex^1.$ A normal triangle dual to $\v$ (resp.\thinspace $\w$) meets each face in a $\v$--type arc (resp.\thinspace $\w$--type arc), and a normal quadrilateral with a corner on $\simplex^1$ meets one of the faces in a $\v$--type and the other in a $\w$--type arc. With respect to the ordering, the corner of a normal disc on $\simplex^1$ is labelled by $(\v, \v),$ $(\v, \w),$ $(\w, \v)$ or $(\w, \w)$ accordingly. Since $S$ and $F$ are compatible, the corners of all quadrilaterals have the same label, say $(\v, \w).$ Read from $\v$ to $\w,$ this gives a word of finite length in these tuples, and if no disc meets a disc of different type, then the word has the form $(\v, \v)^n (\v, \w)^m (\w, \w)^l,$ where $n,m,l \ge 0.$ Given a word not of this form, as before, an inductive argument allows to successively perform regular exchanges which move a pair $(\v, \v)$ to the left or a pair $(\w, \w)$ to the right. Having done this for all edges, it follows that whenever two discs meet, then they are of the same type. The proof is now completed similarly to the proof of Lemma \ref{lem:regular exchange for arcs}.
\end{proof}

\begin{lem}\label{lem:geometric sum is well-defined}
Geometric sum of compatible normal subsets of $\pseudo$ is well-defined and associative.
\end{lem}

\begin{proof}
Let $S$ and $F$ be compatible normal subsets of $\pseudo.$ Then for each 3--simplex $\widetilde{\simplex}^3$ in $\widetilde{\simplex},$ $p^{-1}(S)\cap \widetilde{\simplex}^3$ and $p^{-1}(F)\cap \widetilde{\simplex}^3$ are compatible. It follows that the geometric sum $p^{-1}(S) \uplus p^{-1}(F)$ is well-defined. Since the definition of regular exchange between normal subsets of a 2--simplex is defined without reference to a 3--simplex, it follows that $p\Big(p^{-1}(S) \uplus p^{-1}(F)\Big)$ is a well-defined normal subset of $\pseudo.$
\end{proof}

This has the following, classical, consequence:

\begin{cor}\label{cor:curves in S cap F}
Let $S$ and $F$ be two normal subsets of $\pseudo$ which are in general position, and let $\gamma$ be a simple closed curve in $S\cap F.$ Then $\gamma$ is either 2--sided in both $S$ and $F$ or it is 1--sided in both $S$ and $F.$
\end{cor}

\begin{proof}
If $\gamma$ is non-separating in $S,$ but not in $F,$ then the switch condition at an intersection point of $\gamma$ with the 2--skeleton changes upon a full traverse of $\gamma.$ But this is not possible.
\end{proof}

\begin{lem}[Additive identity]\label{lem:additive identity}
If a normal subset $S$ in $\pseudo$ contains infinitely many normal triangles dual to $\v\in \pseudo^{(0)},$ then $S\uplus B_\v$ is normally isotopic to $S.$
\end{lem}

\begin{proof}
The intersection $S \cap B_\v$ consists of at most finitely many pairwise disjoint simple closed curves and arcs. Since $B_\v$ only contains normal triangles, it follows that $S\uplus B_\v$ and $S$ satisfy the hypothesis of Lemma \ref{lem:uniqueness}.
\end{proof}


\subsection{Boundary curves of spun-normal surfaces}
\label{subsec:Boundary curves of spun-normal surfaces}

Let $S$ be a normal subset in $\pseudo$ which is transverse to $B_\v$ for some $\v \in \pseudo^{(0)}.$ It  may be assumed that $N_\v$ meets only normal triangles in $S$ and that $S$ and $B_\v$ are in general position. Note that $S$ and $B_\v$ are compatible, and that $S \cap B_\v$ is a finite union of pairwise disjoint simple closed curves or properly embedded arcs. Since $B_\v$ is 2--sided in $\pseudo,$ it follows that $S \cap B_\v$ is 2--sided in $S$ and hence, by Corollary \ref{cor:curves in S cap F}, $S \cap B_\v$ is 2--sided in $B_\v.$ This implies that $S \cap B_\v$ can be given a transverse orientation in $B_\v,$ and we will make the following canonical choice. 

Since $\overline{N}_\v$ is the cone over $B_\v$ on $\v,$ it inherits a triangulation which is the cone on $\v$ over the triangulation of $B_\v.$ Moreover, $S\cap \overline{N}_\v$ is a normal subset in $\overline{N}_\v$ with respect to this triangulation. Any normal disc $\normal$ in $\overline{N}_\v$ contained in the 3--singlex $\simplex^3_\v$ in $\overline{N}_\v$ inherits a well-defined transverse orientation by assigning $+1$ to the component of $\Int (\simplex^3_\v) \setminus \normal$ containing $\v.$ These orientations match up to give $S\cap \overline{N}_\v$ a transverse orientation, and $S \cap B_\v$ is given the induced transverse orientation. In particular, this shows that $S\cap \overline{N}_\v$ is 2--sided in $\overline{N}_\v.$

\begin{lem}[Boundary curves are non-separating and 2--sided]\label{lem:characterisation of boundary curves}
Let $S$ be a normal subset in $\pseudo.$ Each curve or arc in $S \cap B_\v$ is 2--sided in $B_\v$ and in $S.$ Up to normal isotopy of $S,$ one may assume that $S \cap B_\v$ is a (possibly empty) union of pairwise disjoint, non-separating, 2--sided simple closed curves or properly embedded arcs on $B_\v.$
\end{lem}

\begin{proof}
It remains to show that $S \cap B_\v$ is non-separating in $B_\v.$ If $S$ contains at most finitely many normal discs dual to $\v,$ then there is a normal isotopy of $S$ making $S$ disjoint from $B_\v,$ and there is nothing to prove. Hence assume that $S$ contains infinitely many normal triangles dual to $\v.$ To simplify notation, let $B=B_\v.$

Let $c$ be a connected component of $S \cap B.$ Perform the regular exchange at $c$ and denote the result by $S(c).$ If $c$ is separating in $B,$ then there are subsets of $S(c)$ corresponding to the components $B_+$ and $B_-$ of $B \setminus c$ and there is a normal homotopy of $S(c)$ which pushes one of these subsets towards $v$ and the other out of $\overline{N}_\v$ so that the result, also denoted by $S(c),$ is transverse to $B$ and meets $B$ in one separating curve fewer than $S.$ Performing regular exchanges along all singular curves in $S(c)$ yields a normal subset $S'(c)$ with the properties that $S'(c)\cap B = S(c)\cap B$ and $S'(c)$ is normally isotopic to $S\uplus B.$ Hence $S'(c)$ is normally isotopic to $S$ by Lemma \ref{lem:additive identity}. Since $S\cap B$ has finitely many connected components, one may, by induction, assume that no component of $S \cap B$ is separating.
\end{proof}

\begin{cor}[Spun-normal surface = normal surface in a manifold]\label{cor:Spun-normal surface = normal surface in a manifold}
If $\pseudo$ is a 3--manifold, then every connected component of a normal subset is a normal surface.
\end{cor}

\begin{proof}
This follows from the above lemma since if $\pseudo$ is a 3-manifold, then every vertex linking surface is a sphere or disc.
\end{proof}


\subsection{The orientable double cover}
\label{subsect:double cover}

If $\pseudo$ is non-orientable, denote by $\widetilde{\pseudo}$ the orientable double cover of $\pseudo$ and by $\widetilde{B}_\v$ a connected component of the pre-image of $B_\v$ for each $\v \in \pseudo^{(0)}.$ A spun-normal surface $S\subset \pseudo$ lifts to a spun-normal surface $\widetilde{S}\subset \widetilde{\pseudo}$ which is invariant under the non-trivial deck transformation. Either $\pseudo$ or $\widetilde{\pseudo}$ is given a fixed orientation. Each boundary component can be given a transverse orientation by choosing a normal vector pointing towards the dual 0--singlex. It is then oriented such that the tuple \emph{(orientation, transverse orientation)} agrees with the orientation of the ambient pseudo-manifold.


\subsection{The boundary curve map}
\label{subsect:boundary curve map}

A transversely oriented normal curve on $B_\v$ defines an element in $H_1(B_\v, \partial B_\v; \mathbb{Z})$ as follows. The triangulation of $B_\v$ can be given the structure of a $\simplex$--complex (see \cite{ha-at}, Section 2.1). Thus, each edge in the triangulation of $B_\v$ is given a well-defined direction, and can be used for a well-defined simplicial homology theory. If $B_\v$ is non-orientable, denote by $\widetilde{B}_\v$ the chosen lift in the double cover; otherwise let $\widetilde{B}_\v=B_\v.$ Recall that $\widetilde{B}_\v$ is oriented. Let $\normala$ be a transversely oriented normal curve on $B_\v,$ and consider its pre-image $\widetilde{\normala}$ in $\widetilde{B}_\v.$ Using the transverse orientation and the orientation of $\widetilde{B}_\v,$ give $\widetilde{\normala}$ an orientation such that the tuple \emph{(orientation, transverse orientation)} associated to $\widetilde{\normala}$ agrees with the orientation of $\widetilde{B}_\v.$ Then homotope $\widetilde{\normala}$ in the direction of the transverse orientation into the 1--skeleton of $\widetilde{B}_\v.$ Using the orientation gives a well-defined element of $Z_1(\widetilde{B_\v}, \partial \widetilde{B_\v}),$ and hence an element of $H_1(\widetilde{B_\v}, \partial \widetilde{B_\v}).$ This maps to a non-torsion element $\overline{\normala} \in H_1({B_\v}, \partial {B_\v})$ since $\normala$ is 2--sided on $B_\v.$ One may therefore view $\overline{\normala} \in H_1({B_\v}, \partial {B_\v}; \mathbb{Z})$ and let $\partial_\v (S) = \sum \overline{\normala}\in H_1({B_\v}, \partial {B_\v}; \mathbb{Z}),$ where the sum is taken over all connected components of $B_\v \cap S.$

\begin{lem}[Boundary curves are well-defined]\label{lem:boundary curves}
Let $S_0$ and $S_1$ be normally isotopic normal subsets of $\pseudo.$ Then $\partial_\v (S_0)=\partial_\v (S_1).$ Moreover, $\partial_\v (S)$ is trivial if and only if $S \cap B_\v$ may be assumed empty.
\end{lem}

\begin{proof}
A dual statement will be proved. It suffices to assume that $B_\v$ is orientable. Let $B_0$ and $B_1$ be two vertex linking surfaces normally isotopic to $B_\v$ and meeting $S$ transversely in normal triangles. The normal isotopy gives a canonical identification $H_1(B_\v; \partial B_\v; \mathbb{Z})\cong H_1(B_k; \partial B_k; \mathbb{Z}).$ Denote the corresponding families of transversely oriented intersection curves by $\mathcal{C}_k.$ Then $B_0$ and $B_1$ bound a submanifold $X$ of $\pseudo$ homeomorphic to $B_\v \times [0,1].$ The surface $S \cap X$ is given a well-defined transverse orientation using the construction in Subsection \ref{subsec:Boundary curves of spun-normal surfaces} which proves that $S \cap N_\v$ is 2--sided. It follows that its set of transversely oriented boundary curves is precisely the union $\mathcal{C}_0 \cup \mathcal{C}_1,$ which therefore determines an element of $B_1({B_\v}, \partial {B_\v}; \mathbb{Z}).$

The argument given in the last paragraph in the proof of Lemma~\ref{lem:characterisation of boundary curves} is now adapted to prove the second part. Let $B=B_\v.$ Assume there is a subset $\mathcal{C}$ of $S \cap B$ with the property that there is a subsurface $B'$ of $B$ with $\partial B' = \mathcal{C}$ and $B'$ is on the positive side of each component of $\mathcal{C}.$ Then the cited argument can be applied with $B_+=B'$ and $B_-=B \setminus B'$ to show that the components in $\mathcal{C}$ can be deleted from $S \cap B$ by a normal isotopy. This completes the proof.
\end{proof}

\begin{cor}
A spun-normal surface $S$ spins into $\v \in \pseudo^{(0)}$ if and only if $\partial_\v (S)\neq 0.$
\end{cor}


\subsection{The spinning construction}

An interesting normal subset in $\overline{N}_\v$ is constructed as follows (this generalises the construction of an infinite normal annulus in \cite{k}). Let $\mathcal{C}$ be a finite collection of pairwise disjoint, transversely oriented, non-separating normal curves on $B_\v$ with the property that the intersection with every normal triangle in $B_\v$ consists of $k$ copies of a normal arc dual to a corner $\normalc_0$ with transverse orientation pointing towards $\normalc_0,$ and $l$ copies of a normal arc dual to $\normalc_1 \neq \normalc_0$ with transverse orientations pointing away from $\normalc_1;$ both $k,l$ are non-negative integers. For instance, the system of transversely oriented normal curves in $S \cap B_\v$ is of this form. 

Let $\{ B_i \},$ $i \in \mathbb{N},$ be a countably infinite family of vertex linking surfaces dual to $v$ in $N_\v$ such that $\cup B_i$ is a normal subset. On each $B_i,$ there is a copy, $\mathcal{C}_i,$ of $\mathcal{C}.$ For each $i$ and for each $c_i \in \mathcal{C},$ delete a small open neighbourhood of $c_i$ in $B_i.$ Since $c_i$ separates its neighbourhood, label the boundary component in the positive direction by $c_i^+,$ and the other by $c_i^-.$ Since there is a normal isotopy of $N_\v$ taking $c_i^-$ to $c_{i+1}^+$ for each $i$ and each $c,$ they can be joined by an annulus in $N_\v.$ Similarly, annuli can be attached between $c$ and $c_1^+$ for each $c \in \mathcal{C}.$ There is a normal isotopy fixing $\overline{N}_\v^{(1)}$ such that the result is a normal subset of $\overline{N}_\v;$ this is said to be obtained by \emph{spinning on $\mathcal{C}.$}

\begin{lem}[Controlled spinning]\label{lem:controlled spinning}
Deleting $S \cap N_\v$ from $S$ and attaching the result of spinning on $S \cap B_\v$ (with the induced transverse orientation) yields a normal subset which is normally isotopic to $S.$
\end{lem}

\begin{proof}
Delete the portion of $S$ in $N_\v$ and replace it by a normal subset in $N_\v$ obtained from the spinning construction on $S \cap B_\v,$ where each component is spun in the $+1$ direction on $B_\v$ according to the convention in the proof of Lemma~\ref{lem:characterisation of boundary curves}. The result meets each 3--singlex in elementary discs and straightening gives a normal subset $S'$ in $P$ with the property that $S$ and $S'$ satisfy the hypothesis of Lemma~\ref{lem:uniqueness}. It follows that $S$ and $S'$ are normally isotopic.
\end{proof}

\begin{cor}\label{cor:topology of spun-normal}
Let $S$ be a spun-normal surface which spins into $\v \in \pseudo^{(0)}.$ Then $\chi(B_\v) \le0.$ If $\chi(B_\v) < 0,$ then $S$ is topologically infinite. If $\chi(B_\v) = 0,$ then $S$ is homeomorphic to the interior of $S \cap (\pseudo \setminus N_\v)$ for suitably chosen $N_\v.$
\end{cor}

\begin{proof}
Assume that $S$ spins into $\v.$ If $\chi(B_\v)>0,$ there there are no 2--sided, non-separating curves or arcs on $B_\v$ and $S$ cannot be spun into $\v.$ Hence $\chi(B_\v) \le 0.$ Since $S$ can be viewed as obtained from the spinning construction, one has 
$$\chi(S) \le \chi(S \cap \pseudo^c) +  \sum_{i=1}^\infty \chi(B_\v)$$
for suitably chosen $\pseudo^c.$ So if $\chi(B_\v) < 0,$ then $S$ is topologically infinite.

If $\chi(B_\v) = 0,$ then $B_v$ may be chosen such that $S \cap B_\v$ contains no collection of curves which is homotopically trivial. It follows that $S \cap \overline{N}_\v$ consists of properly embedded half-open annuli (i.e.\thinspace annuli having one boundary curve removed) or discs with one point on the boundary removed meeting $B_\v$ along a boundary arc.
\end{proof}


\subsection{Kneser's lemma for spun-normal surfaces}

\begin{lem}[Kneser--Haken finiteness]\label{lem:Kneser--Haken finiteness}
Let $\pseudo$ be a triangulated pseudo-manifold (possibly with boundary), and $S$ be a
spun-normal surface such that no two components of $S$ are normally isotopic and no component is a vertex linking surface. Then the number $|S|$ of the components of $S$ satisfies $|S|\le 12 t$ and
\begin{equation}
|S| \le 3 t + \frac{3}{2}\dim H_2(\pseudo^c,\partial \pseudo^c; \mathbb{Z}_2) 
       = 3 t + \frac{3}{2}\dim H_1(\pseudo^c; \mathbb{Z}_2).
\end{equation}
Moreover, if $S$ is 2--sided, then $|S|\le 6 t$ and $|S| \le 3 t + \frac{1}{2}\dim H_2(\pseudo^c,\partial \pseudo^c; \mathbb{Z}_2)$.
\end{lem}

\begin{proof}
Whenever $S$ does not spin into $\v,$ one may add $B_\v$ to $S$ disjointly since no component of $S$ is a vertex linking surface; the result is again denoted by $S.$ The normal discs of $p^{-1}(S)$ divide a 3--simplex in $\widetilde{\simplex}$ into the following types of regions:
\begin{enumerate}
\item \emph{slabs}: trivial $I$--bundles over normal discs, 
\item \emph{thick regions}: truncated tetrahedra and truncated triangular prisms,
\item  \emph{vertex regions}: ``small" tetrahedra which contain a vertex of $\tri$.
\end{enumerate}
In total, there are at most $v + 2t$ components of $\pseudo \setminus S$ which contain the image under $p$ of at least one of the latter types of regions. The remaining components of $\pseudo \setminus S$ are entirely made up of slabs. Each such slab component is either a trivial or a twisted $I$--bundle over a spun-normal surface. If it is a trivial $I$--bundle, then there is a 2--sided surface in $S$ which is the boundary of a twisted $I$--bundle over a 1--sided surface in $S$ (since otherwise there would be two normally isotopic components of $S$). All other slab components are twisted $I$--bundles over a spun-normal surface not in $S$. Since this core surface is not normally isotopic to any of the components of $S$, we may add it disjointly to $S$. Let $F$ be a 1--sided component of $S$. A small regular neighbourhood of $F$ is a twisted $I$--bundle with boundary a 2--sided spun-normal surface, $N,$ in $\pseudo$. If $N$ is not normally isotopic to a component of $S$, we join it to $S$.

Let $\{V_{i}\}$ be the set of all vertex linking surfaces in $S$, $\{F_{i}\}$ be the set of 1--sided surfaces, the set of corresponding 2--sided surfaces be $\{N_{i}\}$, and the set of all remaining 2--sided surfaces be $\{S_{i}\}$. Then $|\{F_{i}\}|=|\{N_{i}\}|.$ Note that $\pseudo^c$ may be chosen such that the intersection of each component of $S$ with it is connected. Let $S^c = S \cap \pseudo^c,$ and similarly for its components. Each $F^c_{i}$ is a 1--sided and therefore non--separating surface in $\pseudo^c.$ Since $F^c_i$ is 1--sided, there is a closed loop in a small regular neighbourhood of $F^c_i$ meeting $F^c_{i}$ transversely in a single point. Using the intersection pairing with $\mathbb{Z}_2$ coefficients, it follows that $F^c_i$ determines a non--zero element of $H_{2}(\pseudo^c,\partial \pseudo^c; \mathbb{Z}_{2})$ and that $\{F_{i}^c\}$ is in bijective correspondence with a linearly independent set of elements of $H_{2}(\pseudo^c,\partial \pseudo^c;\mathbb{Z}_{2})$.

Assume $S'$ is a 2--sided component of $S$ which is not a vertex linking surface and meets no thick region. Then $S'$ bounds an $I$--bundle to either side. If one of them is trivial, then there is a component in $S$ which is normal isotopic to $S'$. Hence both of these $I$--bundles are twisted. But then $\pseudo$ is decomposed into two twisted $I$--bundles glued along their boundaries, and $S=S'$ does not meet any thick region. This is not possible. So each element of $\{N_{i}\}\cup \{S_{i}\}$ meets at least one thick region. There are at most $6t$ normal discs in the boundaries of thick regions. Each surface in $\{N_{i}\}\cup \{S_{i}\}$ meets at least one thick region in at least one disc. This implies that $|\{N_{i}\}| + |\{S_{i}\}| \le 6t.$ Thus $|\{F_{i}\}|\le 6t,$ and the first bound is obtained. We now use a doubling trick from \cite{bach}. Each $N_{i}$ meets thick regions only on one side, and each $S_{i}$ meets at least one thick region on each side. Push each $N_{i}$ off itself away from the twisted $I$--bundle, and call the resulting copy $N'_{i}$, and push each $S_{i}$ off itself and call the resulting copy $S'_{i}$. This can be done such that all surfaces are still pairwise disjoint.  Each surface in $\{N'_{i}\} \cup \{S_{i}\} \cup \{S'_{i}\}$ meets at least one thick region in at least one disc. Thus,
$2|S| = 2|\{V_{i}\}|+ 2|\{F_{i}\}|+ 2|\{N_{i}\}|+2|\{S_{i}\}| 
        = 2|\{V_{i}\}|+ 3|\{F_{i}\}|+ |\{N_{i}\}|+2|\{S_{i}\}|
        \le 2  |\{V_{i}\}| + 3 \dim H_2(\pseudo^c,\partial \pseudo^c; \mathbb{Z}_2) + 6 t$.
Dividing by two and subtracting the vertex linking surfaces gives the second inequality. The stated equality follows since Lefschetz duality and the universal coefficient theorem yield $H_2(\pseudo^c,\partial \pseudo^c; \mathbb{Z}_2)\cong H_1(\pseudo^c; \mathbb{Z}_2).$ 
\end{proof}


\subsection{Normalising properly embedded surfaces}

Haken's approach to normalising properly embedded \emph{closed} surfaces as described in Chapter 3 of \cite{M} applies to ideal triangulations without change. However, it does not apply to properly embedded \emph{non--compact} surfaces; an obstruction for putting such a surface into spun-normal form is, for instance, given by Corollary \ref{cor:topology of spun-normal}. Moreover, it is shown by Kang \cite{k} that there is an ideal triangulation of the complement of the figure eight knot with two ideal 3--simplices, such that no Seifert surface for the knot can be put into spun-normal form (see also Section \ref{The figure eight knot}). A result due to Thurston (see Walsh \cite{w}) states that essential surfaces other than virtual fibres can be normalised in any hyperbolic 3--manifold with torus cusps and ideal triangulation with essential edges. A general theory of normalisation with respect to ideal triangulations is implicit in the work of Brittenham and Gabai \cite{g}.


\section{Matching equations}
\label{sec:Matching equations}

A normal subset $S$ is (up to vertex linking components and normal isotopy) uniquely determined by its quadrilateral discs. Recording the number of quadrilaterals of each type gives the normal $Q$--coordinate of $S.$ This satisfies two necessary conditions: (1) it is admissible, and (2) it satisfies a \emph{$Q$--matching equation} for each edge not contained in the boundary of $\pseudo.$ This equation is given by Tollefson \cite{to} for compact 3-manifolds. The main result of this section (Theorem~\ref{thm:normal unique}) states that, conversely, any admissible solution to the $Q$--matching equation is realised by a spun-normal surface which is unique up to normal isotopy.


\subsection{Normal Q--coordinate}
\label{subsec:Normal Q--coordinates}

There are $t$ 3--singlices in $\pseudo.$ Let $\{ q_1,...,q_{3t}\}$ be the set of all quadrilateral types in $\pseudo.$ A normal subset $S$ in $\pseudo$ meets each 3--singlex in at most finitely many quadrilateral discs of each type. If $x_i$ is the number of normal quadrilaterals of type $q_i,$ then let $N(S) = (x_1,{\ldots} ,x_{3t})$ be the \emph{normal $Q$--coordinate} of $S.$ Thus, $N(S)$ is a well-defined point in $\RR^{3t}$ with the coordinate axes labelled by the quadrilateral types. A point $(x_1,{\ldots} ,x_{3t}) \in \RR^{3t}$ is called \emph{admissible} if each $x_i \ge 0$ and if $q_i, q_j, q_k$ are the three distinct quadrilateral types contained in any 3--singlex, then at most one of $x_i, x_j, x_k$ is non-zero. The cell structure of $S$ allows normal quadrilaterals of at most one quadrilateral type in each ideal 3--singlex, and hence $N(S)$ is admissible. 


\subsection{Convention for oriented 3--simplices}
\label{subsec: convention for tets}

\begin{figure}[h]
\psfrag{D}{{\small $\Delta$}}
\psfrag{z}{{\small $z$}}
\psfrag{z'}{{\small $z'$}}
\psfrag{z"}{{\small $z''$}}
\psfrag{q}{{\small $q$}}
\psfrag{q'}{{\small $q'$}}
\psfrag{q"}{{\small $q''$}}
\begin{center}
  \includegraphics[width=12cm]{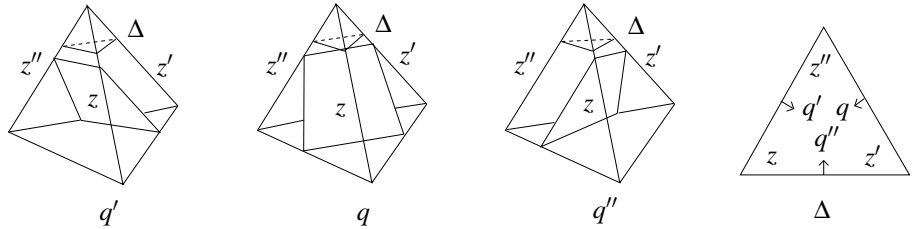}
\end{center}
    \caption{Edge labels, quadrilateral types and $Q$--moduli}
    \label{fig:Delta labels}\label{fig:parameters}\label{fig:quad coods}
\end{figure}

Let $\simplex^3$ be an oriented, regular Euclidean 3--simplex. The edges of $\simplex^3$ are labelled with parameters $z,z',z'',$ such that opposite edges have the same parameter, and the ordering $z,z',z''$ agrees with a right--handed screw orientation of $\simplex^3;$ this is pictured in Figure \ref{fig:parameters}. It follows that the labelling is uniquely determined once the parameter $z$ is assigned to any edge of $\simplex^3.$ The vertices of a normal triangle dual to a vertex of $\simplex^3$ inherit moduli from the edge parameters; this labelling is always viewed from the vertex.

There are three quadrilateral types in $\simplex^3;$ let $q^{(k)}$ denote the quadrilateral type which does not meet the edges labelled $z^{(k)}$ for $k \in \{0,1,2\}.$ The symmetry group of $\simplex^3$ is the alternating group on its four vertices. It contains a normal Kleinian four group which leaves pairs of opposite edges invariant and hence fixes quadrilateral types. The quotient group is $\mathbb{Z}_3$; it may be identified with the group of even permutations of the three quadrilateral types, and is hence generated by the cycle $\pi = (q, q', q'').$ If $p$ and $q$ are quadrilateral types in $\simplex^3,$ let:
\begin{equation*}
s_p(q)=k \ \text{if} \ \pi^k(p)=q \ \text{for} \ k \in \{-1,0,1\}.
\end{equation*}
An isometry $\varphi\co \simplex^3_0 \to \simplex^3_1$ between oriented 3--simplices induces a bijection between quadrilateral types, and one has $s_p(q)= s_{\varphi(p)}(\varphi(q))$ if $\varphi$ is orientation preserving, and $s_p(q)= -s_{\varphi(p)}(\varphi(q))$ otherwise. The above extends to normal quadrilateral types in an oriented 3--singlex since the quotient map from a 3--simplex to a 3--singlex induces a bijection between normal quadrilateral types.


\subsection{Degree and abstract neighbourhood}

The \emph{degree of a 1--singlex} $\simplex^1$ in $\pseudo,$ $\deg(\simplex^1),$ is the number of 1--simplices in $\widetilde{\simplex}$ which map to $\simplex^1.$ Given a 1--singlex $\simplex^1$ in $\pseudo,$ there is an associated \emph{abstract neighbourhood $B(\simplex^1)$} of $\simplex^1$ which is a ball triangulated by $\deg (\simplex^1)$ 3--simplices with the property that there is a well-defined simplicial quotient map $p_{\simplex^1}\co B(\simplex^1)\to \pseudo$ taking $\widetilde{\simplex}^1$ to $\simplex^1.$

If $\simplex^1$ has at most one pre-image in each 3--simplex in $\widetilde{\simplex},$ then $B(\simplex^1)$ is obtained as the quotient of the collection $\widetilde{\simplex}_{\simplex^1}$ of all 3--simplices in $\widetilde{\simplex}$ containing a pre-image of $\simplex^1$ by the set $\Phi_{\simplex^1}$ of all face pairings in $\Phi$ between faces containing a pre-image of $\simplex^1.$ There is an obvious quotient map $b_{\simplex^1}\co B(\simplex^1)\to \pseudo$ which takes into account the remaining identifications on the boundary of $B(\simplex^1).$

We now describe $B(\simplex^1)$ if $\simplex^1$ has more than one pre-image in some 3--simplex. In this case $\deg(\simplex^1)>1.$ For each 3--simplex $\widetilde{\simplex}^3_j,$ if $\simplex^1$ has $k$ pre-images in $\widetilde{\simplex}^3_j,$ take $k$ copies of this 3--simplex, $\widetilde{\simplex}^3_{j,1},...,\widetilde{\simplex}^3_{j,k},$ each with an isometry $i_{j,h} \co \widetilde{\simplex}^3_j \to \widetilde{\simplex}^3_{j,h}.$ To each pre-image $\widetilde{\simplex}^1$ of $\simplex^1$ in $\widetilde{\simplex}^3_j$ assign one of the maps $i_{j,h}$ and call $i_{j,h}(\widetilde{\simplex}^1)$ the \emph{marked edge} of $\widetilde{\simplex}^3_{j,h}.$ This gives a collection $\widetilde{\simplex}_{\simplex^1}$ of $\deg (\simplex^1)$ 3--simplices.

A face pairing $\varphi$ of faces of $\widetilde{\simplex}^3_j$ and $\widetilde{\simplex}^3_h$ gives rise to a face pairing between faces of $\widetilde{\simplex}^3_{j,i}$ and $\widetilde{\simplex}^3_{h,k}$ if and only if it identifies the edges of $\widetilde{\simplex}^3_j$ and $\widetilde{\simplex}^3_h$ corresponding to the marked edges of $\widetilde{\simplex}^3_{j,i}$ and $\widetilde{\simplex}^3_{h,k}.$ Denote the resulting set of face pairings by $\Phi_{\simplex^1}.$

The quotient map is denoted by $p_{\simplex^1}\co \widetilde{\simplex}_{\simplex^1}  \to \widetilde{\simplex}_{\simplex^1} / \Phi_{\simplex^1}=B(\simplex^1).$ If $\deg(\simplex^1)>1$ or if $\simplex^1\subset \partial \pseudo,$ then the triangulation of $B(\simplex^1)$ is simplicial.

We now describe the map $b_{\simplex^1}$ in the case where $\simplex^1$ has more than one pre-image in some 3--simplex. Since $\deg(\simplex^1)>1,$ the given triangulation of $B(\simplex^1)$ is simplicial and hence for each simplex $\widetilde{\simplex}^3_{j,h}$ in $B(\simplex^1)$ the inverse $p_{\simplex^1}^{-1}$ is well-defined. Since the face pairings defining $B(\simplex^1)$ arise from the face pairings in $\Phi,$ the map $p_{\simplex^1}\co B(\simplex^1)\to \pseudo$ defined by $p\circ i_{j,h}^{-1}\circ p_{\simplex^1}^{-1} \co \widetilde{\simplex}^3_{j,h} \to \pseudo$ for each 3--simplex is well-defined.

The unique edge in $B(\simplex^1)$ which is the image of all marked edges is termed the \emph{axis} of  $B(\simplex^1).$ One of the endpoints of the axis is termed the \emph{north pole} and the other the \emph{south pole}. The set of all 1--singlices not containing a pole are termed \emph{the equator} of $B(\simplex^1).$ See Figure \ref{fig:slopes}(a).

\begin{lem}
Let $S$ be a normal subset of $\pseudo.$ Then $p_{\simplex^1}^{-1}(S)$ is a normal subset of $B(\simplex^1).$
\end{lem}

\begin{proof}
The set $p^{-1}(S)$ is a normal subset of $\widetilde{\simplex}$ with the property that whenever $\normala \subset p^{-1}(S)$ is a normal arc on a 2--simplex in the range of some face pairing $\varphi \in \Phi,$ then $\varphi(\normala) \subset p^{-1}(S).$ The result now follows from the definitions of $\widetilde{\simplex}_{\simplex^1}$ and $\Phi_{\simplex^1}.$
\end{proof}

\begin{lem}
Any normal subset in $B(\simplex^1)$ is a (possibly infinite) union of pairwise disjoint properly embedded discs in $B(\simplex^1).$
\end{lem}

\begin{proof}
This follows from the fact that $B(\simplex^1)$ is a manifold and Corollary \ref{cor:Spun-normal surface = normal surface in a manifold}.
\end{proof}


\subsection{Slopes of quadrilateral types}

\begin{figure}[t]
\psfrag{z}{{\small $z$}}
\psfrag{z'}{{\small $z'$}}
\psfrag{z"}{{\small $z''$}}
\psfrag{q}{{\small $q$}}
\psfrag{q'}{{\small $q'$}}
\psfrag{q"}{{\small $q''$}}
\psfrag{e}{{\small $e'$}}
\psfrag{B(e)}{{\small $B(e)$}}
\psfrag{equator}{{\small equator}}
\psfrag{a}{{\small $$}}
\psfrag{b}{{\small $$}}
\begin{center}
  \subfigure[Abstract neighbourhood]{
    \begin{minipage}[b]{5cm}
      \centering
      \includegraphics[height=3.5cm]{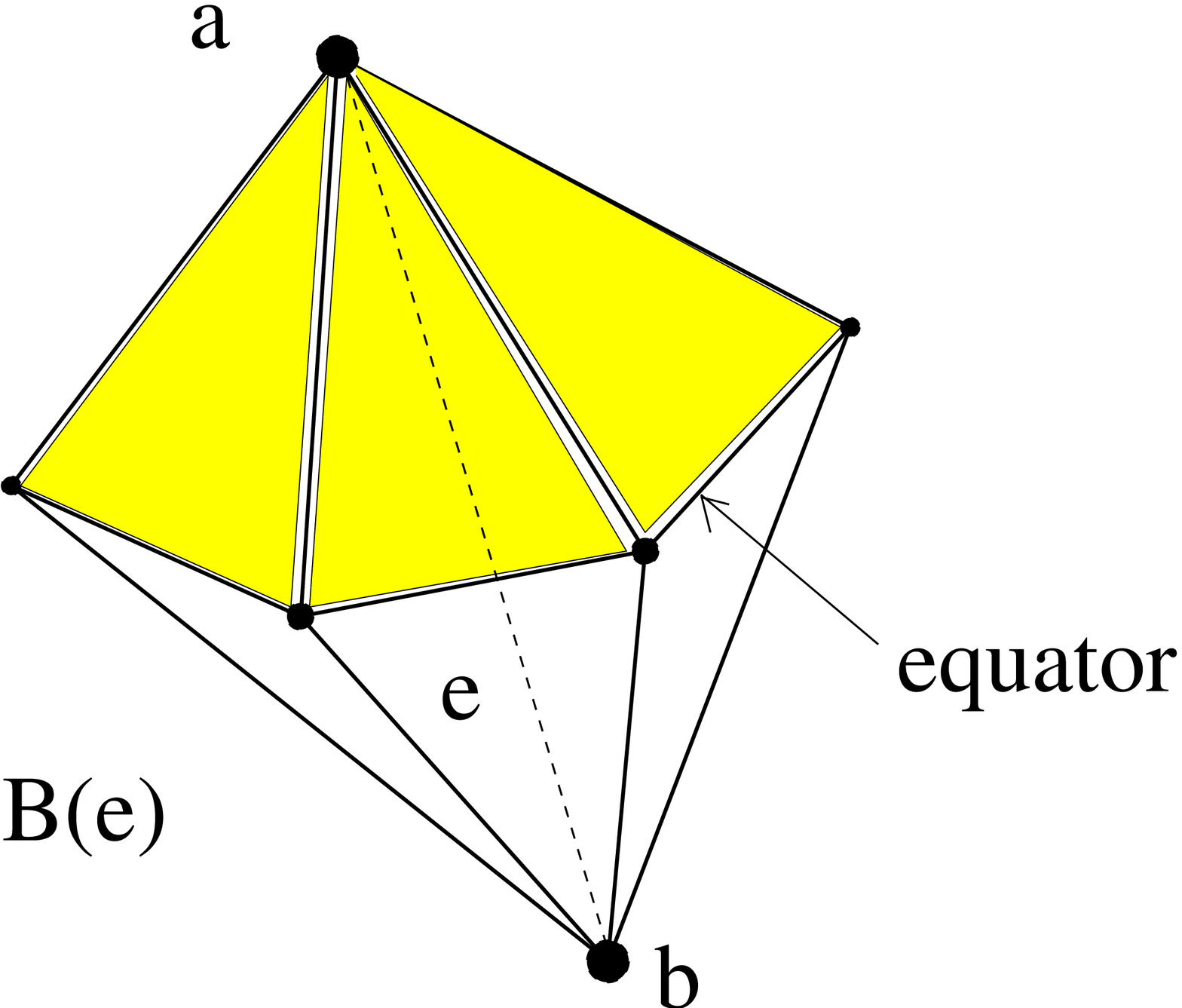}
    \end{minipage}} 
 \subfigure[positive (left) and negative (right) slope]{
    \begin{minipage}[b]{7.3cm}
      \centering
      \includegraphics[width=2.4cm, height=3.5cm]{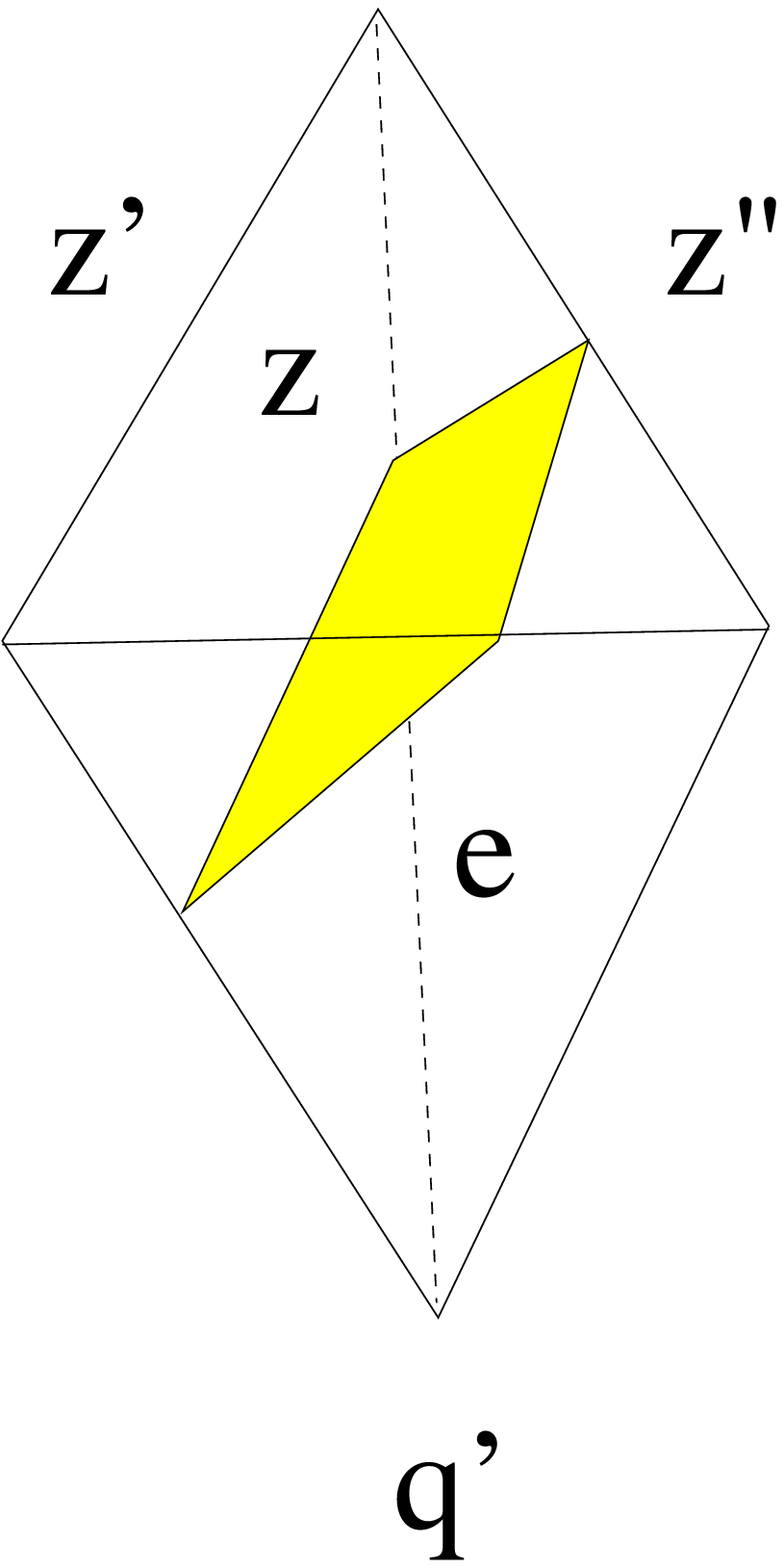}
      \qquad
      \includegraphics[width=2.4cm, height=3.5cm]{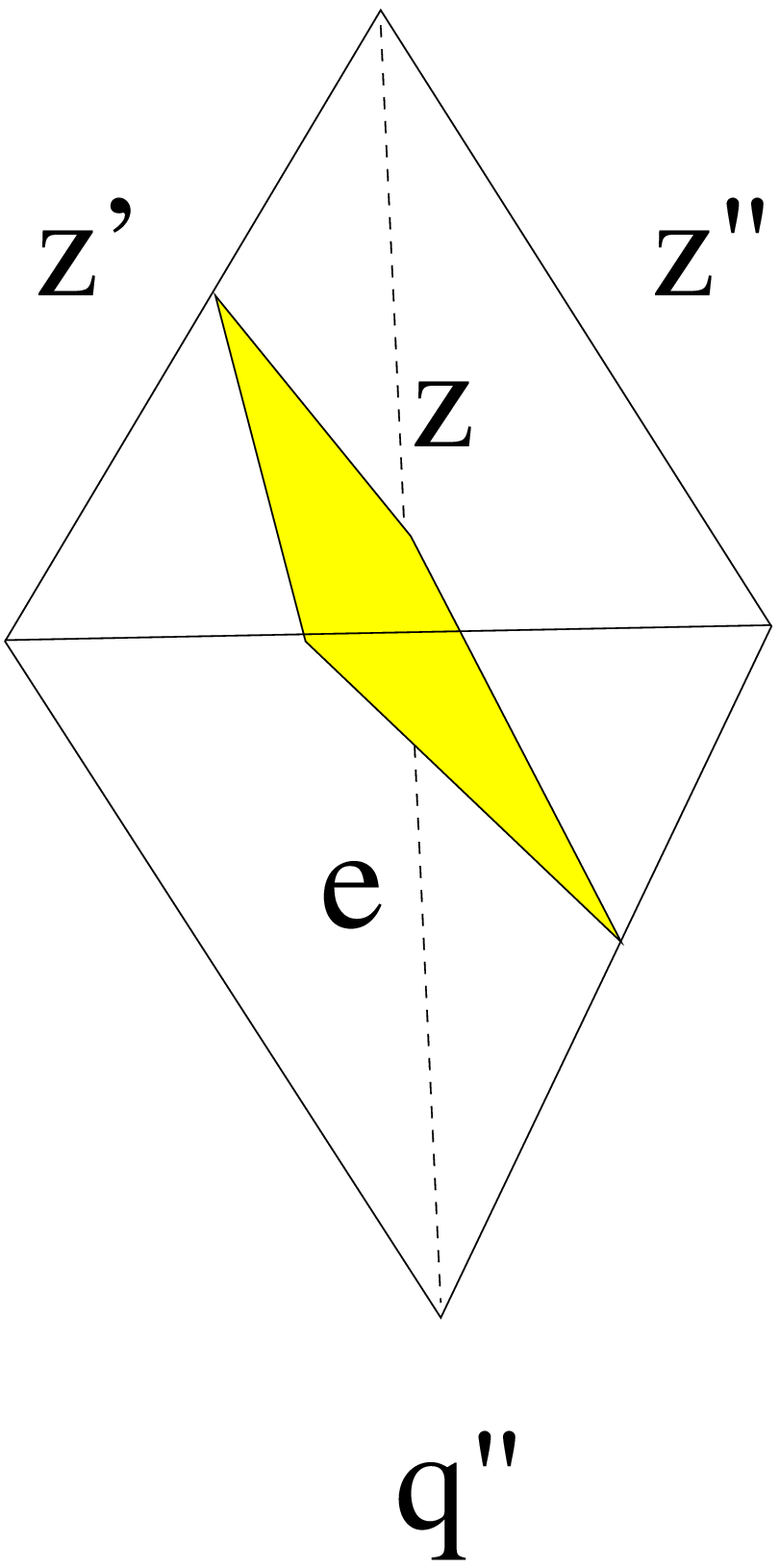}
    \end{minipage}}
\end{center}
    \caption{Abstract neighbourhood and slopes of quadrilaterals}
     \label{fig:Q_corner_gl}\label{fig:slopes}
\end{figure}

Let $\simplex^1$ be a 1--singlex not contained in $\partial \pseudo.$ The abstract neighbourhood $B(\simplex^1)$ of a 1--singlex $\simplex^1$ in $\pseudo$ is a ball and hence the 3--singlices in $B(\simplex^1)$ may be oriented coherently. In particular, the convention of Subsection \ref{subsec: convention for tets} applies to each oriented 3--singlex in $B(\simplex^1).$ To simplify notation, write $e=\simplex^1.$ If $p$ is the type of a normal quadrilateral in the oriented abstract neighbourhood $B(e)$ which does not meet the axis $e',$ let $s_e(q) = s_p(q)$ for any quadrilateral type $q$ in the same 3--singlex as $p.$ This gives a function $s_e$ defined on the set of all quadrilateral types in $B(e)$, and the value $s_e(q)$ is termed \emph{the slope of $q$ w.r.t.\thinspace $e.$} The two quadrilateral types in a 3--singlex meeting $e'$ have slopes of opposite signs, and a quadrilateral type which does not meet $e'$ has slope zero. See Figure \ref{fig:Q_corner_gl}(b), where the shown 3--singlex is oriented using a right-handed screw orientation.

The orientation of $B(e)$ induces an orientation of its hemispheres, and the equator will be given the orientation induced from the southern hemisphere.

If $D \subset B(e)$ is a properly embedded disc such that the poles are contained in different components of $\partial B(e) \setminus \partial D,$ then $\partial D$ is given the orientation induced from the component containing the south pole. If $P$ is a point in the intersection of $\partial D$ and the equator, let $s(P)=1$ if $\partial D$ crosses the equator from the southern to the northern hemisphere at $P$, and let $s(P)=-1$ otherwise. Thus $0=\sum s(P),$ where the sum is taken over all points in the intersection of $\partial D$ and the equator; the empty sum is defined to be zero throughout. 

Now assume that $D$ is a connected normal surface in $B(e)$ which meets $e'.$ Then $D$ is a disc; it is the join of $\partial D$ to $D \cap e'.$ In particular, $\partial D$ is a normal curve separating the poles on $\partial B(e)$ and $D$ meets each 3--singlex in exactly one normal disc. The boundary of a normal triangle in $D$ does not meet the equator, and a quadrilateral disc of type $q$ in $D$ contributes two normal arcs to its boundary which meet in a point $P$ at the equator. The above conventions imply that  $s(P)=s_e(q).$ Thus, $\sum s_e(q)=0,$ where the sum is taken over all quadrilateral types in $D.$ If $D$ is any normal subset in $B(e)$ which does not meet $e'$, this equation is trivially satisfied since a quadrilateral disc meets the equator if and only if it has a vertex on $e'.$

Note that the equation is well-defined up to sign since reversing the orientation of $B(\simplex^1)$ changes the sign of the slope of each quadrilateral meeting $\simplex^1.$ If $\simplex^1$ is contained in $\partial \pseudo,$ then define $s_{\simplex^1}(q)=0$ for each quadrilateral type in $B(\simplex^1).$

\begin{figure}[t]
\psfrag{a}{}
\psfrag{b}{}
\psfrag{c}{}
\psfrag{d}{}
\psfrag{e}{{\small $\sum s_{\simplex^1}(q) x_q=0$}}
\psfrag{f}{{\small $\sum s_{\simplex^1}(q) x_q\neq 0$}}
    \begin{center}
      \includegraphics[width=10cm]{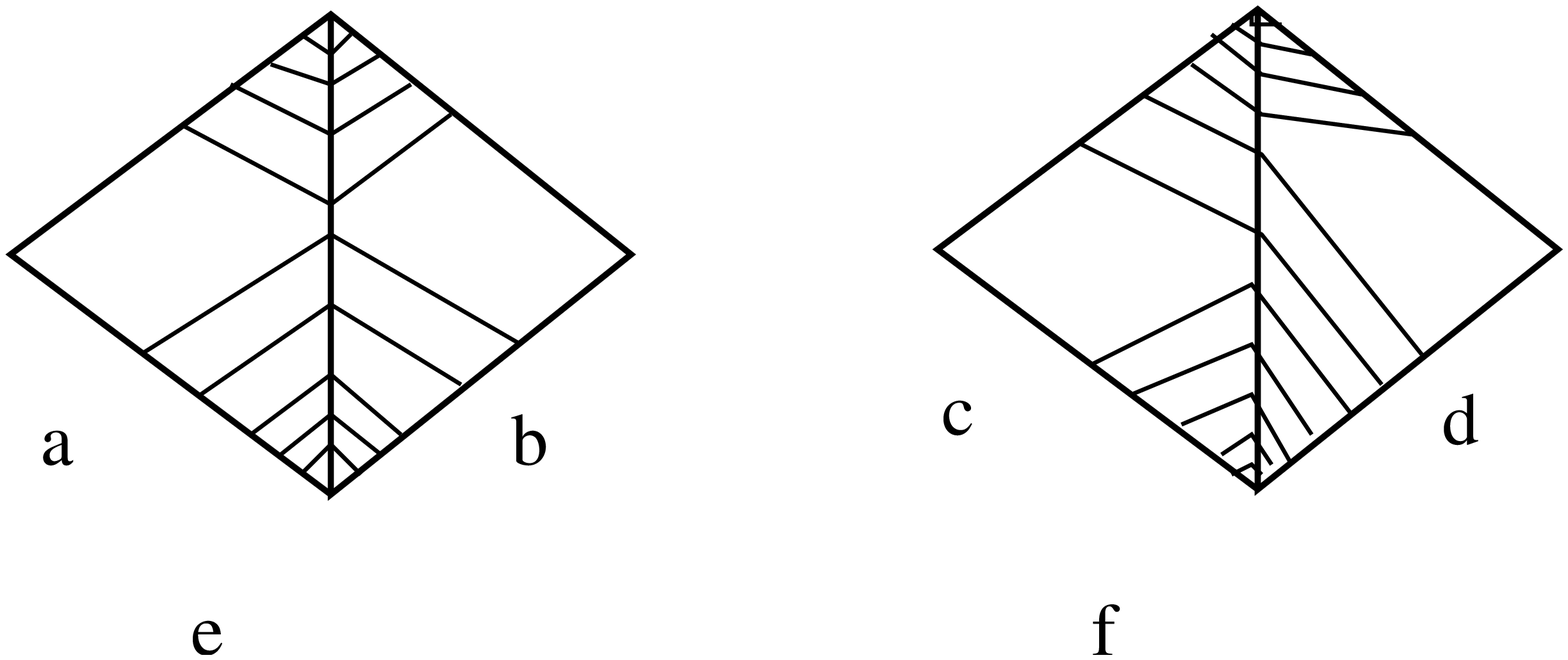}
    \end{center}
\caption{$B(\simplex^1)$ cut open along $p_{\simplex^1}(\widetilde{\simplex}^2)$}
\label{fig:lemma23}
\end{figure}

\begin{lem}\label{lem:prep normal unique}
A normal subset $\tilde{S} \subset \tilde{\simplex}_{\simplex^1}$ which contains infinitely many normal triangles of each type can be normally isotoped such that $p_{\simplex^1}(\tilde{S})$ is a normal subset of $B(\simplex^1)$ if and only if $\sum s_{\simplex^1}(q) x_q=0,$ where $x_q$ is the number of normal quadrilaterals of type $q$ in $p_{\simplex^1}(\tilde{S}).$
\end{lem}

\begin{proof}
It follows from the above discussion that the equation is a necessary condition, and it remains to show that it is sufficient. Given two 3--simplixes in $\tilde{\simplex}_{\simplex^1},$ the pattern on each face is the same: there are three countably infinite stacks of parallel copies of normal arcs. Fix the discs of $\tilde{S}$ in one 3--simplex $\widetilde{\simplex}^3.$ If there are no face pairings for faces of $\widetilde{\simplex}^3,$ then $B(\simplex^1)=\widetilde{\simplex}^3$ and there is nothing to prove.

Next assume that there is precisely one face pairing, $\varphi_1,$ involving faces of $\widetilde{\simplex}^3.$ Then the equation is trivial. Let $\widetilde{\simplex}^3_2$ be a simplex with one face $\widetilde{\simplex}^2$ in the range or domain of $\varphi_1.$ Then there is an ideal normal isotopy of $\widetilde{\simplex}^3_2$ which fixes all subsimplices not in $\widetilde{\simplex}^2$ such that $\varphi_1$ matches the normal arcs on $\widetilde{\simplex}^2$ bijectively with the normal arcs on its counterpart. If there are no other face pairings involving $\widetilde{\simplex}^3_2,$ we are done. Otherwise there is exactly one other face pairing, and the above procedure can be iterated. It terminates after $\deg (\simplex^1)-1$ steps.

Hence assume that there are two face pairings for the faces of every 3--singlex. Then the equation is non-trivial. Ignoring one of the face pairings, $\varphi_0,$ for faces of $\widetilde{\simplex}^3,$ one can proceed as in the previous paragraph. Denote by $\widetilde{\simplex}^3_k$ the last simplex in the iteration, and assume that $\varphi_0 (\widetilde{\simplex}^2) = \widetilde{\simplex}^2_k,$ where $\widetilde{\simplex}^2$ is a face of $\varphi_0.$ It follows from the construction that $\tilde{S} \cap \widetilde{\simplex}^3_k$ can only be altered by a normal isotopy which fixes all edges of $\widetilde{\simplex}^3_k$ other than the non-marked edges of $\widetilde{\simplex}^2_k;$ see Figure~\ref{fig:lemma23}. To obtain a normal subset of $B(\simplex^1),$ the normal arcs dual to the vertices not contained on the marked edges on $\widetilde{\simplex}^2$ and $\widetilde{\simplex}^2_k$ have to be identified in pairs, and this determines a unique bijection between the normal corners on the edges which are not marked. This bijection determines a unique normal isotopy of $\tilde{S} \cap \widetilde{\simplex}^3_k$ such that equivalent corners have the same image under $p_{\simplex^1}.$ Note that $p_{\simplex^1}(\tilde{S}) \cap \partial B(\simplex^1)$ is a normal subset.

If $p_{\simplex^1}(\tilde{S})$ is a normal subset, then $p_{\simplex^1}(\tilde{S}) \cap \partial B(\simplex^1)$ is a countably infinite union of circles and the equation $\sum s_{\simplex^1}(q) x_q=0,$ is satisfied. Assume that $p_{\simplex^1}(\tilde{S})$ is not a normal subset. Then some component of $p_{\simplex^1}(\tilde{S}) \cap \partial B(\simplex^1)$ crosses the equator an odd number of times, and, since $p_{\simplex^1}(\tilde{S}) \cap \partial B(\simplex^1)$ is properly embedded in $\partial B(\simplex^1)$ minus the poles, it follows that $\sum s_{\simplex^1}(q) x_q\neq 0.$ This proves the lemma.
\end{proof}


\subsection{The Q--matching equation}

For each 1--singlex $\simplex^1$ in $\pseudo,$ the \emph{total slope with respect to $\simplex^1,$ $s_{\simplex^1}(q),$} of a quadrilateral type $q$ in $\pseudo$ is defined to be the sum over all pre-images of $q$ in $B(\simplex^1)$ of the signs of their slopes. If $S$ is a normal subset of $\pseudo,$ then it follows from Lemma \ref{lem:prep normal unique} that the normal $Q$--coordinate $N(S)=(x_1,...,x_{3t})$ of $S$ satisfies a linear equation for $\simplex^1,$ called the \emph{$Q$--matching equation of $\simplex^1:$}
\begin{equation*} 
0 = \sum_{i=1}^{3t} s_{\simplex^1}(q_i) x_i.
\end{equation*}
This equation is well-defined up to sign since reversing the orientation of $B(\simplex^1)$ changes the sign of the slope of each quadrilateral meeting $\simplex^1.$

\begin{thm}\label{thm:normal unique}
Let $\pseudo$ be a triangulated pseudo-manifold. For each admissible integer solution $N$ of the $Q$--matching equations there exists a spun-normal surface $S$ in $\pseudo$ with no vertex linking components such that $N(S)=N.$ Moreover, $S$ is unique up to normal isotopy.
\end{thm}

\begin{proof}
It suffices to show existence; uniqueness follows from Lemmata~\ref{lem:vertex linking surface} and \ref{lem:uniqueness}.

Recall that there is a bijective correspondence between normal disc types in $\widetilde{\simplex}$ and $\pseudo.$ Given an admissible integer solution $(x_1,...,x_{3t}),$ place $x_i$ pairwise disjoint normal discs in $\widetilde{\simplex}$ each of which maps to a normal disc of type $q_i$ in $\pseudo.$ Then place infinitely many normal triangles of each type in each 3--simplex in $\widetilde{\simplex}$ such that the resulting collection of normal discs is a normal subset, denoted by $\widetilde{S},$ of $\widetilde{\simplex}.$ It suffices to show that there is a normal isotopy $i$ of $\widetilde{\simplex}$ such that $p(i(\widetilde{S}))$ is a normal subset of $\pseudo.$

The normal isotopy is uniquely determined by defining it on the 1--skeleton. The normal subset $\widetilde{S}$ determines a normal subset $\widetilde{S}_{\simplex^1}$ in $\widetilde{\simplex}_{\simplex^1}.$ It follows from Lemma \ref{lem:prep normal unique} that this may be normally isotoped so that $p_{\simplex^1}(\widetilde{S}_{\simplex^1})$ is normal in $B(\simplex^1).$ For each 1--simplex $\widetilde{\simplex}^1$ in $\widetilde{\simplex}$ there exists a unique marked edge $\widetilde{\simplex}^1_m$ in some $\widetilde{\simplex}_{h,k}.$ Define $i$ on $\widetilde{\simplex}^1$ to be an normal isotopy taking the normal corners in $\widetilde{S} \cap {\widetilde{\simplex}^1}$ to the corresponding normal corners in $i_{h,k}^{-1}(\widetilde{S}_{\simplex}^1 \cap \widetilde{\simplex}^1_m).$ Since each 1--simplex in $\widetilde{\simplex}$ corresponds to a unique marked edge, this is well defined and extends to the desired normal isotopy.
\end{proof}

\begin{rem}
There are proofs for two special cases of Theorem \ref{thm:normal unique} in the literature. The proofs by Tollefson \cite{to} (for compact manifolds) and Kang \cite{k} (for topologically finite manifolds with torus ends) construct an explicit surface. A different proof for 3--manifolds with torus cusps is sketched by Weeks in the documentation of {\tt SnapPea} \cite{we}. The above proof is inspired by the latter.
\end{rem}


\subsection{Projective solution space}
\label{comb:Projective solution space}

Let $B$ be the coefficient matrix of the Q--matching equations. Considering $B$ as a linear transformation $\RR^{3t}\to\RR^{e},$ the set of all solutions to the $Q$--matching equations is $\ker B.$ This will be denoted by $Q(\tri).$ One often considers the \emph{projective solution space} $PQ(\tri),$ which consists of all elements of $Q(\tri)$ with the property that the sum of the coordinates equals one. This is a convex polytope, and its vertices are called \emph{vertex solutions}. Since the entries in $B$ are integers, the vertices of $PQ(\tri)$ are rational solutions, and hence the subset of all rational points is dense in $PQ(\tri).$ Such a polytope is termed \emph{rational}. The vector in $\RR^{3t}$ with each coordinate equal to one is contained in $Q(\tri).$ Thus, $\dim_\RR PQ(\pseudo; \tri) = \dim_\RR Q(\pseudo; \tri) -1.$

The subset of all admissible solutions in $PQ(\tri)$ is denoted by $\N (\tri).$ If $N \in \N (\tri)$ is rational, then there is a countable family of spun-normal surfaces $S_i$ without vertex linking components such that $N(S_i) = \alpha_i N$ for some $\alpha_i\in\RR.$ A spun-normal surface is a \emph{minimal representative} for $N,$ if the corresponding scaling factor is minimal. Since the admissible solutions can be found by setting in turn two coordinates from each tetrahedron equal to zero, it follows that $\N (\tri)$ is a finite union of convex rational polytopes, each of which is contained in a $(t-1)$--dimensional rational polytope. Addition of solutions within a convex cell of $\N (\tri)$ corresponds to the geometric sum operation.


\subsection{Spun-normal branched immersions}
\label{subsec:bins}

\begin{definition}[Spun-normal branched immersion]\label{defn:Spun-normal branched immersion}
Let $\pseudo$ be a pseudo-manifold (possibly with boundary). Let $S$ be obtained from a compact (not necessarily connected or closed) surface by removing some boundary components. Suppose that $S$ has a cell structure with $2$--cells which are either triangles and quadrilaterals, and that $f\co S\rightarrow \pseudo$ is a piecewise linear map such that the interior of every $2$--cell of $S$ is mapped homeomorphically onto the interior of a normal disc in $\pseudo$ satisfying the following extra conditions.
\begin{enumerate}
\item $f$ is normally isotopic to an immersion in the complement of the $0$--skeleton of $S;$
\item $f$ is transverse to the 2--singlices in $\pseudo;$
\item the set of accumulation points of $f(S) \cap \pseudo^{(1)}$ is contained in $\pseudo ^{(0)};$
\item if $\{ x_i\} \subset f(S)$ has accumulation point $x\in \pseudo,$ then $x \in f(S)$ or $x \in \pseudo ^{(0)}.$
\end{enumerate}
Then $S$ contains finitely many quadrilaterals and $f$ is proper. The map $f$ is called a \emph{spun-normal branched immersion}, and if $f$ is in general position, then $f(S)$ is termed a \emph{branched immersed spun-normal surface} in $\pseudo.$ Let $x(f) \in \mathbb{Z}^{3t}$ be the point defined as follows: the coefficient for each normal quadrilateral type is the number of 2--cells in $S$ which map to discs of this type.
\end{definition}

\begin{pro}\label{pro:normal branched immersions}
Let $\pseudo$ be a 3--dimensional pseudo-manifold with triangulation $\tri.$ If $f\co S\rightarrow \pseudo$ is a spun-normal branched immersion, then $x(f) \in Q(\tri).$ Every non-zero point in $Q(\tri)$ with non-negative integral coordinates is represented by a (not necessarily unique) spun-normal branched immersion with finitely many branch points.
\end{pro} 

It follows from the definition that $x(f) \in Q(\tri);$ the remainder of the proposition is proved in Subsection \ref{subsec:bins proof}.


\subsection{Traditional normal surface theory}

To the $3t$ quadrilateral coordinates $(q_{1},...,q_{t}'')$ adjoin $4t$ further coordinates, one for each normal triangle type. If $S$ is a closed normal surface in $\pseudo,$ let $N_{\Delta}(S) \in \RR^{7t}$ be its \emph{normal coordinate}. The normal coordinate of a closed normal surface satisfies three \emph{compatibility equations} for each 2--singlex not contained in $\partial \pseudo,$ arising from the fact that the total numbers of normal arcs ``on either side'' have to match up. Denote by $C(\tri)$ the space of all solutions to the compatibility equations.

\begin{figure}[h]
\psfrag{v}{{\small $v$}}
\psfrag{g}{{\small $\gamma$}}
\psfrag{t}{{\small $t_{0}$}}
\psfrag{q}{{\small $q_{0}$}}
\psfrag{p}{{\small $q_{1}$}}
\psfrag{s}{{\small $t_{1}$}}
\psfrag{+}{{\small $+$}}
\psfrag{=}{{\small $=$}}
\psfrag{qk}{{\small $q_k''$}}
\psfrag{q1}{{\small $q_1'$}}
\psfrag{qa}{{\small $q_{1}''$}}
\psfrag{q2}{{\small $q_2'$}}
\psfrag{qb}{{\small $q_{2}''$}}
\psfrag{q3}{{\small $q_3'$}}
\psfrag{qa}{{\small $q_{1}''$}}
\psfrag{q3-qb=t2-t3}{{\small $q_3'-q_2''=t_2-t_3$}}
\psfrag{q2-qa=t1-t2}{{\small $q_2'-q_1''=t_1-t_2$}}
\psfrag{q1-qk=tk-t1}{{\small $q_1'-q_k''=t_k-t_1$}}
\begin{center}
\subfigure[Compatibility equations]{
    \includegraphics[width=5.5cm]{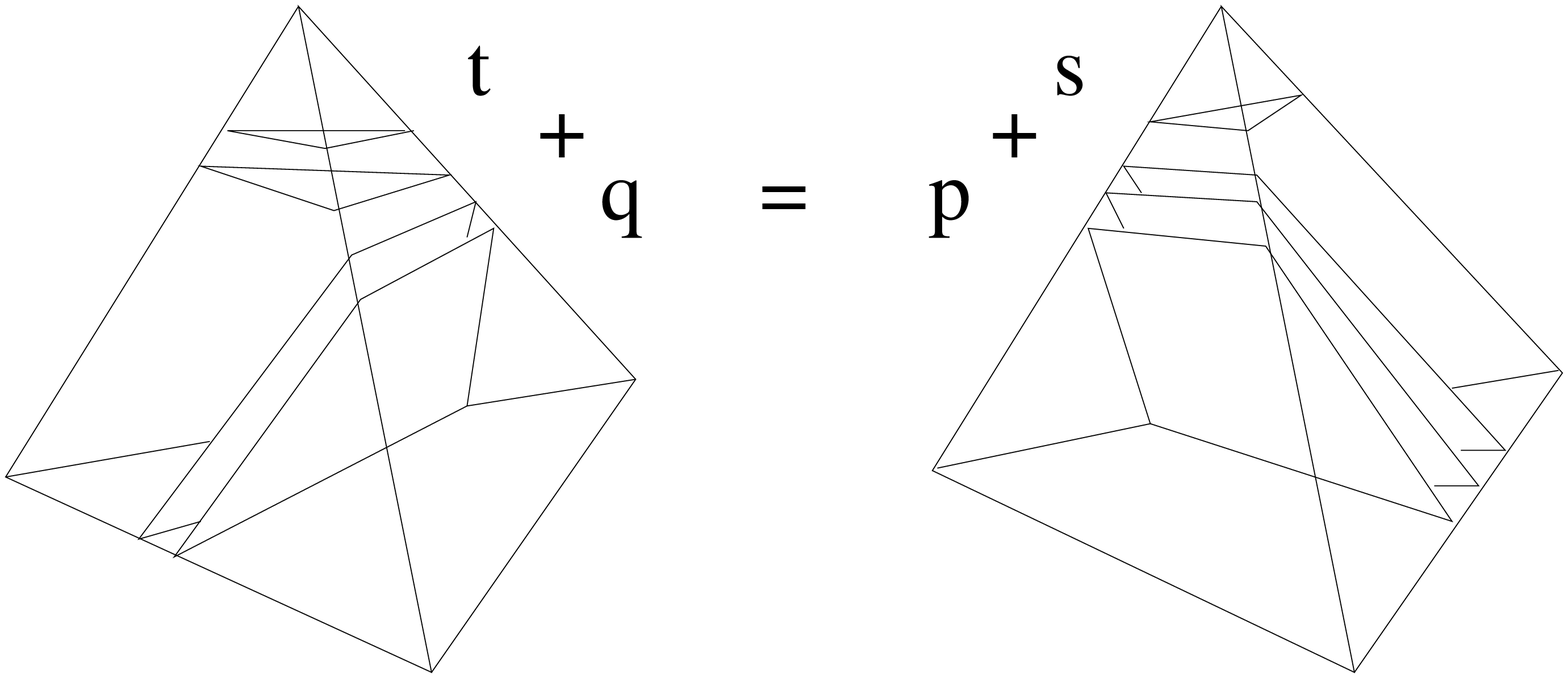}}
\qquad
\subfigure[Compatibility and $Q$--matching]{
    \includegraphics[width=5.5cm]{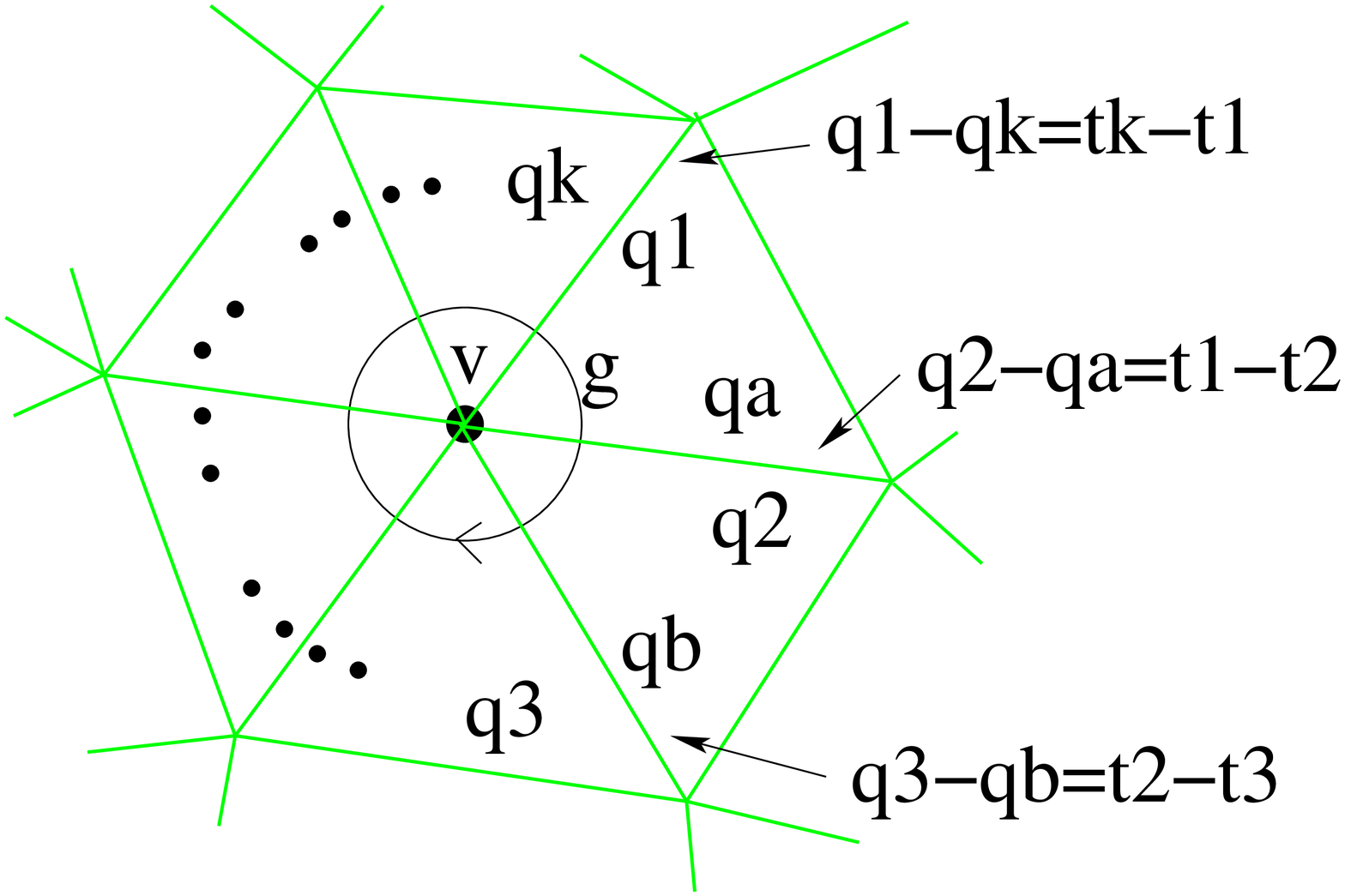}}
\end{center}
    \caption{Normal surface coordinates}
    \label{fig:trad_coods}
\end{figure}

Since the compatibility equations can be viewed as conditions associated to normal arc types, there is precisely one such equation for each edge in the induced triangulation of $\partial \pseudo^c$ not contained in its boundary. Let $\normalc \in B_\v^{(0)}$ be contained in the interior of $B_\v.$ Then each element of $C(\tri)$ satisfies one compatibility equation for each normal arc in $B_\v$ which ends in $\normalc$ (see Figure \ref{fig:trad_coods}(b)). Pick one of the arcs, and write the compatibility equation in the form $q_{i+1}'-q_{i}''=t_{i}-t_{i+1}.$ Then, proceeding to the next arc, one obtains $q_{i+2}'-q_{i+1}''=t_{i+1}-t_{i+2},$ etc. Summing all these equations gives on the left hand side (up to sign) the
$Q$--matching equation of the 1--singlex in $\pseudo$ containing $\normalc,$ and the right
hand side equals zero.  Hence there is a well--defined linear map 
$$pr\co C(\tri) \to Q(\tri)$$ 
defined by projection onto the quadrilateral coordinates.

A canonical basis for $C(\tri)$ consisting of $t$ \emph{tetrahedral solutions} and $e$ \emph{edge solutions} is given in \cite{kr}. Let $\simplex^3$ be a 3-singlex in $\pseudo,$ and let $\{ \normalq_i, \normalt_j | i=1,...,3, j=1,...,4\}$ be the set of normal discs types in $\simplex^3.$ The \emph{tetrahedral solution} associated to $\simplex^3$ is:
\begin{equation*}
W_{\simplex^3} = \normalt_1+  \normalt_2 +  \normalt_3 +  \normalt_4 - \normalq_1-\normalq_2-\normalq_3.
\end{equation*}
For each edge $\widetilde{\simplex}^1$ contained in a 3-simplex $\widetilde{\simplex}^3$ in $\widetilde{\simplex}$ there is a unique normal quadrilateral type $\normalq(\widetilde{\simplex}^1)$ in $\widetilde{\simplex}^3$ disjoint from it, and there are two normal triangle types $\normalt_1(\widetilde{\simplex}^1)$ and $\normalt_2(\widetilde{\simplex}^1)$ meeting it. The \emph{edge solution} associated to a 1--singlex $\simplex^1$ in $\pseudo$ is:
\begin{equation*}
W_{\simplex^1} = p \Big( \sum_{\widetilde{\simplex}^1\subset p^{-1}(\simplex^1)} 
\normalt_1(\widetilde{\simplex}^1)+\normalt_2(\widetilde{\simplex}^1)-\normalq(\widetilde{\simplex}^1) \Big).
\end{equation*}

Using signed intersection numbers with edges, it follows from work in \cite{kr} that these solutions form a basis for $C(\tri)$ (even though this is not stated in \cite{kr} in this generality):

\begin{pro}[Kang--Rubinstein]\label{pro: basis NS}
Let $\pseudo$ be a pseudo-manifold with triangulation $\tri.$ The set of all tetrahedral and edge solutions is a basis for $C(\tri)$ as a vector space over $\RR.$ In particular, $C(\tri)$ has dimension $t+e.$
\end{pro}


\subsection{Convention for orientable pseudo-manifolds}
\label{comb:conventions}

If $\pseudo$ is oriented, it is possible to fix a convention for the quadrilateral labels which determines the $Q$--matching equations without reference to the abstract neighbourhoods as follows. The orientation of $\pseudo$ is pulled back to $\widetilde{\simplex}.$ Assign the parameter $z_i$ to any edge of $\tilde{\simplex^3}_i,$ and label the remaining edges according to the convention in Section \ref{subsec: convention for tets} in the unique resulting way with $z_i,z'_i, z''_i,$ and denote the quadrilateral types in $\tilde{\simplex^3}_i$ by $q_i,q'_i,q''_i.$  There is a 1--1 correspondence between quadrilateral types in $\tilde{\simplex^3}_i$ and in $\simplex^3_i=p(\tilde{\simplex^3}_i),$ and the quadrilateral types in $\simplex^3_i$ are denoted by the same symbols.

The orientation of $B(\simplex^1)$ can be chosen such that the quotient map $p_{\simplex^1}\co B(\simplex^1) \to \pseudo$ is orientation preserving. It can be deduced from Figure \ref{fig:slopes} that if a 1--simplex in $p^{-1}(\simplex^1)$ has parameter $z_i,$ then the normal quadrilaterals of type $q'_i$ have positive slope on $\partial B(\simplex^1),$ and the normal quadrilaterals of type $q''_i$ have negative slope. Label the 1--singlices in $\pseudo$ by $\simplex^1_1,...,\simplex^1_e.$ It follows that if $a_{ij}$ is the number of pre-images of $\simplex^1_j$ with label $z_i,$ then the contribution to the Q--matching equation is $a_{ij} (q'_i - q''_i).$ Defining $a'_{ij}$ and $a''_{ij}$ accordingly, one obtains the \emph{$Q$--matching equation of $\simplex^1_j$}:
\begin{equation}\label{eq:Q-matching}
    0 = \sum_{i=1}^t (a''_{ij} - a'_{ij}) q_i + (a_{ij} - a''_{ij}) q'_i + (a'_{ij} - a_{ij}) q''_i,
\end{equation}
and the coefficient matrix of the system of Q--matching equations is given by:
\begin{equation}
    B = \begin{pmatrix}
           a''_{11}-a'_{11} & a_{11}-a''_{11} & a'_{11}-a_{11} & 
{\ldots}
& a'_{t1}-a_{t1} \\
\vdots         &         &        &          & \vdots   \\
a''_{1e}-a'_{1e} & {\ldots}        &        &          &
a'_{te}-a_{te} \\
         \end{pmatrix}.
\end{equation}


\section{Boundary matters}
\label{Boundary stuff}

It was shown in the first section that a spun-normal surface is uniquely determined by the quadrilateral discs in its cell decomposition, and that it has well-defined boundary curves. This section introduces an algebraic version of the boundary curve map which is defined on the whole solution space to the $Q$--matching equations. The discussion is restricted to closed pseudo-manifolds for pragmatic reasons: more notation and words are required if $\pseudo$ is not closed.


\subsection{The boundary curve map revisited}
\label{Transverse orientations}

Let $\normalt$ be a normal triangle in $B_\v$ and $\simplex^3$ be the 3--singlex which contains $\normalt.$ Let $\normala$ be a normal arc in $\partial \normalt,$ and $\simplex^2$ be the face of $\simplex^3$ containing $\normala.$ There is a unique quadrilateral type $q$ in $\simplex^3$ such that $q$ and $\normalt$ have the same arc type on $\simplex^2;$ see Figure \ref{fig:Delta labels}. Let $q$ be the \emph{$Q$--modulus of $\normala$ (with respect to $\normalt$)}, and give $\normala$ a transverse orientation (with respect to $\normalt$) by attaching a little arrow pointing into the interior of $\normalt.$ This construction is dual to the labelling of the vertices. The orientation conventions of Subsections \ref{subsect:double cover} and \ref{subsect:boundary curve map} will be used throughout this section. Lift the labelling to $\widetilde{B}_\v$ if $B_\v$ is non-orientable; otherwise write $\widetilde{B}_\v=B_\v.$

Let $\gamma$ be an oriented path in $\widetilde{B}_\v$ which is disjoint from the 0--skeleton, has endpoints in the 1--skeleton, and whose interior meets the 1--skeleton transversely. To $\gamma$ one can associate a linear functional $\nu(\gamma)$ in the quadrilateral types by taking the positive $Q$--modulus of an edge if it crosses with the transverse orientation, and by taking the negative $Q$--modulus if it crosses against it (where each edge in $\widetilde{B}_\v$ is counted twice --- using the two adjacent triangles). Evaluating $\nu(\gamma)$ at a solution $N$ to the $Q$-matching equations, gives a real number $\nu_{N}(\gamma).$ 

\begin{figure}[t]
\psfrag{z}{{\small $z$}}
\psfrag{z'}{{\small $z'$}}
\psfrag{z"}{{\small $z''$}}
\psfrag{q}{{\small $q$}}
\psfrag{q'}{{\small $q'$}}
\psfrag{q"}{{\small $q''$}}
\psfrag{q1}{{\small $q_1$}}
\psfrag{q2}{{\small $q_2$}}
\psfrag{q3}{{\small $q_3$}}
\psfrag{q4}{{\small $q_4$}}
\psfrag{qk}{{\small $q_k$}}
\psfrag{qk-1}{{\small $q_{k-1}$}}
\psfrag{v}{{\small $v$}}
\psfrag{g}{{\small $\gamma$}}
\begin{center}
  \subfigure[Quadrilateral slopes and transverse orientation]{
    \includegraphics[width=6.3cm, height=2.8cm]{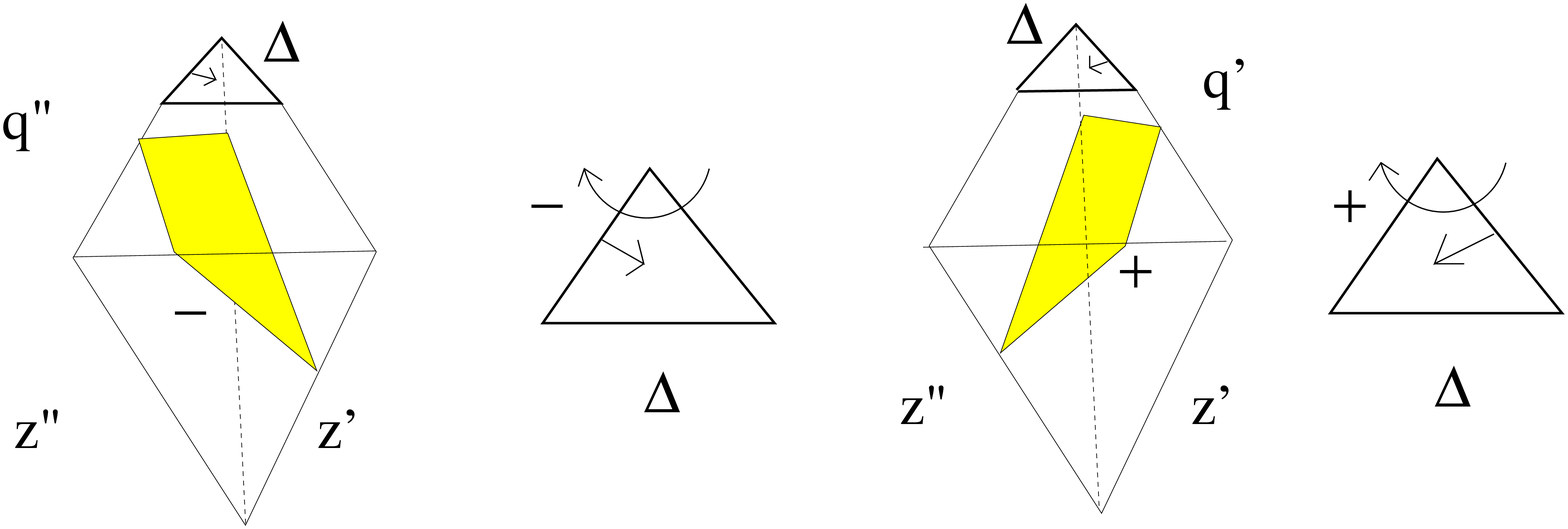}}
  \quad 
  \subfigure[$0=\nu(\gamma) = \sum_{i=1}^{k} (-1)^i q_i$ is the
  $Q$--matching equation.]{
    \includegraphics[width=5.3cm]{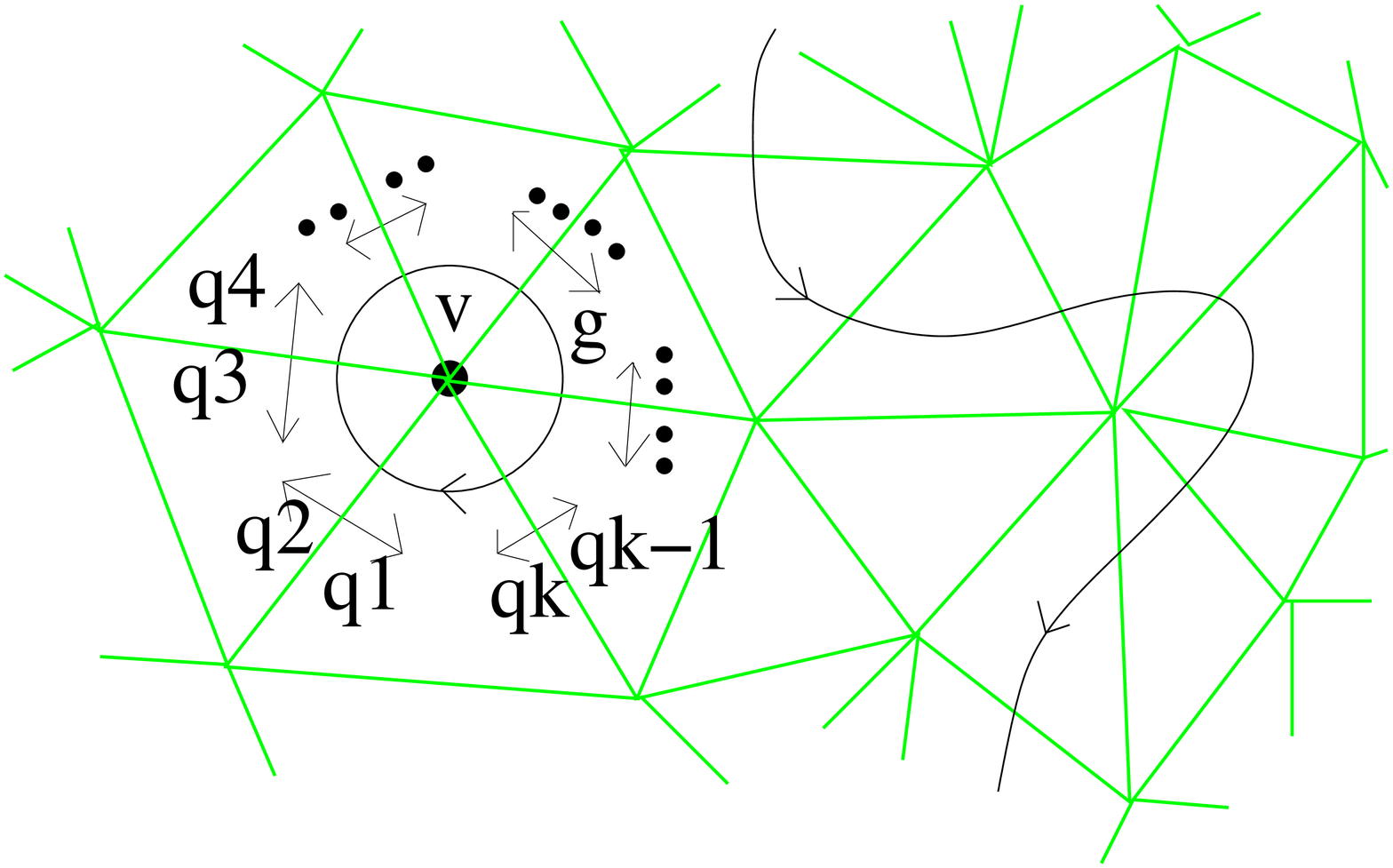}}
\end{center}
    \caption{$Q$--moduli of edges in $B_\v$ and slopes of 
      normal quadrilaterals}
    \label{fig:Q_duality}
\end{figure}

Let $\normalc$ be a vertex in the triangulation of $B_\v$ contained in the 1--singlex $\simplex^1$ in $\pseudo,$ and consider a pre-image $\widetilde{\normalc}$ of $\normalc$ in $\widetilde{B}_\v.$ Then the the linear functional $\nu(\gamma)$ associated to a small circle $\gamma$ with clockwise orientation around $\widetilde{\normalc}$ gives the $Q$--matching equation of $\simplex^1$ by setting $\nu(\gamma)=0.$ This can be deduced from Figure \ref{fig:Q_duality}(b).

\begin{lem} \label{boundary:homo}
Let $\gamma$ be an oriented, closed loop in $\widetilde{B}_\v$ and $N \in Q(\tri).$ The number $\nu_{N}(\gamma ) \in \RR$ depends only on the homotopy class of $\gamma$ and defines a homomorphism $\nu_{N, \v} \co \pi_1(\widetilde{B}_\v) \to (\RR , +).$
\end{lem}

\begin{proof}
Since $\nu_{N}(\gamma)=0$ if $\gamma$ is a small circle about a vertex of $\widetilde{B}_\v,$ it follows that $\nu_{N}$ is well-defined for homotopy classes of loops which intersect the 1--skeleton transversely away from the 0--skeleton.
\end{proof}

\begin{figure}[t]
\psfrag{a}{{\small $\l_1$}}
\psfrag{b}{{\small $\m_1$}}
\psfrag{c}{{\small $\l_i$}}
\psfrag{d}{{\small $\m_i$}}
\psfrag{e}{{\small $\l_{k}$}}
\psfrag{f}{{\small $\m_{k}$}}
    \begin{center}
      \includegraphics[width=10cm]{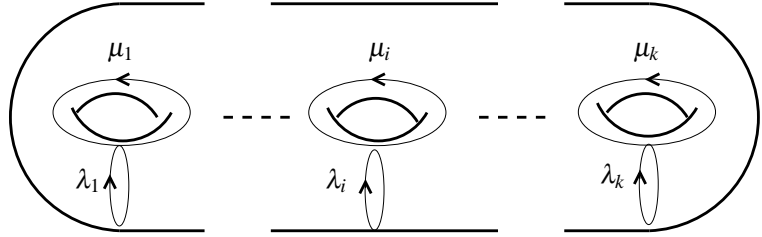}
    \end{center}
\caption{The shown surface $\widetilde{S}$ of genus $k,$ $k\ge 0$ even or odd, is placed in $\RR^3$ such that it is invariant under reflection in any of the three coordinate planes. It follows that the antipodal map $(x,y,z)\to(-x,-y,-z)$ restricts to an orientation reversing homeomorphism with quotient $S$ so that the quotient map $\widetilde{S}\to S$ is a covering of degree two.}
\label{fig:surface gens}
\end{figure}

Assume that $\widetilde{B}_\v$ is a closed, orientable surface of genus $g_\v.$ Corresponding to a basis of $H_1(\widetilde{B}_\v)$ choose the system of $2g_\v$ closed oriented curves shown in Figure~\ref{fig:surface gens}.  Then one obtains a linear map $\nu_\v \co Q(M) \to\RR^{2{g_\v}}$ defined by
\begin{equation}
\nu_\v (N) =   (-\nu_N(\l_1), \nu_N(\m_1),{\ldots} ,-\nu_N(\l_{g_\v}), \nu_N(\m_{g_\v})).
\end{equation}
Let $g=\sum g(\widetilde{B}_\v),$ and let $\nu \co Q(M) \to\RR^{2{g}}$ be the map defined by $\nu = \oplus \nu_\v.$

There is the following canonical isomorphism $\mathbb{Z}^{2g_\v}\cong H_1(\widetilde{B}_\v; \mathbb{Z}).$ Let $\gamma$ be an oriented closed curve on $\widetilde{B}_\v.$ Then $(x_1, y_1, ..., x_{g_\v}, y_{g_\v}) \in \mathbb{Z}^{2g_\v}$ maps to $[\gamma]\in H_1(\widetilde{B}_\v; \mathbb{Z})$ if and only if $\iota(\gamma, \l_i)=-x_i$ and $\iota(\gamma, \m_i)=y_i$ for each $i,$ where $\iota$ denotes the algebraic intersection number. This convention identifies the standard system of curves with the standard basis of $\mathbb{Z}^{2g_\v}.$ This extends to a unique homomorphism 
$$h_\v\co \RR^{2g_\v} \cong H_1(\widetilde{B}_\v; \mathbb{R})\to H_1({B}_\v; \mathbb{R}).$$ 
It follows from the definitions that $h_\v (\mathbb{Z}^{2g_\v} ) \subseteq  H_1({B}_\v; \mathbb{Z}).$ 

\begin{definition}[Boundary curve map]
Let $\partial_\v= h_\v\circ \nu_{\v}$ and $\partial \co Q(\pseudo) \to \oplus_\v H_1(B_\v ; \RR)$ be the map defined by $\partial = \oplus_\v \partial_\v.$
\end{definition}

The map $\nu_\v,$ and hence the boundary curve map, has the following geometric property:

\begin{pro}\label{pro:boundary curve map gives boundary curves}
If $N$ is the normal $Q$--coordinate of the spun-normal surface $S,$ then $h_\v\circ \nu_\v(N)=\partial_\v(S).$
\end{pro}

\begin{proof}
It suffices to assume that $\pseudo$ is oriented. Let $S$ be a spun-normal surface in $\pseudo,$ and denote by $S'$ the subcomplex consisting of all normal quadrilaterals in $S.$ Orient the edges of normal quadrilaterals such that the orientation of the normal arc, the corresponding normal vector $n_{\v}$ pointing towards the dual vertex and the transverse orientation form a positively oriented basis (see Figure \ref{fig:Q_orientation}). It follows that whenever quadrilateral disc $\tilde{\normalq}$ has positive slope with respect to $\widetilde{\simplex}^1,$ then its oriented edges point away from $\tilde{\normalq} \cap \widetilde{\simplex}^1,$ and if it has negative slope then they point towards the intersection point. It follows that at each point of intersection of $S$ with a 1--singlex, the number of oriented edges pointing towards it equals the number pointing away from it. In particular, one obtains an oriented 1--chain $c$ of quadrilateral edges in $\partial S'$ along which subcomplexes made up of triangles meet $S'$ in $S.$

The normal vector $n_{\v}$ and the orientation of $\pseudo$ determine an orientation of each vertex linking surface $B_\v.$ Give each transversely oriented curve in $B_\v \cap S$ the orientation induced from the orientation of $B_\v.$ Then $[B_\v \cap S]=\partial_\v(S)\in H_1(B_\v; \mathbb{Z}).$ It suffices to show that $c$ can be viewed as a union of oriented closed curves $\cup c_i$ with the property that each $c_i$ is homotopic in $S\cup B_\v$ to an oriented simple closed curve $h_i$ in $B_\v$ such that $[B_\v \cap S]=\sum [h_i].$

\begin{figure}[t]
\begin{center}
 \subfigure[Orientation of quadrilateral sides]{
  \includegraphics[width=5cm]{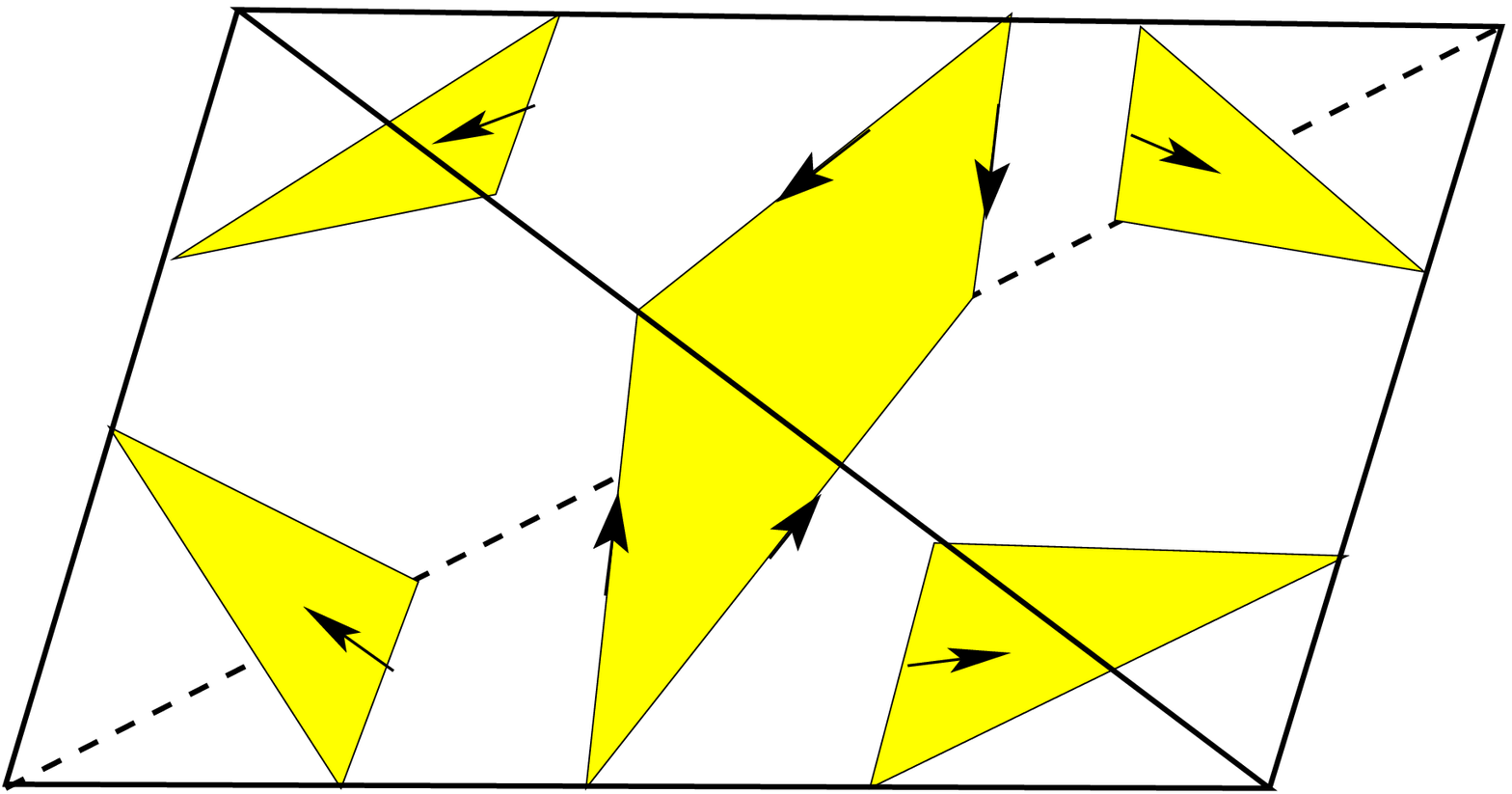}
  }
  \qquad
 \subfigure[Construction and labelling of oriented arcs]{
  \includegraphics[width=6cm]{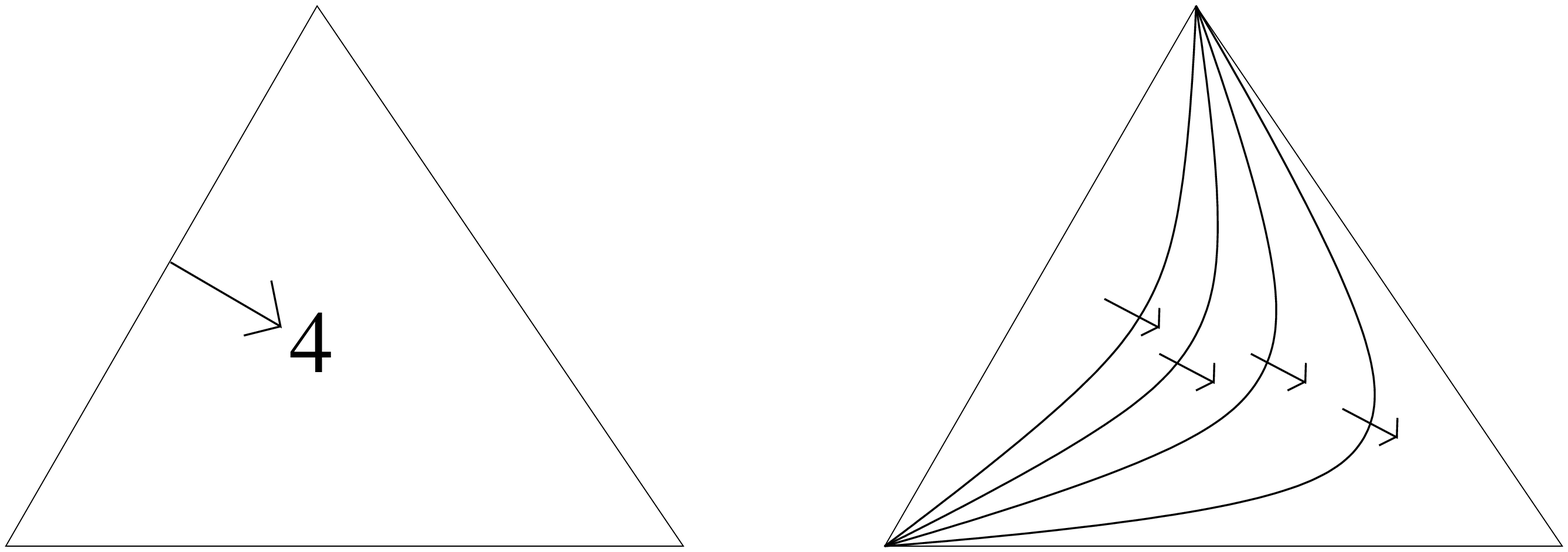}
  }  
\end{center}
  \caption{Oriented 1--chains of quadrilateral edges}
  \label{fig:Q_orientation}\label{fig:Q_orient_arcs}
\end{figure}

The following construction determines $\cup c_i$ and $\cup h_i.$ Given $N=N(S),$ associate a modulus and a transverse orientation to 1--singlices in $B_\v$ as follows. Identify $Q$--moduli with the values given by $N.$ Let $s$ be an edge in $B_\v$ and let $\normalt$ and $\normalt'$ be the triangles in $B_\v$ meeting in $s,$ and $n$ and $n'$ be the associated $Q$--moduli. If $n=n',$ then give $s$ the modulus zero and no transverse orientation. If $n > n',$ then give $s$ the modulus $n-n'$ and the transverse orientation inherited from $\normalt,$ and if $n < n',$ then give $s$ the modulus $n'-n$ and the transverse orientation inherited from $\normalt'.$ The modulus of $s$ determines the number of quadrilateral edges of the same arc type as $s$ which are identified with edges of normal triangles in $S,$ and the transverse orientation determines the 3--singlex containing the corresponding quadrilaterals (see Figure \ref{fig:Q_examples}).

\begin{figure}[t]
\psfrag{v}{{\small $v$}}
\psfrag{a}{{\small $\v$}}
\psfrag{b}{{\small $B(\simplex^1)$}}
\psfrag{e}{{}}
\begin{center}
  \subfigure[Example 1]{
      \includegraphics[width=3.5cm]{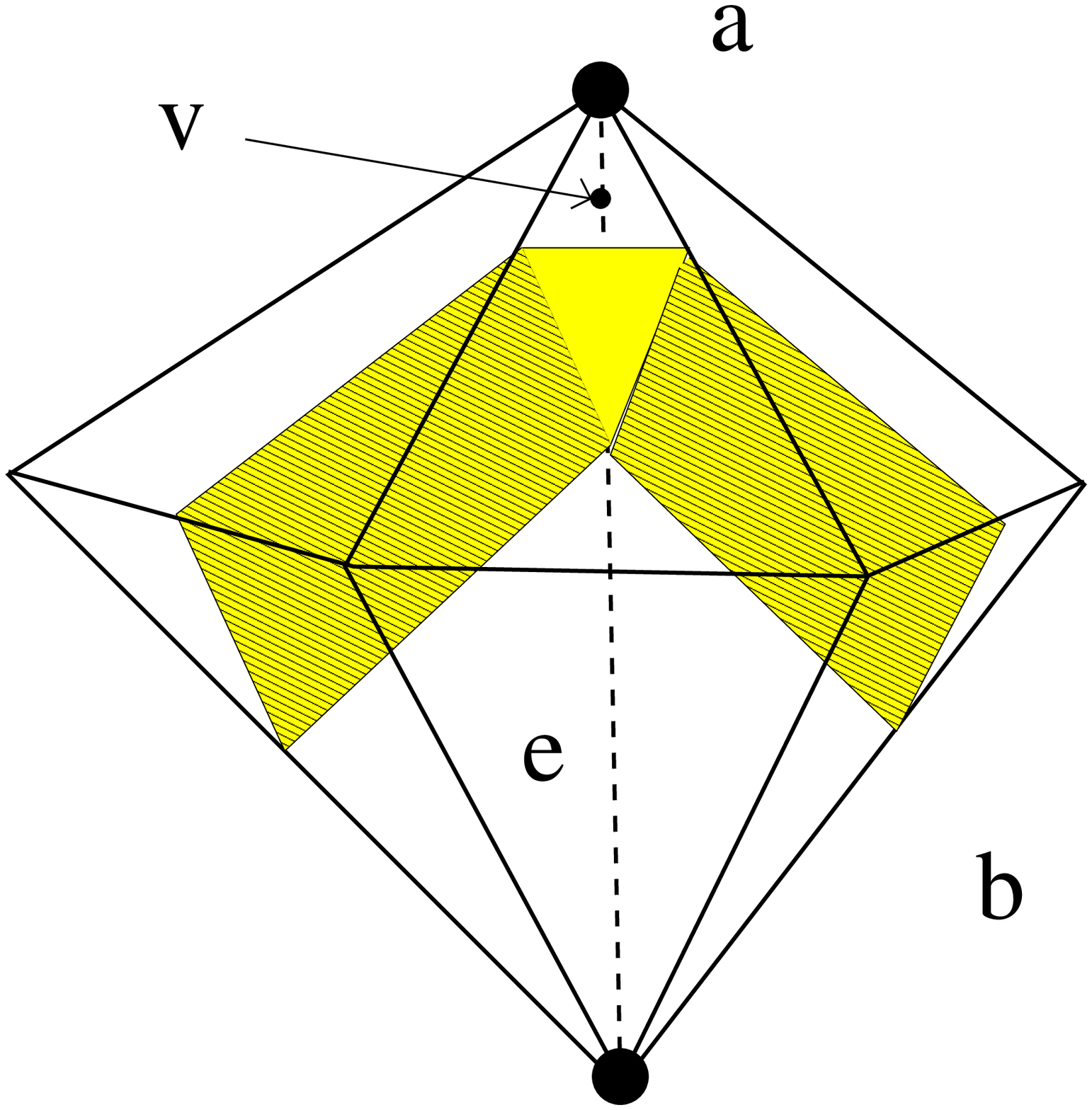}}
      \qquad\qquad
  \subfigure[Example 1]{
     \begin{minipage}[b]{5cm}
        \centering
        \includegraphics[height=1.5cm]{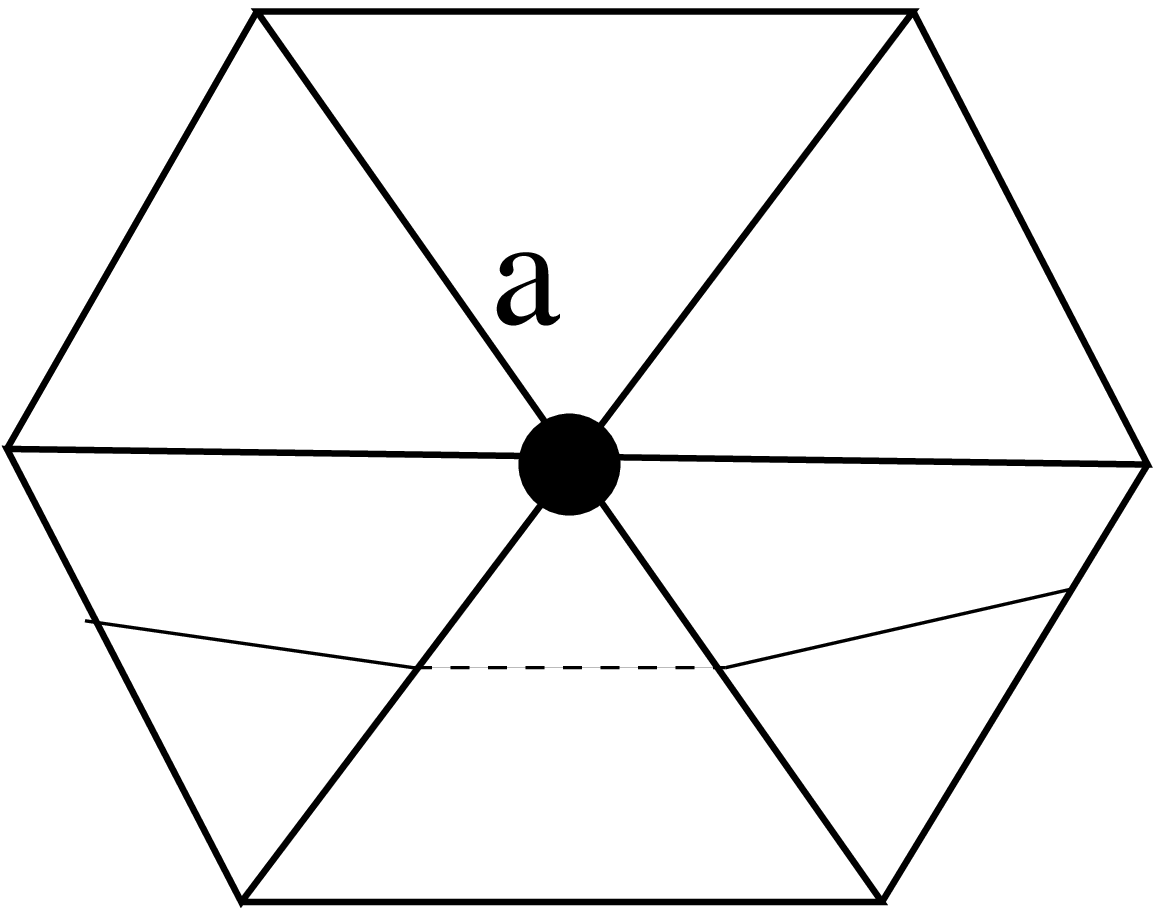}
        \qquad
        \includegraphics[height=1.5cm]{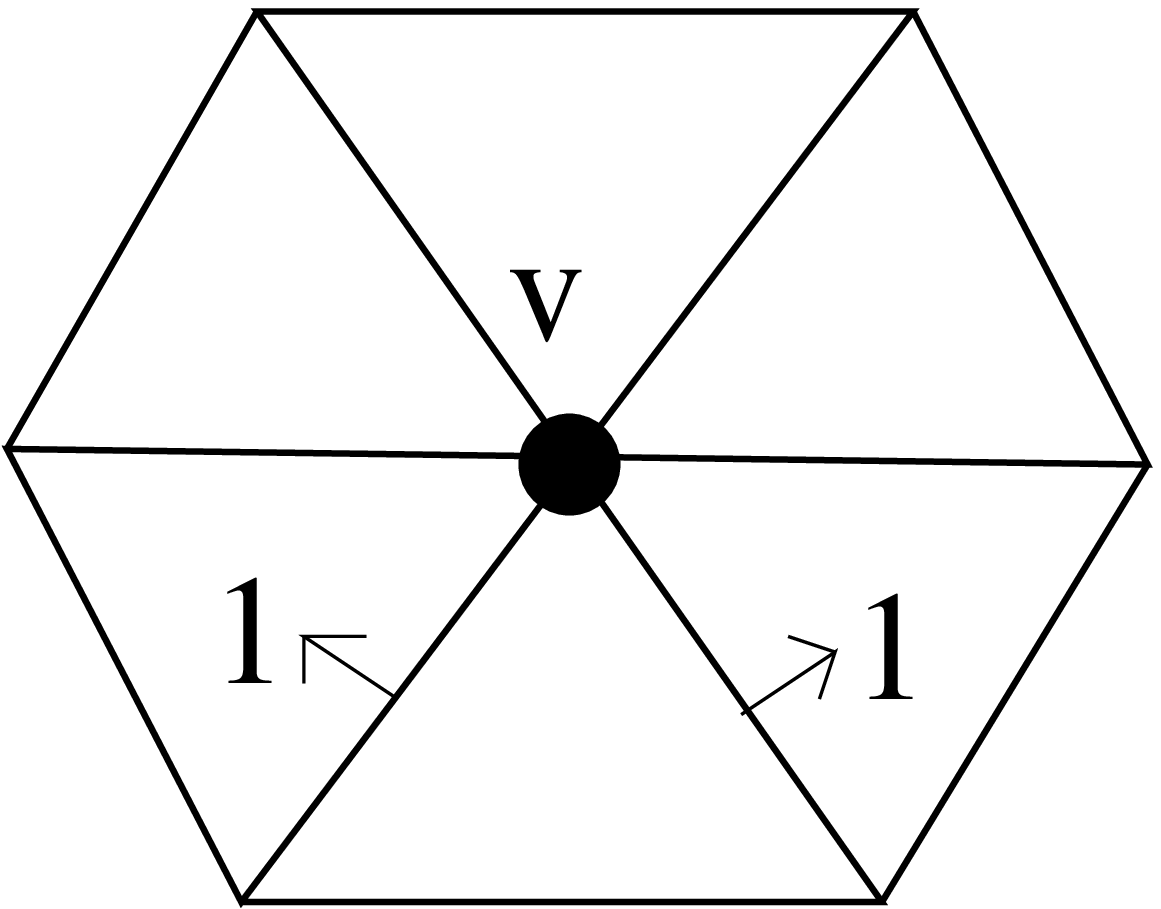}
         \\ \vspace{0.5cm}
        \includegraphics[height=1.5cm]{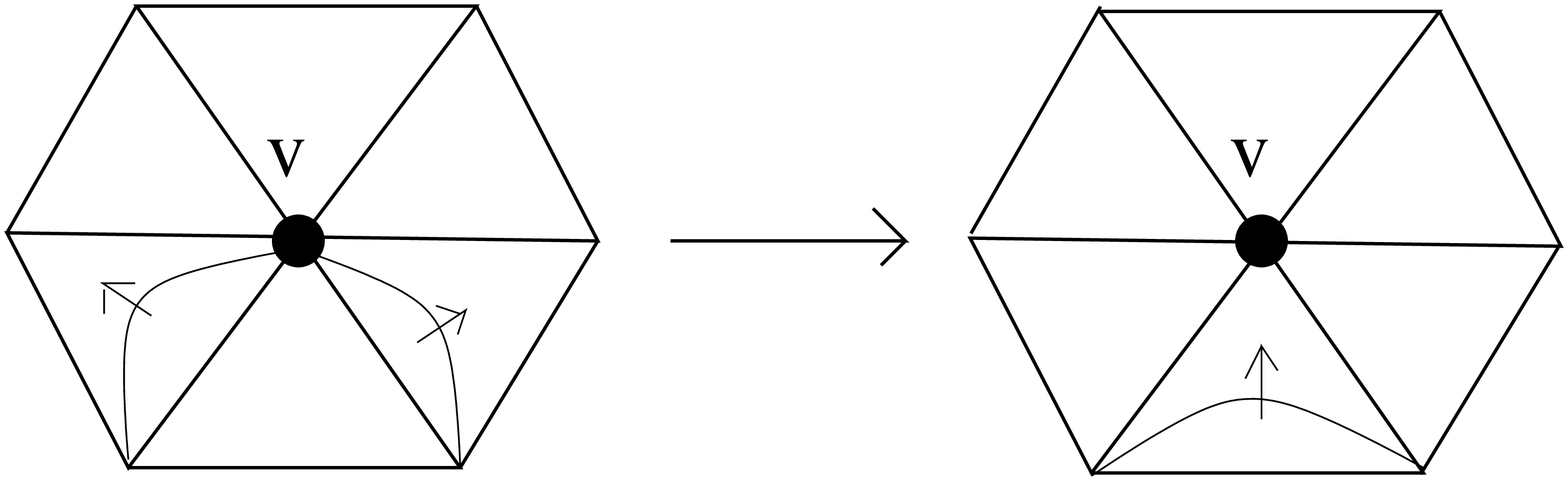}
     \end{minipage}}
   \\
   \subfigure[Example 2]{
              \includegraphics[width=3.5cm]{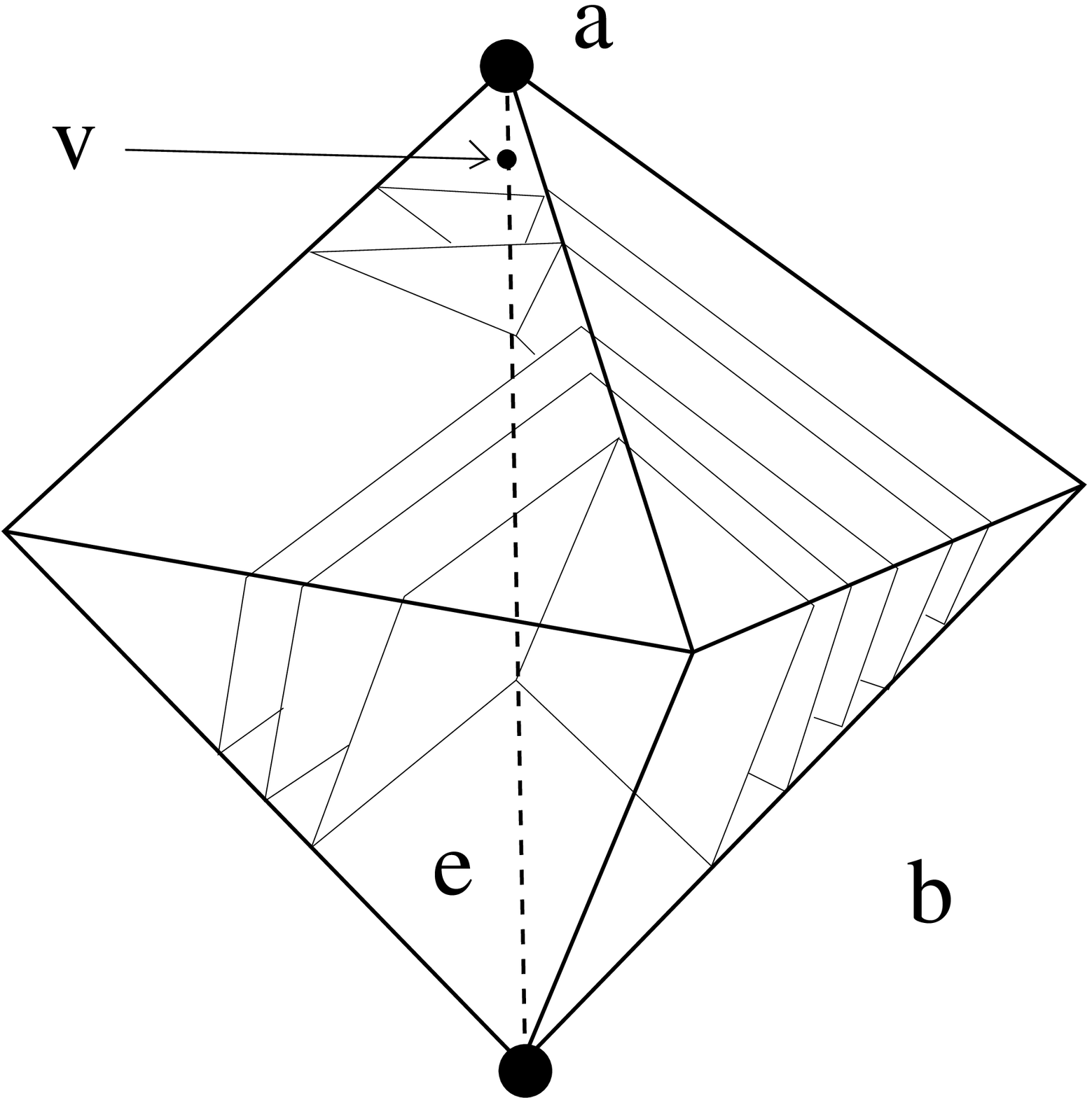}}
      \qquad\qquad
   \subfigure[Example 2]{
     \begin{minipage}[b]{5cm}
        \centering
        \includegraphics[height=1.5cm]{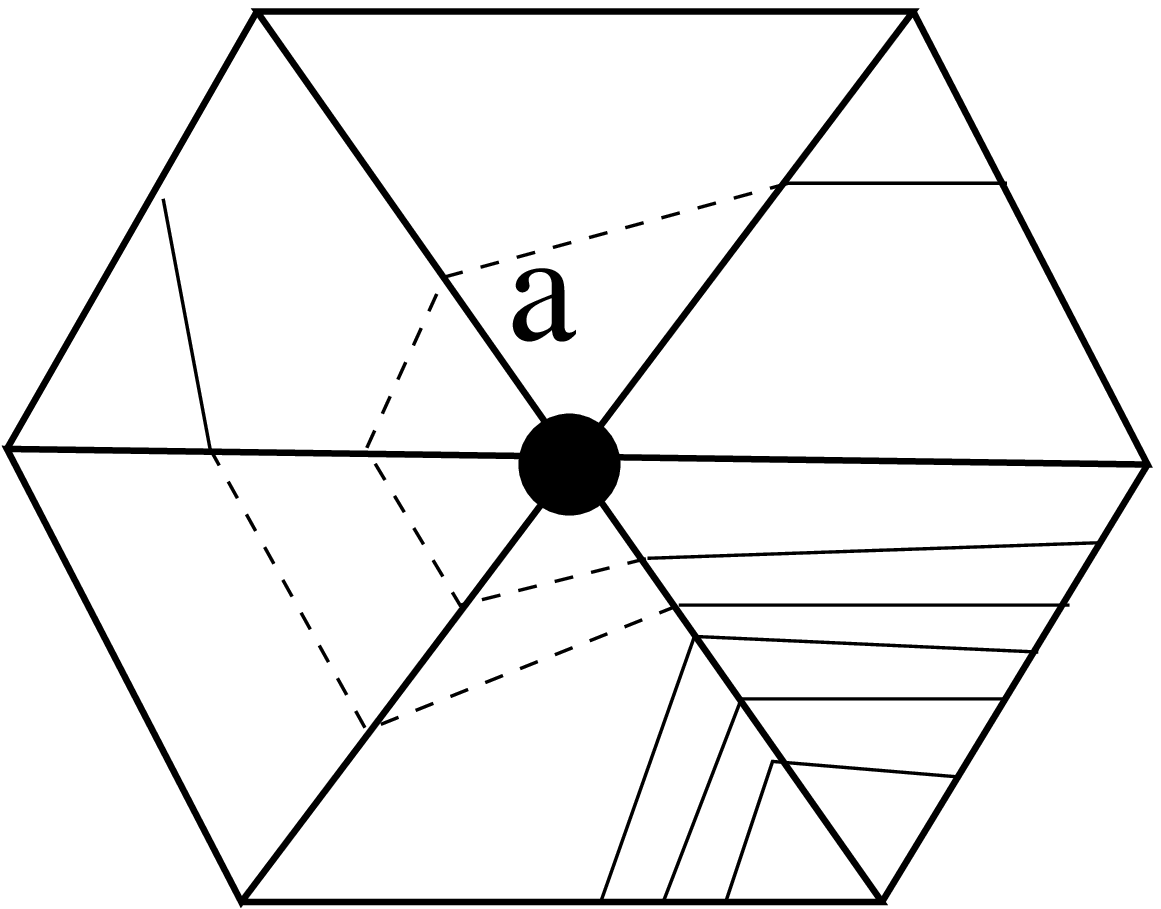}
        \qquad
        \includegraphics[height=1.5cm]{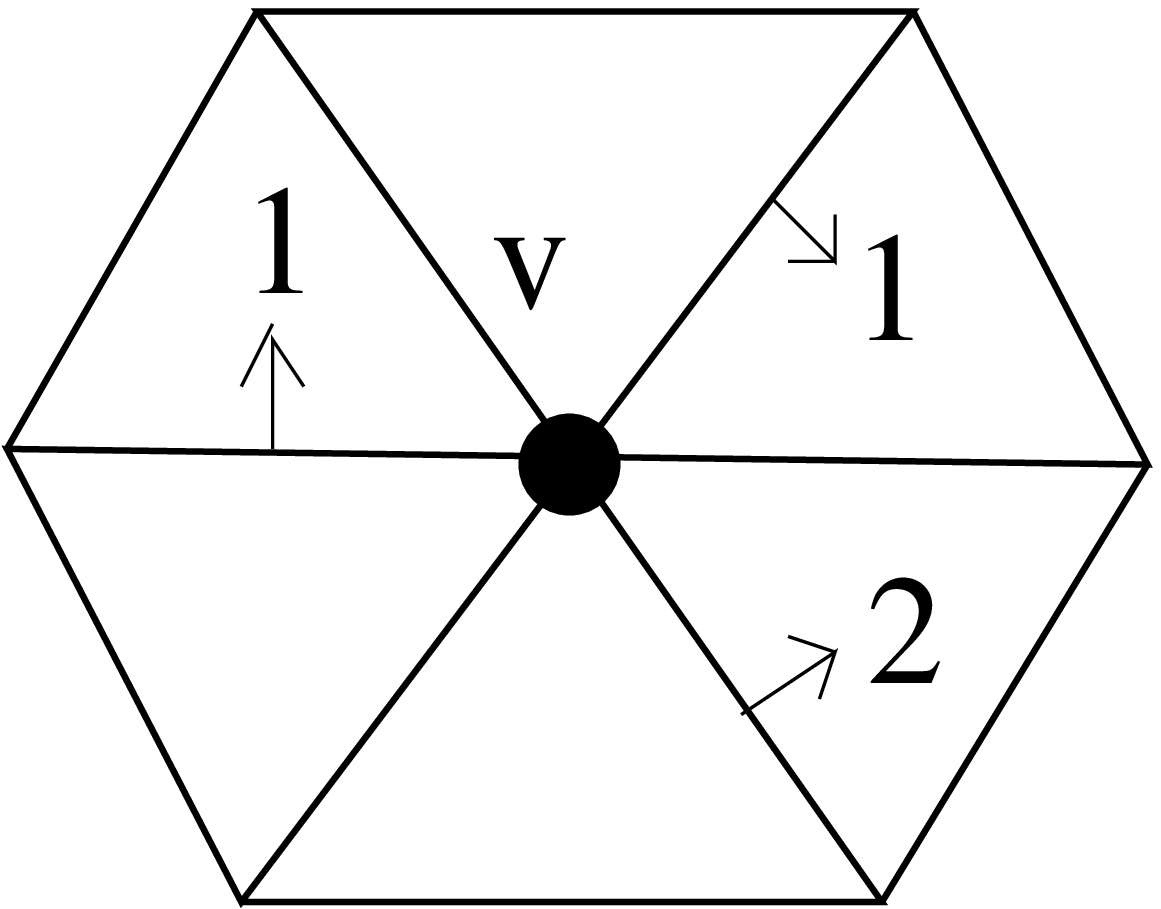}
        \\ \vspace{0.5cm}
        \includegraphics[height=1.5cm]{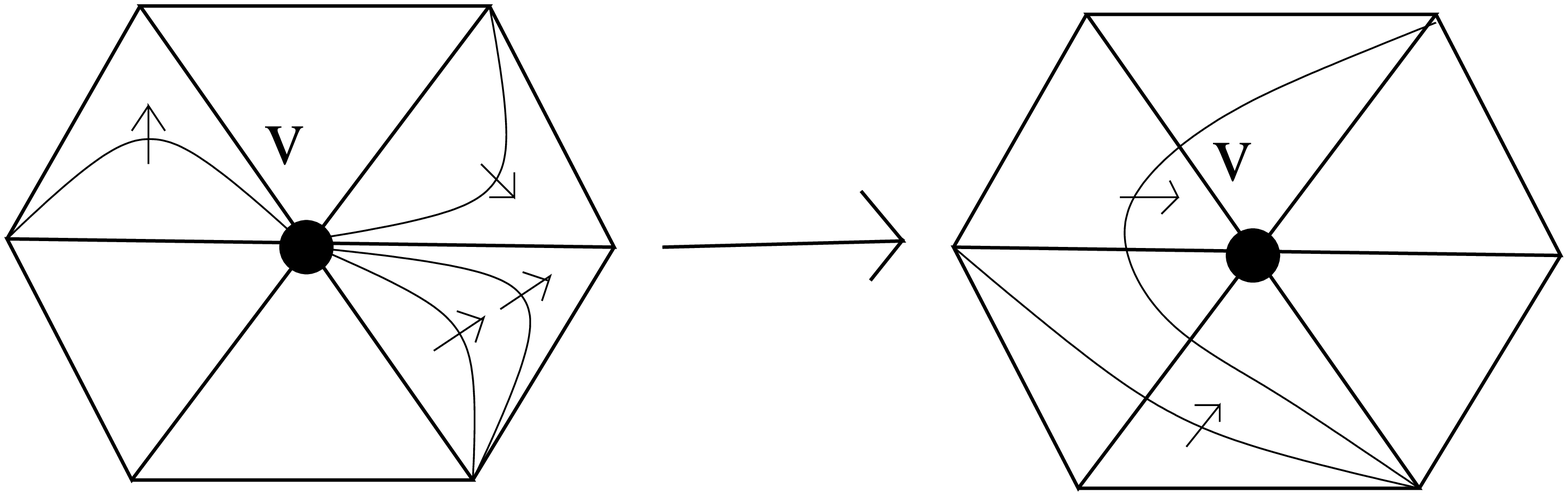}
     \end{minipage}}
\end{center}
\caption[Two examples]{Two examples; figures (b) and (d) show: (Top left)
        The pattern on the upper 
  hemisphere of $B(\simplex^1);$ (Top right) The edge labelling on $B_\v;$ (Bottom)
  Introducing arcs and resolving intersection points at $v$}
     \label{fig:Q_examples}
\end{figure}

Let $c'$ be the union of edges in $B_\v$ with non--zero modulus. A unique collection of transversely oriented simple closed curves on $B_\v$ is derived from $c'$ as follows. If $s$ is an edge in $c'$ with modulus $n$ and transverse orientation pointing into $\normalt,$ place $n$ arcs in $\normalt,$ each of which has endpoints identical to $s,$ any two of which only meet in their endpoints, and give all of them the induced transverse orientation (pointing towards the vertex of $\normalt$ opposite $s;$ see Figure \ref{fig:Q_orient_arcs}). These arcs are in a natural bijective correspondence with edges of quadrilaterals in $S$ meeting edges of triangles. A bijection between endpoints of arcs at a vertex $v$ of $B_\v$ is now defined which corresponds to the identifications of quadrilaterals and triangles in $S.$ Please refer to Figure \ref{fig:Q_examples}(b) and \ref{fig:Q_examples}(d). Let $H_\v$ be the union of all triangles in $B_\v$ containing $v,$ and let $A$ be the set of arcs which have not been paired yet. An
\emph{outermost pair} $\gamma$ and $\delta$ satisfies the following two criteria (this is a variant of constructions in \cite{to, k}):
\begin{enumerate}
\item If the endpoints of $\gamma$ and $\delta$ are identified at $v,$ then the transverse orientations of $\gamma$ and $\delta$ match.
\item The disc cut out from $H_\v$ by $\gamma \cup \delta$ which the transverse orientation of $\delta$ points away from does not contain any arcs in $A.$
\end{enumerate}
If an outermost pair is found, identify their endpoints and isotope the union away from $v$ (opposite to their transverse orientation). Then repeat the above by replacing $A$ by $A - \{ \gamma, \delta\}.$ Thus, after gluing and isotopy, one obtains a unique disjoint collection $\mathcal{C}$ of simple closed curves, each with a transverse orientation. Orient these curves using the above convention for orientation of $Q$--edges. It follows from the construction that:
\begin{enumerate}
\item Each $c_i \in \mathcal{C}$ is homotopic in $S \cup B_\v$ to a unique closed (possibly not simple) curve $c_i''$ of quadrilateral edges in $\partial S'$ (e.g.\thinspace introduce labels for the arcs specifying a 1-- and a 2--simplex of $B_\v$ and a normal quadrilateral). In particular, $c_i''$ inherits a well--defined orientation from $c_i$ which agrees with the orientation of the quadrilateral edges.
\item The element in $H_1(B_\v; \mathbb{Z})$ represented by $\mathcal{C}$ is $\partial_\v(S).$
\end{enumerate}
This completes the proof.
\end{proof}


\subsection{Spun-normal branched immersions}
\label{subsec:bins proof}

\begin{proof}[Proof of Proposition \ref{pro:normal branched immersions}]
Let $N=(x_i)$ be a point in $Q(\tri)$ with non-negative integral coordinates. If $\pseudo$ is not cosed, one may double $\pseudo$ along its boundary and construct a spun-normal branched immersion corresponding to the doubled solution (noting that doubling reverses signs of corners). The restriction of this to either half gives the desired map. Hence assume that $\pseudo$ is closed. At first, a normal branched immersion is constructed which is not in general position. For each quadrilateral type $q_i$ place $x_i$ copies of this type in $\pseudo$ such that all vertices are barycentres of 1--singlices; this is viewed as a map $f_q\co S_q \to \pseudo,$ where $S_q$ is an abstract union of $\sum x_i$ pairwise disjoint quadrilateral discs. Choose each $B_\v$ such that its intersection with the 1--skeleton consists of barycentres.

Using the construction in the proof of Proposition \ref{pro:boundary curve map gives boundary curves}, a finite family, $\mathcal{C}_\v,$ of pairwise disjoint transversely oriented simple closed curves in $B_\v$ can be constructed with the property that each curve in $\mathcal{C}_\v$ is homotopic to a unique curve in the 1--skeleton of $B_\v;$ denote the resulting homotopy taking $\mathcal{C}_\v$ into the 1--skeleton of $B_\v$ by $h_\v\co B_\v \times I \to B_\v.$

Note that the spinning construction can be applied to any finite collection of pairwise disjoint, transversely oriented, non-separating (not necessarily normal) curves on $B_\v.$ After performing the straightening, the result is a union, $S''_\v,$ of normal discs in $\overline{N}_\v$ and normal discs in $\pseudo$ which meet $B_\v$ in normal cells; $S''_\v \cap N_\v$ is properly embedded in $N_\v.$ Applying this spinning construction to $\mathcal{C}_\v$ yields a subset $F_\v$ of $\overline{N}_\v.$ Then $h_\v$ can be extended in a neighbourhood of $B_\v$ in $N_\v$ such that $h_\v(F_\v, 1)$ meets $\pseudo$ in normal triangles and $F_\v \cap N_\v= h_\v(F_\v, 0)\cap N_\v$ is isotopic to $h_\v(F_\v, 1) \cap N_\v.$ Let $S'_\v$ be the union of all connected components of $h_\v(F_\v, 1)$ which are not vertex linking. There is a triangulated (possibly non-compact) surface $S_\v$ and a simplicial isomorphism $f_\v \co S_\v \to S'_\v.$

Lift $S'=f_q(S_q) \cup_\v  f_\v (S_\v)$ to $\widetilde{\simplex}.$ Let $\normala$ be a normal arc in $\pseudo$ with endpoints barycentres of 1--singlices, and denote its pre-images in $\widetilde{\simplex}$ by $\normala_0$ and $\normala_1.$ It follows from the definition of moduli that there are as many normal discs in $p^{-1}(S')$ containing $\normala_0$ as there are normal discs containing $\normala_1.$ Moreover, these numbers are finite. One can therefore define an arbitrary bijection between these two sets which determines a unique bijection between the corresponding boundary arcs in $S_q \cup_\v  S_\v.$ Do this for all such normal arcs and denote the resulting quotient by $S.$ Then $S$ is a triangulated (possibly non-compact) surface $S$ and there is a well-defined normal branched immersion $f\co S \to \pseudo$ obtained from gluing the maps $f_q$ and $f_\v.$ The set of branched points is contained in the pre-image of all barycentres; hence there are finitely many branched points. It follows from the construction that $x(f)=N,$ and the map $f$ can be homotoped into general position if desired. 
\end{proof}


\subsection{Analysis of the linear maps and Proof of Theorem \ref{thm:kr}}

If $B_\v$ is non-orientable, then it may be assumed that the non-trivial deck transformation $\widetilde{\pseudo} \to \widetilde{\pseudo}$ induces an involution $\sigma$ on $\widetilde{B}_\v$ which coincides with the antipodal map in Figure~\ref{fig:surface gens}. Whence $\sigma \m_i = \m_{g_\v+1-i}^{-1}$ and $\sigma \l_i = \l_{g_\v+1-i}.$ Denote the resulting involution on $H_1(\widetilde{B}_\v)$ by $\sigma$ also. It follows that $\nu_N(\m_i) = \nu_N(\sigma \m_i) = - \nu_N(\m_{g_\v+1-i}),$ and $\nu_N(\l_i) = \nu_N(\sigma \l_i) = \nu_N(\l_{g_\v+1-i}).$ This observation yields the following upper bound:
\begin{lem}
$\dim (\im\nu) \le 2\chi(\pseudo)-v_n$
\end{lem}
\begin{proof}
If $B_v$ is orientable, then $2 g_\v =2- \chi(\widetilde{B}_\v) = 2 - \chi(B_\v).$ 

If $B_v$ is non-orientable, then $g_\v = 1 - \chi(B_\v).$ It follows from the description of $\tau$ that the kernel of $\nu_{\cdot, \v}$ has dimension at least $g_\v$ if $g_\v$ is even. If $g_\v=2k+1$ is odd, then $\nu_N(\m_{k+1}) = \nu_N(\sigma \m_{k+1}) = - \nu_N(\m_{k+1})$ implies $\nu_N(\m_{k+1})=0,$ and the kernel of $\nu_{\cdot, \v}$ has again dimension at least $g_\v.$ 

Thus:
\begin{align*}
\dim (\im \nu) \le &\sum_{\text{or'ble}} 2 - \chi(B_\v) + \sum_{\text{non-or'ble}} 1 - \chi(B_\v)\\
& = 2 v_o + v_n - \chi(\partial \pseudo^c) =2\chi(\pseudo)-v_n,
\end{align*}
since $\pseudo$ is closed and hence $\chi(\partial \pseudo^c) = 2 v - 2 \chi(\pseudo).$
\end{proof}

\begin{lem}\label{lem: pr}
  One has $\ker\nu = \im pr.$ Moreover, if $N \in\ker \nu,$ then $pr^{-1}(N)=
  L + span\{N_{\Delta}(B_{1}),...,N_{\Delta}(B_{v})\}$ for some $L \in
  C(\tri).$ 
  
  In particular, the image of $pr\co C(\tri) \to Q(\tri)$ has dimension $t+e-v.$
\end{lem}
\begin{proof}
First assume that  $B_\v$ is orientable. Let $N \in Q(\tri).$ Assume that the induced triangulation $\tri_\v$ of $B_\v$ supports the triangle coordinates $t_{0},...,t_{k}.$ If a value is assigned to $t_{0},$ then the $q$--moduli of the common edges uniquely determine values for the adjacent triangle coordinates such that the matching equations are satisfies. It needs to be shown that these assignments are well defined globally if and only if $N \in\ker \nu.$
  
Indeed, choose a closed, oriented, simple path $\gamma$ in $B_\v,$ disjoint from the 0--skeleton, transverse to the 1--skeleton and meeting 2--cells in at most one normal arc. Let $\normalt$ be a 2--cell of $B_\v$ which $\gamma$ passes through, and assume that $t_{0}$ is the associated triangle coordinate. If $t_0$ is given any value, then the adjacent coordinate $t_{i}$ of the next triangle that $\gamma$ passes through must take the value $t_{i}= (q_{0}-q_{i})+t_{0},$ where $-(q_{0}-q_{i})$ is the oriented sum of the q--moduli associated to the edge crossed by $\gamma.$ In this way, uniquely determined values are given to the coodinates of triangles traversed by $\gamma,$ and upon returning to the initial triangle, one has the equation $t_{0}= -\nu_{N}(\gamma) + t_{0},$ which is satisfied if and only if $\nu_{N}(\gamma)=0.$ Similarly, values can be assigned to the remaining triangle coordinates supported by $\tri_\v.$ This procedure is well defined, if and only if $\nu_{N,\v}\co \pi_1(B_\v)\to (\RR, +)$ is the trivial homomorphism. Changing the initial value assigned to $t_{0}$ changes the triangle coordinates by a multiple of $N_{\Delta}(B_\v).$

Now assume that $B_\v$ is non-orientable. Then the above can be done for the lift of the triangulation of $\pseudo$ to the orientable double cover $\widetilde{\pseudo}.$ Let $\tilde{t}_0$ and $\tilde{t}_1$ be normal triangles on $\widetilde{B}_\v$ which are exchanged by the non-trivial deck transformation $\sigma.$ Let $\gamma$ be a closed normal path in $\widetilde{B}_\v$ passing through $\tilde{t}_0$ and $\tilde{t}_1$ which is invariant under $\sigma.$ Then $\gamma$ can be written as the union of arcs $a$ and $\sigma(a)$ with $\nu_N(a)=\nu_N(\sigma(a)).$ Then $0 = \nu_N(\gamma) = \nu_N(a) + \nu_N(\sigma(a)) = 2 \nu_N(a).$ Whence the above procedure assigns the same values to $\tilde{t}_0$ and $\tilde{t}_1.$ It follows that the normal triangles constructed in $\widetilde{\pseudo}$ can be isotoped to be invariant under $\sigma.$
\end{proof}

\begin{lem}\label{lem: rank B}
The matrix $B$ has rank $e-v_o,$ and we have $\dim Q(\tri) =3t-e+v_o= 2t + \chi(\pseudo)-v_n$ and $\dim \im \nu = 2\chi(\pseudo)-v_n.$
\end{lem}
\begin{proof}
It will first be shown that the rank of $B$ is at most $e-v_o.$ The rank of $B$ equals $3t - \dim \ker (B) = e - \dim\ker B^T.$ It therefore suffices to find $v_o$ linearly independent elements in $\ker B^T.$ The dual system of equations described by the transpose has one equation for each quadrilateral type and one variable, $x_{\simplex^1},$ for each 1--singlex $\simplex^1$ in $\pseudo.$ The equation associated to $q$ is $0 = \sum s_{\simplex^1}(q)x_{\simplex^1},$ where the sum is taken over all 1--singlices in $\pseudo$ and $s_{\simplex^1}(q)$ is the total slope of $q$ w.r.t.\thinspace $\simplex^1$ introduced in Section \ref{subsec:Normal Q--coordinates}. 

Assume that $\pseudo$ is orientable. In this case, $q$ has the same slope at opposite normal corners, and opposite signs at normal corners sharing a normal arc. Let $S$ be a vertex linking surface. Then the assignment $x_{\simplex^1}(S) = \#(S\cap \simplex^1)$ defines an element in $\ker B^T.$ Since each triangle disc is uniquely determined by its intersection with the 1--simplices, it follows that the set of $v$ solutions thus obtained from the $v$ vertex linking surfaces is linearly independent, giving rank at most $e-v.$

Hence assume that $\pseudo$ is non-orientable. Considering its lift to the orientable double cover, one finds that if the edge of a quadrilateral disc in $\pseudo$ is normally isotopic into an orientable vertex linking surface, then the signs at the respective corners are opposite. This implies that the above argument can be applied to an orientable component of $\partial \pseudo^c.$ It does not, however, work for non-orientable components. This gives rank at most $e-v_o$ in general.

Thus, $\dim Q(\tri) \ge 3t-e+v_o= 2t + \chi(\pseudo)-v_n.$ It needs to be shown that this is also an upper bound. One has $2t + \chi(\pseudo)-v_n \le \dim Q(\tri) = \dim \ker \nu + \dim \im \nu = \dim \im pr + \dim \im \nu = t+e-v + \dim \im \nu.$ Thus $2\chi(\pseudo)-v_n \ge \dim \im \nu \ge  2\chi(\pseudo)-v_n.$ This forces equality and hence the conclusion.
\end{proof}

\begin{lem}
$\ker \partial = \ker \nu$
\end{lem}

\begin{proof}
It follows from the definition that $\ker \nu \subset \ker \partial.$ Assume that $\partial(N)=0.$ If $B_\v$ is orientable, then $h_\v$ is an isomorphism which implies $\nu_\v(N)=0.$ Hence assume that $B_\v$ is non-orientable. Since $\dim \im \nu = 2\chi(\pseudo)-v_n,$ it suffices to show that whenever $\delta$ represents the homotopy class of an oriented 1--sided simple closed curve in $B_\v$ with the property that $\nu_\v(N)$ determines $\tilde{\delta},$ then $\iota (\gamma, \tilde{\delta})=0$ for each standard generator. Indeed, since $\delta$ is 1--sided, we have $\sigma (\tilde{\delta}) = \tilde{\delta}.$ Thus, $\iota (\tilde{\delta}, \gamma) = \nu_N (\gamma) = \nu_N(\sigma( \gamma)) = \iota (\tilde{\delta}, \sigma(\gamma)) =  \iota (\sigma(\tilde{\delta}), \sigma(\gamma)) = - \iota (\tilde{\delta}, \gamma).$ Whence $\iota (\tilde{\delta}, \gamma)=0$ for each standard generator, which gives $\nu_\v(N)=0.$
\end{proof}

\begin{proof}[Proof of Theorem \ref{thm:kr}]
The first part of the statement follows from the above lemmata in conjunction with the discussion in Subsection \ref{comb:Projective solution space}. Since $\ker \partial = \ker \nu,$ it follows that $\dim \im \nu = \dim \im \partial,$ and hence that $\partial$ is surjective. Since $\nu$ is defined over the integers, its restriction to integer lattice points in $Q(\tri)$ has image in $\mathbb{Z}^{2g},$ whence $\partial$ has image of finite index in $\oplus_\v H_1(B_\v, \mathbb{Z}).$
\end{proof}


\subsection{Intersection numbers and Proof of Theorem \ref{thm:dim of N(M)}}

Assume that $\pseudo$ is oriented. The material in this subsection is based on the observation that one can determine the algebraic intersection number of the oriented boundary curves of two spun-normal surfaces from their normal $Q$--coordinates. The intersection pairing on $H_1({B}_\v; \RR)$ pulls back to a bi-linear, skew-symmetric pairing $\star$ on $\RR^{2g_\v}.$ Taking sums gives a pairing $\star$ on $\RR^{2g}$ with the property that for any $N, L \in Q(\tri),$ we have:
\begin{equation*}
\nu(N)\star \nu(L)=
\sum_{\v \in \pseudo^{(0)}}  \nu_\v (N) \star \nu_\v (L) = \sum_{\v \in \pseudo^{(0)}} \iota \Big( \partial_\v(N) , \partial_\v(L) \Big).
\end{equation*}
Tho obtain a corresponding bi-linear, skew-symmetric form on $Q(\tri),$ let
\begin{equation*}
    C = \begin{pmatrix}
       0 & 1 & -1 \\
       -1 & 0 & 1 \\
       1 & -1 & 0
       \end{pmatrix},
\end{equation*}
and let $C_t$ be the $(3t \times 3t)$ block diagonal matrix with $t$ copies of
$C$ on its diagonal. Then for any $N, L \in Q(\tri),$ define:
\begin{equation}\label{eq:pairing2}
   \langle N, L \rangle= N^{T} C_{n} L.
\end{equation}

\begin{figure}[t]
  \psfrag{+1}{{\small $\pm 1$}} \psfrag{-1}{{\small $\mp 1$}}
\begin{center}
  \includegraphics[width=9cm]{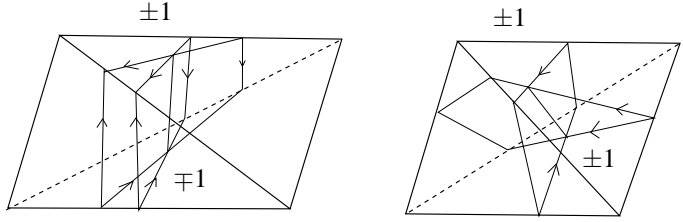}
\end{center}
  \caption{Orientation of quadrilateral sides and intersection number}
  \label{fig:intersection}
\end{figure}

\begin{lem}
Let $\pseudo$ be an oriented, closed pseudo-manifold. For all $N, L \in Q(\tri)$:
\begin{equation*}
\langle N , L\rangle = \nu(N) \star \nu(L).
\end{equation*}
\end{lem}

\begin{proof}
One can compute $\langle N, L \rangle$ by summing intersection numbers obtained for normal quadrilaterals over all 3--singlices as follows (see Figure \ref{fig:intersection}). If $N$ and $L$ have different non--zero $Q$--coordinates in some 3--singlex, the associated normal quadrilaterals are naturally transverse, and the contribution to $\langle N, L \rangle$ is (up to sign) twice the product of the $Q$--coordinates, and the sign is determined by the order of the surfaces. If they have the same non--zero $Q$--coordinate in some tetrahedron, the contribution is zero, independent of the chosen transverse intersection. Given the map $\nu$ and the orientation conventions, this gives precisely the intersection pairing in homology of vertex links.
\end{proof}

\begin{proof}[Proof of Theorem \ref{thm:dim of N(M)}]
If $R$ is a maximal convex polytope in $\N (\tri),$ then it is the intersection of a $(\dim R +1)$--dimensional vector space $R'$ with the unit simplex. The set $\N (\tri)$ is obtained by setting in turn two coordinates from each 3--singlex equal to zero. Thus $R'$ is the intersection of $Q(\tri)$ with a $t$--dimensional subspace of $\RR^{3t},$ and its dimension is at least $(2t + \chi(\pseudo)-v_n)+t-3t=\chi(\pseudo)-v_n.$
  
Assume that $\pseudo$ is oriented. For any $N,$ $L \in R',$ one has $\nu(N)\star \nu(L)=\langle N, L \rangle= 0.$ The image $\nu (R')$ therefore lies in a self--annihilating subspace of $\RR^{2g}.$ Its dimension is thus at most $g= \chi(P),$ giving $\dim R' \le \chi(\pseudo) + \dim (R' \cap \ker \nu).$ Thus $\dim R = \dim R' -1 \le \chi(\pseudo) + \dim (R' \cap \ker \nu)-1 = \chi(\pseudo) + \dim (R \cap \ker \nu).$

If $\pseudo$ is non-orientable, the previous paragraph applies to $\widetilde{\pseudo},$ and we have $\chi(\widetilde{\pseudo}) = 2 \chi(\pseudo) - v_n.$ This completes the proof.
\end{proof}


\section{Examples}
\label{sec:Examples}

A general procedure to determine the boundary slope of a spun-normal surface on a vertex linking torus is given. Then the figure eight knot complement and the Gieseking manifold are discussed. These examples can be found in \cite{kr}.


\subsection{Boundary slopes on tori}

If $B_\v$ is a torus, then the boundary curve map gives a convenient way of determining the boundary curves of a spun-normal surface $S$ from its normal $Q$--coordinate $N=N(S).$ In this case $B_\v$ may be chosen such that it meets $S$ in a finite family of parallel simple closed curves, all with the same transverse orientation, and it meets $N_\v$ in the same number of parallel open annuli. The homology class of such a curve is termed the \emph{slope} of $S$ on $B_\v.$ Since $N_\v \setminus \{\v\}$ is orientable, one can choose generators $\{\l, \m \}$ for $\pi_1(B_\v)$ such that together with a normal vector $n$ pointing into $N_\v,$ $\{\l, \m, n \}$ is a positively oriented basis for $TN_\v.$ This corresponds to a standard meridian--longitude pair of a knot complement. Since $\nu_N$ is a homomorphism, the following holds under the above assumptions.
\begin{enumerate}
\item If $\nu_N(\m) = \nu_N(\l) =0,$ then $S$ is disjoint from $B_\v.$\footnote{This is also shown in the file
    normal\_surface\_construction.c in \cite{we}.}
\item If $\nu_N(\m) \ne 0$ or $\nu_N(\l) \ne 0,$ then let $d>0$ denote the greatest common divisor of the numbers $|\nu_N(\m)|$ and $|\nu_N(\l)|.$ Put $p = - \nu_N(\l)/d$ and $q = \nu_N(\m)/d.$  Then $s = \m^p\l^q$ has $\nu_N(s)=0$ and hence it is a boundary slope of $S.$ Furthermore, $S$ has $d$ boundary curves on $B_\v.$
\end{enumerate}

Now assume that each vertex linking surface is a torus; label the surfaces $B_\v$ by $T_1,...,T_v.$ Choose generators $\{\l_i, \m_i \}$ for each $T_i$ as above. The oriented boundary curves of a spun-normal surface $S$ with normal $Q$--coordinate $N = N(S)$ are then determined by the vector:
\begin{equation}\label{comb:normal curves}
    (-\nu_N(\l_1), \nu_N(\m_1),{\ldots} ,-\nu_N(\l_v), \nu_N(\m_v)) 
\in \Z^{2v},
\end{equation}
and boundary curves of linear combinations of compatible surfaces can be determined using linear combinations of vectors of the form (\ref{comb:normal curves}) since $\nu$ is additive. The intersection pairing has the following form in this case:
\begin{equation}\label{eq:pairing}
\nu(N) \star \nu(L) =  \frac{1}{2}\sum_{i=1}^v \nu_N(\m_i)\nu_L(\l_i)-\nu_N(\l_i)\nu_L(\m_i).
\end{equation}


\subsection{The figure eight knot complement}
\label{The figure eight knot}

\begin{figure}[t]
\psfrag{0}{{\small $0$}}
\psfrag{1}{{\small $1$}}
\psfrag{2}{{\small $2$}}
\psfrag{3}{{\small $3$}}
\psfrag{4}{{\small $4$}}
\psfrag{5}{{\small $5$}}
\psfrag{6}{{\small $6$}}
\psfrag{7}{{\small $7$}}
\psfrag{a}{{\small $z$}}
\psfrag{b}{{\small $z'$}}
\psfrag{c}{{\small $z''$}}
\psfrag{w}{{\small $w$}}
\psfrag{x}{{\small $w'$}}
\psfrag{y}{{\small $w''$}}
    \begin{center}
      \includegraphics[width=8cm]{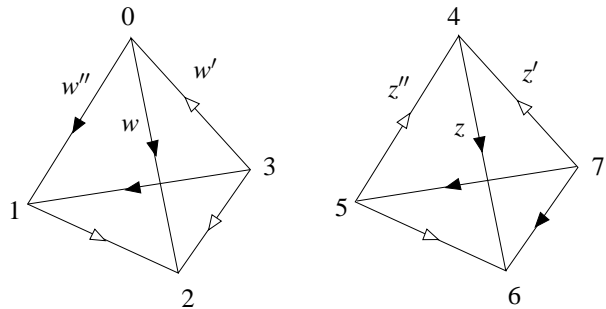}
    \end{center}
    \caption{An ideal triangulation of the figure eight knot complement}
    \label{fig:fig8_tets}
\end{figure}

Let $M$ denote the complement of the figure eight knot. An oriented, ideal triangulation of $M$ is encoded in Figure \ref{fig:fig8_tets}. Since $M$ is oriented, the convention given in Subsection \ref{comb:conventions} will be used, with the quadrilateral types dual to $w^{(k)}$ and $z^{(k)}$ denoted by $p^{(k)}$ and $q^{(k)}$ respectively. One computes:
\begin{equation*}
B = \begin{pmatrix} -1 & -1 & 2 & -1 & -1 & 2 \\
                          1 & 1 & -2 & 1 & 1 & -2 \end{pmatrix},
\end{equation*}
which determines a single $Q$--matching equation:
\begin{equation*}
0 = p + p' - 2 p''+q + q' - 2 q''.
\end{equation*}
This implies that the space $PQ(\tri )$ is four--dimensional. 

\begin{figure}[h]
\psfrag{0}{{\small $0$}}
\psfrag{1}{{\small $1$}}
\psfrag{2}{{\small $2$}}
\psfrag{3}{{\small $3$}}
\psfrag{4}{{\small $4$}}
\psfrag{5}{{\small $5$}}
\psfrag{6}{{\small $6$}}
\psfrag{7}{{\small $7$}}
\psfrag{a}{{\small $z$}}
\psfrag{b}{{\small $z'$}}
\psfrag{c}{{\small $z''$}}
\psfrag{d}{{\small $w$}}
\psfrag{e}{{\small $w'$}}
\psfrag{f}{{\small $w''$}}
\begin{center}
  \includegraphics[width=8cm]{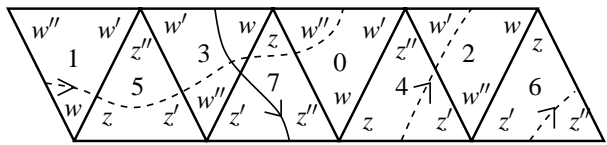}
\end{center}
    \caption{The induced triangulation of the vertex linking torus, where the sides of the rectangle are identified by translations parallel to its sides and triangle $i$ is dual to vertex $i$ in Figure~\ref{fig:fig8_tets}. The shown elementary curves are the standard meridian (solid) and longitude (dashed).}
    \label{fig:torus triang fig8}
\end{figure}

The induced triangulation of the vertex linking surface is shown in Figure \ref{fig:torus triang fig8} and used to determine the linear functionals associated to the standard peripheral curves:
\begin{align*}
\nu (\l ) &=  2 p +2 p' - 4 p'',\\
\nu (\m ) &= -p'+p''-q+q''.
\end{align*}
The normal $Q$--coordinates of closed normal surfaces satisfy $\nu (\l ) = \nu (\m )=0,$ whence any closed (embedded) normal surface is vertex linking. It follows that $\N (\tri )$ is zero--dimensional. A direct calculation reveals that there are four projective classes of admissible solutions; all have minimal representative a once--punctured Klein bottle. Their normal $Q$--coordinates and boundary slopes are listed in Table \ref{tab:surfaces fig8}. This calculation in particular shows that no spun-normal surface is a Seifert surface for the knot. 

\begin{table}[h]
\begin{center}
\begin{tabular}{ c | r | r | r}
solution & $\nu (\m )$  & $\nu (\l )$ & slope \\
\hline
(2,0,0,0,0,1) & 1 &4  &        --4  \\
(0,2,0,0,0,1) & --1 & 4  &       4   \\
(0,0,1,2,0,0) & --1 & --4  & --4 \\
(0,0,1,0,2,0) &  1&  --4  &     4 \\
\end{tabular}
\end{center}
\caption{Normal surface in the figure eight knot complement}
\label{tab:surfaces fig8}
\end{table}


\subsection{The Gieseking manifold}
\label{subsec:Gieseking}

The Gieseking manifold, $M',$ is double covered by the figure eight knot complement, $M,$ and the covering transformation is encoded by the involution $(05)(14)(26)(37)$ on the vertices in Figure \ref{fig:fig8_tets}. The resulting ideal triangulation $\tri'$ of $M'$ has one ideal 3--singlex and one ideal 1--singlex. Denote the three quadrilateral types in $M'$ by $r,r',r'',$ where $r^{(i)}$ lifts to $p^{(i)}.$ The $Q$--matching equation is $r+r'-2r''=0.$ It can be worked out from the triangulation or by observing that the induced involution on quadrilateral types in $M$ is $(p \ q')(p'\ q)(p''\ q'').$ Thus, $\dim PQ(\tri')=1$ and $\N(\tri')=\emptyset.$

The boundary curve map is defined via the induced triangulation of the double cover of the Klein bottle; using the generators from the above section, one has: $\nu (\l ) = - 2 r - 2 r' + 4 r'' = 0$ and $\nu (\m ) = -2r'+2r''.$ Generators $\l', \m'$ can be chosen for $H_1(B_\v)\cong \mathbb{Z}\oplus\mathbb{Z}_2$ such that the map $H_1(\widetilde{B}_\v) \to H_1(B_\v)$ is given by $\l\to \l'$ and $\m\to(\m')^2=0.$ The composition $$Q(\tri') \to \mathbb{Z}^2 \to H_1(\widetilde{B}_\v; \RR)\to H_1(B_\v; \RR)$$ is then
$$N \to (-\nu_N(\l), \nu_N(\m))=(0, \nu_N(\m)) \to \nu_N(\m) \l \to  \nu_N(\m) \l'.$$
Since $\nu (\m ) = -2r'+2r'',$ it follows that the map $\partial \co Q(\tri') \to H_1(B_\v; \RR)$ is surjective. Its restriction to integral points in $Q(\tri')$ has image of index two in $H_1(B_\v; \mathbb{Z}),$ which gives a subgroup of index four in $H_1(B_\v).$



\address{Department of Mathematics and Statistics, The University of Melbourne, VIC 3010, Australia} 
\email{tillmann@ms.unimelb.edu.au} 
\Addresses
\end{document}